\documentclass[opre,sglanonrev]{informs4}
\RequirePackage{tgtermes}
\RequirePackage{newtxtext}
\RequirePackage{newtxmath}
\RequirePackage{bm}
\RequirePackage{endnotes}

\OneAndAHalfSpacedXII

\usepackage{algorithm}
\usepackage{algpseudocode}
\usepackage{tikz}
\usepackage{comment}
\usepackage{subcaption}
\usepackage{xcolor}
\usepackage{booktabs}
\usepackage{hyperref}
\usepackage{natbib}
 \bibpunct[, ]{(}{)}{,}{a}{}{,}%
 \def\bibfont{\small}%

\newcommand{\BN}{\mathbb{N}}

\newcommand{\BR}{\mathbb{R}}

\newcommand{\CR}{\mathcal{R}}
\newcommand{\CC}{\mathcal{C}}

\newcommand{\addrefnum}{\stepcounter{equation}\tag{\theequation}}

\EquationsNumberedThrough    

\TheoremsNumberedBySection
\ECRepeatTheorems  %

\begin{document}

\RUNAUTHOR{Di, Andrad\'{o}ttir, Ayhan}

\RUNTITLE{{Pricing in Queues with Abandonments}}

\TITLE{{Pricing in Queues with Abandonments:\\ Optimal Policies and Practical Heuristics}}

\ARTICLEAUTHORS{

\AUTHOR{Jieqi Di, Sigr\'{u}n Andrad\'{o}ttir, Hayriye Ayhan}
\AFF{H.\ Milton Stewart School of Industrial and Systems Engineering, Georgia Institute of Technology, Atlanta, GA 30332-0205, U.S.A.}

} 

\ABSTRACT{%
We investigate the optimal pricing strategy in a service-providing framework, where customers can leave the system prior to service completion.  In this setting, a price is quoted to an incoming customer based on the current number of customers in the system.  When the quoted price is lower than the price the incoming customer is willing to pay (which follows a fixed probability distribution), then the customer joins the system and a reward equal to the quoted price is earned.  A cost is incurred upon abandonment and a holding cost is incurred for customers waiting to be served.  Our goal is to determine the pricing policy that maximizes the long-run average profit. 
Unlike traditional queueing systems without abandonments, we show that the optimal quoted prices do not always increase with the queue length in this setting. We fully characterize the possible structure of the optimal dynamic pricing policy
and provide conditions guaranteeing that the optimal policy is increasing in the number of customers in the system.
Moreover, we introduce two heuristics that simplify the optimal dynamic pricing policy. Both heuristics admit customers until the number of customers in the system reaches a certain threshold.  The cutoff-static policy charges all admitted customers a fixed price while the two-price policy charges one price when the arriving customer can enter service immediately and another price if the customer needs to wait. By selecting the price(s) and threshold that maximize the long-run average profit, both heuristics achieve near optimality in general and the two-price policy provides more robustness compared to the cutoff-static policy. 
}%

\FUNDING{This research was supported by [2127778, NSF] and [2348409, NSF].}


\maketitle

\section{Introduction}

{
Pricing and revenue management in stochastic systems are challenging yet valuable areas of study. 
The unpredictability of customer behavior, particularly when customers may leave the system after joining but prior to completing service, introduces further complexity into these problems. Such abandonments may occur while waiting or during service, driven by impatience, external disruptions, or system failures. This phenomenon occurs across diverse industries. In call centers, customers who are placed on hold may become impatient and leave \citep{Gan.et.al.2003,Aksin.et.al.2017}. In healthcare, abandonments occur when patients leave without being seen due to long waits \citep{Green.et.al.2006,Batt&Terwiesch2015,Bolandifar.et.al.2019,RAVID.et.al.2025}, or in the case of organ transplant lists, when patients pass away (abandon) while waiting \citep{Su&Zenios2004}. Customer abandonments are also prevalent in grocery checkout lines \citep{Buell2021} and ride-hailing settings \citep{Yu.et.al.2022, Wang.et.al.2024}.} 

One type of abandonments in the existing literature is called rational abandonment, in which the customer reneging is strategic. Namely, {delay-time sensitive} customers are actively making the choice of leaving or staying while in the system. The papers focusing on this type of abandonment tend to study the abandonment strategies \citep{Hassin&Haviv2012, Afeche&Vahid2015,Ata&Peng2018,Abouee2022}. 

{However, customers may also abandon due to external disruptions or system failures. To capture these passive abandonments, which are also known as forced abandonments in \citet{BENELLI2021426} or abandonments due to exogenous reasons in Section 3.4 of \citet{Hassin&Snitkovsky2020}, we interpret abandonments as random events and assume that customers' abandonment times after joining the system follow a predetermined probability distribution (see also \citealt{Whitt&Ward2006, Bassamboo&Randhawa2016, Dong&Ibrahim2021}). The unpredictable nature of random abandonment leads to ambiguous congestion information, even when the queue length is observable. Given that congestion information is often obscure, customers may primarily rely on the quoted price in deciding whether or not to join the system, implying that customers are purely price-sensitive.}

{The prevalence of abandonments in real-life service systems motivates us to study the pricing strategies in these systems in order to improve profits.}
We consider a queueing system whose arrivals follow a Poisson process and the service times are exponentially distributed. A price is quoted to any arriving customer based on the current number of customers in the system. Viewing that quoted price, the customer makes the decision to join the system or balk. 
The customer can abandon while waiting or during service, with associated costs to the system, and the abandonment time is exponentially distributed. There is also a per unit time, per customer holding cost. 
To maximize the long-run average profit for this model, we employ a Markov decision process to determine the optimal quoted prices. This strategy, known as dynamic pricing, adjusts the price according to the number of customers in the system. 

{
A focal application of our pricing model is related to the growth of expedited service mechanisms in modern platforms \citep{Roet-Green&Shetty2022}.
In these systems, customers pay a per-use premium to access dedicated high-priority service operating with strict priority over base service. The customers joining the expedited service no longer consider the base service as an option. All premium customers receive identical fast treatment but may still abandon due to impatience, disruptions, or system failures. Since the expedited service begins immediately upon payment by placing customers in the premium queue, payments are typically non-refundable even if customers later abandon. Namely, the payment is for preferential treatment already enjoyed upon joining. In ride‑hailing and on‑demand delivery, options like Uber Priority \citeyearpar{uber_priority}  
and DoorDash Express Delivery \citeyearpar{doordash_priority}
place customers who pay an additional fee for expedited service into such a premium lane with more drivers or couriers (potentially with a faster service).
Customers may abandon before driver assignment (while waiting) or afterward (during service) without receiving refunds. 
Similarly, in cloud computing, per-request high-priority modes admit jobs to strict-priority GPU pools reserved for premium customers who pay the priority fee. Jobs may cancel or time out before GPU upload or during processing, with no refunds issued.
This expedited mechanism extends to appointment-based services for automotive or equipment repair (3M Repair \& Service \citeyear{3M_repair}, PEV Works \citeyear{pev_works}). Rather than scheduling distant appointments, service providers can offer expedited systems that assign customers to the next available server. Customers demanding speedy service who prefer not to wait for scheduled appointments can pay a premium to enter this priority system. Their abandonment remains subject to no refund, as the premium already secured their priority service status.
}

{To the best of our knowledge, this is the first paper that characterizes the optimal dynamic pricing policy exactly in queues with abandonments. Our major contributions are threefold.} 
\begin{enumerate}
\item \textit{Our paper reveals that customer abandonment fundamentally transforms the structure of optimal pricing policies}. From our numerical experiments, we observe that in queueing systems with abandonments, the optimal quoted prices do not always increase as the number of customers in the system grows. This behavior contrasts sharply with queueing systems without abandonments, where the optimal quoted price is known to increase with system congestion \citep{Low1974}. We prove the optimal policy always has a uni-modal structure (first decreasing and then increasing), and we also provide a condition guaranteeing that the optimal quoted price increases with the number of customers in the system. These results not only highlight a fundamental difference between systems with and without abandonments, {but also enable the development of a simplified algorithm for efficiently determining the optimal dynamic pricing policy.}
\item \textit{Our paper provides managerial guidance for determining when optimal prices should rise versus fall as system congestion increases, enabling companies to incorporate abandonment considerations into their pricing schemes.} In short, if providing immediate service is less costly than having customers abandon while waiting, then reducing congestion via increasing the prices is profitable, creating a monotone increasing pricing structure. However, when queue abandonment is less costly than providing immediate service, counterintuitively, it can be profitable to increase congestion by reducing the prices, which leads to the non-monotone pricing. {We also demonstrate that our insights match real-life scenarios in expedited service.} 
\item \textit{Our paper introduces practical pricing heuristics that overcome the implementation barriers of traditional dynamic pricing methods}. Due to the curse of dimensionality, finding the optimal dynamic pricing policy can be cumbersome. Also, dynamic pricing can be complex and difficult to implement because the price discrimination may cause dissatisfaction and cast doubts on social fairness. Therefore, we propose two heuristics, the cutoff-static policy and the two-price policy, that combine static pricing and admission control, offering a more practical and robust implementation compared to optimal dynamic pricing. {We provide theoretical upper bounds on the optimality gap between the optimal dynamic policy and the heuristics. We also provide numerical experiments that suggest both heuristics are near optimal in general.}
\end{enumerate}

The rest of this paper is organized as follows. Section 2 reviews the related literature. Section 3 introduces the formulation of the underlying queueing system and the setup of the pricing problem. {Section 4 provides the structure of the optimal dynamic pricing policy and develops an efficient policy iteration algorithm to find the optimal policy. Section 5 establishes the best cutoff-static and two-price policies, provides methods to find both heuristics efficiently, and characterizes theoretical upper bounds on the gap between the optimal dynamic policy and the heuristics.}
Section 6 provides numerical results from randomly generated datasets {that demonstrate the near optimal performance of the heuristics}. Conclusions are drawn in Section 7. {The detailed proofs for Sections 4 and 5 and additional tables and figures are provided in the Electronic Companion.}

\section{Literature Review} \label{sec.literature.review}

{In this section, we categorize the related literature into two parts. First, we review the literature on queueing systems with abandonments, including admission control in queues with abandonments in Section \ref{subsec.lit.review.queue.abandon}. Next, we examine works about pricing problems in queueing systems in Section \ref{subsec.lit.review.pricing}. Although each of these topics has attracted considerable attention, limited research has been conducted on pricing in queueing systems with abandonments.}

\subsection{Queueing System with Abandonments} \label{subsec.lit.review.queue.abandon}
{The abandonment mechanism in our work is closely related to the Erlang-A model, first introduced in \citet{Palm1957}, and has been of interest for a long time \citep{Ancker1963,Baccelli1984}.} More recently, \citet{Jouini_Dallery_2007}, \citet{Boxma2011}, \citet{Brandt2013}, and \citet{Legros2018} have done exact analysis on performance measures for queueing systems with abandonments under different service time and abandonment time distributions.
Due to the interest in large-scale cloud computing centers and call centers, {many authors considered the asymptotic analysis of queues with abandonments under heavy traffic assumptions \citep{Garnett.et.al.2002,Whitt.2006, Liu&Whitt2011,Dai2013}.} \citet{WARD2012} identified the approximate asymptotic behavior of performance measures in the conventional heavy-traffic regime when the inter-arrival times, service times, and abandonment times share general distributions. {\citet{Dong.et.al.2015} carried out the performance analysis in the Quality-and-Efficiency Driven (QED) regime for a modified Erlang-A model.}

These exact and asymptotic analyses built the ground for studying the control of queues with abandonments, {including staffing in \citet{Rouba2018}, scheduling in \citet{Atar.et.al.2010,Dong&Ibrahim2024} and information announcement in \citet{afeche.et.al.2024}. Here, we emphasize the admission control problem, which shares similarity with pricing. Any admission control policy can be viewed as a pricing policy with two-prices, in which one of the prices is low enough for all customers to join the system while the other price is expensive enough to prevent any incoming customers from entering the system. \citet{Zayas-caban2020} considered a finite horizon admission control problem in queues with abandonments, two classes of customers, and periodically varying parameters. \citet{AYHAN2022} tackled a similar problem for single-class queues in an infinite time horizon setting and demonstrated that a threshold-type policy is optimal. } 
\citet{Runhua.2024} extended the infinite time horizon problem into a multi-class customer setting and completely characterized the structure of the optimal policy. The three papers mentioned above involve Markovian queueing systems. For queueing systems with general service and abandonment times, \citet{Ward&Kumar2008} provided an asymptotically optimal admission control policy. 
All the works reviewed so far involve queues with abandonments but none has discussed the pricing problem on these systems.

\subsection{Pricing in Queueing Systems} \label{subsec.lit.review.pricing}

{We break down the literature related to pricing in queues based on how the state of the queue impacts customers' choices to join the system. If the state of the queue does not directly impact customer joining decision, the customers are purely price-sensitive. In contrast, the customer can be delay-sensitive, where the state of the queue, as a signal of customers' waiting time, may impact customers' joining decision. Customers may observe the state of the queue to predict their waiting time ({observable queues}), or not observe the state of the queue but have an expectation based on past experience of what their waiting time will be (unobservable queues). }

{In our model, due to the variability in the abandonment process, the current queue length may not be a reliable or precise predictor of the individual customer's final waiting time. This ambiguity leads the customer to primarily rely on the price signal. Therefore, we adapt a pure price-sensitive model.} \citet{Low1974} presented a dynamic pricing model in which the quoted prices change when the number of customers in the system changes and proved that the optimal prices are monotone in the queue length. Using Low's pricing model, \citet{ZIYA2006214} studied static pricing for queueing systems with general arrival and service distributions. {\citet{Maoui_Ayhan_Foley_2007,Maoui.et.al.2009} identified the static and dynamic pricing that maximize the long-run average reward for queueing systems with holding costs, and \citet{AKTARANKALAYCI2009120} studied the sensitivity of pricing policies in a similar queueing model.} \citet{Wang.et.al.2019} considered dynamic pricing on tandem queues with finite buffers using a Markov decision process and concluded that the optimal policy has a three-way monotone structure. 
{\citet{Ata&Shneorson2006} extended Low's model by studying the dynamic pricing and capacity sizing problem and showed that the optimal arrival and service rates have monotone structures. \citet{LEE2014527} analyzed pricing and capacity sizing problems asymptotically in $GI/GI/1$ queue under heavy-traffic assumptions.}
A similar dynamic pricing model is also introduced in inventory control \citep{Gallego&Ryzin1994,Federgruen&Heching,Chen&Chao2020,YAN2024935}.

{For delay-sensitive customers, \citet{Naor1969} introduced observable queues and analyzed static pricing in this setting.   \citet{Chen&Frank2001} generalized Naor's results to dynamic pricing.} Using Markov decision processes, \citet{Yoon&Lewis2004} and \citet{Cil2011} studied dynamic pricing for single station observable queues with multiple classes of customers. \citet{Maglaras&Zeevi2003}, \citet{Maglaras2006}, and \citet{Kim&Randhawa2017} analyzed pricing problems asymptotically in observable queues under heavy-traffic assumptions using fluid approximations or diffusion models. {\citet{afeche&Baris2013} considered the pricing problem in observable queues with Bayesian learning. \citet{Edlson&Hildebrand1975} is the first paper that considered unobservable queues, where they studied pricing problems for revenue and social welfare maximization. \citet{Afeche2013} and \citet{Afeche2016} studied a joint static pricing, lead-time quotation, and scheduling problem for an unobservable queue with multiple types of delay costs. \citet{Haviv&Randhawa2014} characterized the optimal static price for revenue maximization in an unobservable queue with unknown demand.}

The existing literature also includes comparisons between different pricing strategies. 
\citet{Besbes.et.al.2021} compared the long-run average gain of optimal static and dynamic pricing in queueing systems with reusable resources and showed that static pricing guarantees a certain percentage of the profit achieved by optimal dynamic pricing. \citet{Bergquist&Elmachtoub2023} introduced a pricing strategy that resembles static pricing and showed that this strategy is near optimal in an $M/M/1$ queue and {\citet{adam&jiaqi2025} considered the related problem in an Erlang Loss model.} \citet{Kim&Randhawa2017} proposed a two-price policy for $M/M/1$ queues and showed that under a heavy traffic assumption, the gap between the two-price policy and the optimal dynamic policy is within $O(N^{1/3})$, where $N$ is the system capacity. \citet{Varma&etal2023} adopted a two-price policy for a two-sided matching queue and proved an optimality gap result similar to that of \citet{Kim&Randhawa2017}. {\citet{wenxin2025} considered a stock-dependent two-price policy for reusable resources and showed that the performance loss can be as low as $\mathcal{O}((\log c)^2)$. Examples of other important works which analyze pricing in queueing models that are tailored for more specific applications are \citet{afeche&baron2019,Lin.et.al.2023}, \citet{Ming.hu.2026} and \citet{Wang.et.al.2026}. }
Although the papers above indicate that pricing in queueing systems has been studied extensively, the literature on pricing in queues with abandonments is scarce.

To the best of our knowledge, \citet{Lee&Ward2019} is the only paper addressing pricing problems in queues with random abandonments. They studied pricing and capacity sizing problems in an $M/GI/1+GI$ queue and investigated the asymptotic behavior of optimal static pricing and the loss in profit when customer abandonments are ignored.
In comparison to their work, we introduce a Markov decision process to solve for the optimal dynamic pricing exactly without the heavy-traffic assumption {for multi-server systems}. We characterize the structure of the optimal dynamic pricing policy, {which is not always monotone,} and propose an algorithm to find it efficiently under mild assumptions guaranteeing that the optimal dynamic pricing policy is monotone. We also develop two pricing heuristics that are easier to implement but remain near optimal compared to optimal dynamic pricing.

\section{Model Description} \label{sec.model_description}

We consider a single-station Markovian queueing system. We assume that the customers arrive with respect to a Poisson process of rate $\Lambda>0$ and customers joining the system get served by any one of the $m \geq 1$ identical servers whose service times are independent exponential random variables with rate $\mu>0$. The queueing system can hold at most $N$ customers, where $m\leq N<\infty$, and arriving customers finding $N$ customers in the system have to leave. 
The customers who are waiting in the queue can leave after an exponential amount of time with rate $\theta_q> 0$. The customers who are being served can also abandon after an exponential amount of time with rate $\theta_s\geq 0$.

When an incoming customer arrives to the system, a price is quoted to that customer based on the current number of customers in the system. Each incoming customer has an evaluation, drawn independently for different customers from a cumulative distribution function $F$, with $F(p)$ denoting their probability of not joining the system for any given price $p$. We assume that the evaluation distribution $F(\cdot)$ is continuous and has a positive probability density function $f(\cdot)$ on the range of acceptable prices $[a,b]$, where $0\leq a<b< \infty$. Let $\Bar{F}(p) = 1 - F(p)$ indicate the probability of the customer joining the system when we quote a price $p$. Let $p_n\in[a,b]$ be the quoted price when there are $n$ customers in the system. Then, the arrival rate in state $n$ can be expressed as
$ \lambda_n =\Lambda\times \Bar{F}(p_n).$
If the customer joins the system, a reward equal to the quoted price is earned.
Furthermore, there is a per unit time, per customer holding cost of $c_h\geq 0$ for customers in the system and a cost of $c_q\geq 0$ ($c_s\geq 0$) if a customer abandons while waiting in the queue (during service). Our objective is to determine the optimal dynamic policy for the quoted price that maximizes the long-run average profit.

Since the cumulative distribution function $F(p)$ is strictly increasing and continuous on the range $(a,b)$,  
we can define $\bar{F}^{-1}$ for any value $y\in[0,1]$ as
$
\bar{F}^{-1}(y) = \inf\{ x\in\BR : \bar{F}(x)\leq y \}.
$
{Based on this definition of $\bar{F}^{-1}(\cdot)$, the quoted price $p_n$ corresponding to an arrival rate $\lambda_n$ at state $n$ is given as $p_n=\Bar{F}^{-1}({\lambda_n}/{\Lambda})$, which means that the service provider should quote the price $p_n$ to induce an arrival rate of $\lambda_n$ at state $n$. Given this relation, finding the optimal pricing policy is equivalent to finding the optimal arrival rate policy. Thus, we treat the arrival rates as the primary decision variables, which follows the convention in the revenue management literature \citep{Gallego&Ryzin1994}. Since $\bar{F}(\cdot)$ is strictly decreasing, then the quantile function $\bar{F}^{-1}(\cdot)$ is non-increasing, which implies that a decrease in the optimal arrival rate corresponds to an increase in the optimal quoted price, and vice versa.}

Recall that $\Lambda$ is the maximum possible arrival rate and $\lambda_n\in[0,\Lambda]$ is the arrival rate at state $n$. Since there is no decision to be made when $n=N$, $\lambda_N=0$. Let $A_n$ be the set of allowable actions in state $n$, where $A_n = [0,\Lambda]$ for $n=0,1,\ldots,N-1$ and $A_N=\{0\}$. Let $\Vec{\mathbf{\lambda}}=(\lambda_0,\ldots,\lambda_N)\in [0,\Lambda]^{N}\times\{0\}$ be an arrival rate policy and $\Pi$ be the set of all possible arrival rate policies. Denote $\{X^{\Vec{\lambda}}(t):t\geq 0\}$ as the stochastic process representing the number of customers in the system at time $t$ under policy $\Vec{\lambda}\in\Pi$. Then, $\{X^{\Vec{\lambda}}(t):t\geq 0\}$ is a continuous-time Markov chain with state space $\mathcal{S}=\{0,1,\ldots,N\}$.

Since the buffer size is finite, all the transition rates in the system are no larger than $m (\mu+\theta_s) + (N-m) \theta_q$. Then, the continuous-time Markov chain $\{X^{\Vec{\lambda}}(t)\}$ is uniformizable. Without loss of generality, let uniformization constant be 
$
1 = \Lambda + m (\mu+\theta_s) + (N-m) \theta_q
$
by scaling the time unit. Notice that changing the time unit will affect the per unit time, per customer holding cost. With a slight abuse of notation, we redefine $c_h$ as the holding cost after the change of time unit. 
 
From uniformization, we reduce $\{ X^{\vec{\lambda}}(t)\}$ to a discrete-time Markov decision problem.
Let $r_n(\lambda_n)$ be the immediate reward when action $\lambda_n$ is chosen in state $n$. We denote the reward vector $\mathbf{r}^{\vec{\lambda}}$ under an arrival rate policy $\vec{\lambda}$ as $\mathbf{r}^{\vec{\lambda}} = (r_0(\lambda_0),\ldots,r_N(\lambda_{N}))$. Let the one-step probability from state $n$ to state $j$ when action $\lambda_n$ is chosen in state $n$ be  $\mathbf{P}(j|n, \lambda_n)$. Denote $\mathbf{P}^{\vec{\lambda}} = [\mathbf{P}(j|n, \lambda_n)]_{n,j}$ as the transition probability matrix.
We also define a short-hand notation for the death rate in state $n$ as 
\[
\gamma_n = \min\{n,m\}(\mu+\theta_s)+\max\{n-m,0\}\theta_q \text{ for } n=0,\ldots,N.
\]
{\begingroup\renewcommand{\baselinestretch}{0.8}\normalsize
Then, we have the transition probabilities as
$$
\mathbf{P}(j|n, \lambda_n) =
\begin{cases}
\lambda_n & \text{ for } n=0,\ldots,N-1 \text{ and } j=n+1,\\
\gamma_n & \text{ for } n=1,\ldots,N \text{ and } j=n-1,\\
1-\lambda_n-\gamma_n & \text{ for } n=0,\ldots,N \text{ and } j=n,\\
0 &\text{   otherwise}.
\end{cases}
$$
\endgroup}
We then provide the reward functions of this Markov decision process as
$$r_n(\lambda_n) = \lambda_n p_n - nc_h - \min\{n,m\}c_s\theta_s -\max\{n-m,0\}c_q\theta_q, $$
where $p_n= \Bar{F}^{-1}({\lambda_n}/{\Lambda})$ for $n=0,\ldots,N-1$, and 
$
r_N(\lambda_N) = r_N(0) = - Nc_h - m c_s\theta_s -(N-m)c_q\theta_q.
$

Now for any given policy $\Vec{\lambda}\in\Pi$, we provide the steady-state distribution and the long-run average gain under that policy. Let $P_n(\Vec{\lambda})$ be the steady-state probabilities of being in state $n$ for $n=0,\ldots, N$. We have 
$
P_n(\Vec{\lambda})=P_0(\vec{\lambda})\times a_n(\Vec{\lambda}),
$
where $a_n$ and $P_0$ are defined as follows: for any $\Vec{\lambda}\in\Pi$, $a_0(\Vec{\lambda})=1$,
\[
a_n(\Vec{\lambda}) = \dfrac{\lambda_{n-1}}{\gamma_n}\times a_{n-1}(\Vec{\lambda}) \text{ for } 1\leq n\leq N, \text{ and }
P_0(\Vec{\lambda}) = \dfrac{1}{\sum_{n=0}^N a_n(\Vec{\lambda})}.
\]
Furthermore, the long-run average gain $g^{\vec{\lambda}}$ under a given policy $\Vec{\lambda}$ is $g^{\Vec{\lambda}} = \sum_{n=0}^N r_n(\lambda_n) \times P_n(\vec{\lambda}).$

Note that due to non-idling servers and abandonment, the Markov decision process is unichain. Hence, the Bellman equations for all $n=0,\ldots,N$ are given as
\[
g^* + h(n) = \max_{\lambda\in A_n}\Bigg\{ r_n(\lambda)+\sum_{j=0}^N \mathbf{P}(j|n,\lambda) h(j)  \Bigg\}, \addrefnum \label{Bellman.equation}
\]
where $g^*\in\BR$ is the optimal long-run average gain and $h(n)\in\BR$ is the bias for state $n$ (see Sections 8.2.1 and 8.4.1 in \citet{Puterman.1994} for details). Furthermore, we know from Theorem 8.4.7 in \citet{Puterman.1994} that there exists an optimal Markovian stationary deterministic policy for the arrival rate at each state $n$, which we define as
\[
\lambda_n^* = \inf \Bigg( \argmax_{\lambda\in A_n}\Big\{ r_n(\lambda)+\sum_{j=0}^N \mathbf{P}(j|n,\lambda) h(j)  \Big\} \Bigg). \stepcounter{equation}\tag{\theequation} \label{def.opt.rate}
\]
Recall that when $n=N$, we have $A_N=\{0\}$, which means $\lambda_N^*=0$. In the later sections, we refer to this optimal dynamic policy $\vec{\lambda}^*$ as the optimal arrival rate policy.

\section{Structure of the Optimal Dynamic Policy} \label{sec.monotonicity}

\noindent In this section, we examine the structural properties of optimal dynamic policies. Understanding the structure aids in developing efficient algorithms to determine the optimal policy.
One important structural property is monotonicity. In our pricing context, monotonicity implies that the optimal quoted price increases as the system becomes more congested. Correspondingly, when we consider the arrival rate as the decision variable, monotonicity would mean that the optimal arrival rate decreases as the number of customers in the system increases. To formalize this concept, we provide the following definition.

\begin{definition} \label{def.mono}
An arrival rate policy $\Vec{\lambda}\in [0,\Lambda]^N\times \{0\}$ has a monotone decreasing structure when $\lambda_n\geq \lambda_{n+1}$ for all $n=0,\ldots, N-1$. {A pricing policy has a monotone increasing structure if the corresponding arrival rate policy has a monotone decreasing structure.}
\end{definition}

{From a managerial perspective, when customers are only price sensitive, it is intuitive to increase the price and reduce the number of joining customers as the system gets more congested.}
A larger queue incurs higher holding costs and increases the risk of customer abandonments, potentially leading to significant congestion penalties. 
{Therefore, when the system is congested, raising the quoted price reduces the rate at which new customers join the system to ease the congestion. }
The decreasing structure of the arrival rate policy captures this intuition mathematically. 

However, numerical experiments suggest there are cases where the optimal arrival rate policy is not monotone decreasing. Figure 1 provides such an example where the chosen system parameters are $\Lambda=5, \mu=2, N=10,m=1,\theta_s=1,\theta_q=2, c_s=40, c_q=10,c_h=10,$ and $F(\cdot)\sim U(20,50)$. In Figure 1(a), the optimal arrival rate policy first increases and then decreases as the number of customers in the system increases. To see this tendency more clearly, Figure 1(b) omits the last state in which the optimal arrival rate is always 0.

\begin{figure}[H]
	\centering
	\caption{The plot of an optimal arrival rate policy for each state}
	\begin{minipage}{0.45\textwidth}
		\centering
		\includegraphics[width=\columnwidth]{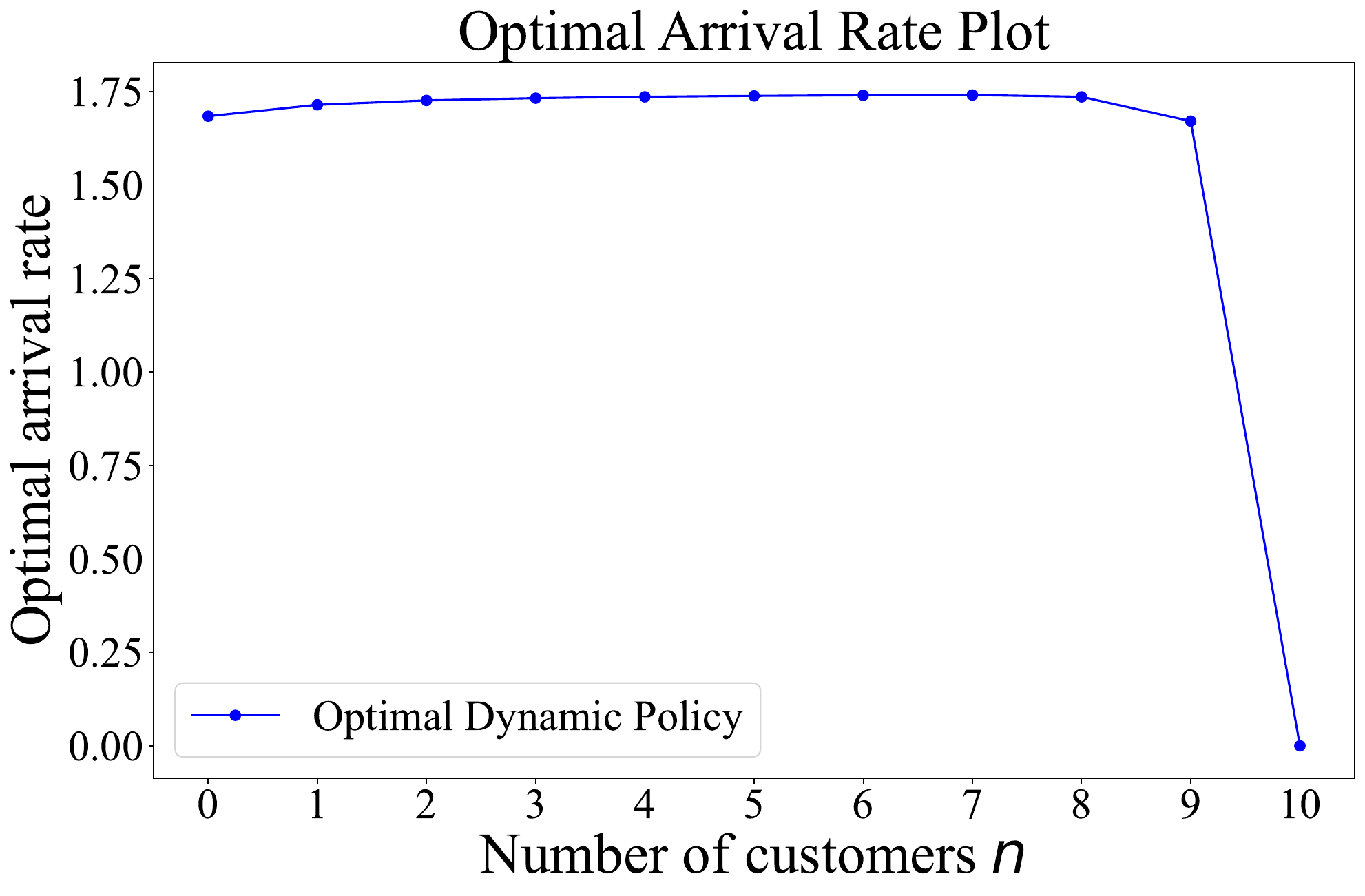}
		\subcaption{}
	\end{minipage}
	\begin{minipage}{0.45\textwidth}
		\centering
		\includegraphics[width=\columnwidth]{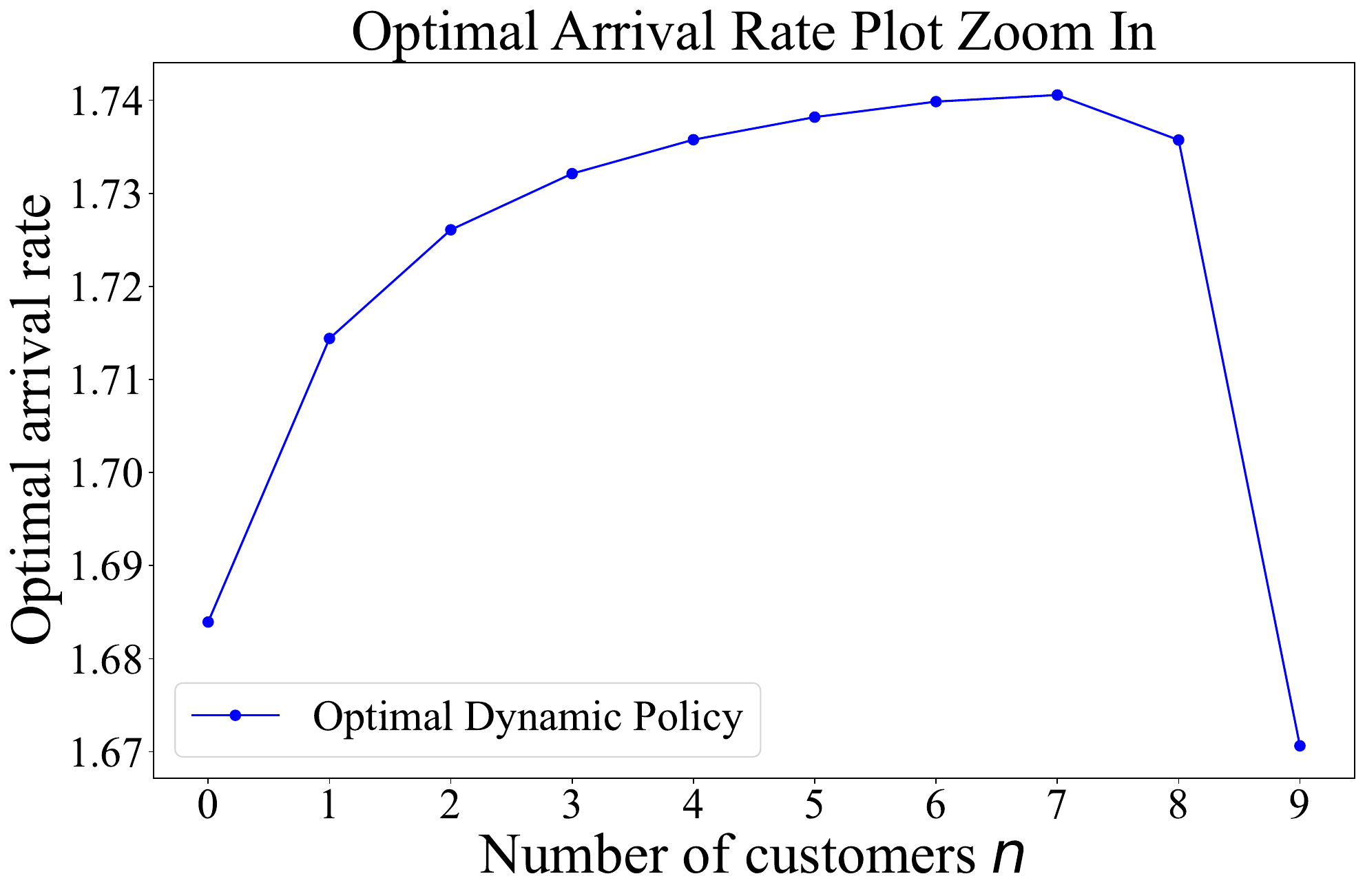}
		\subcaption{}
	\end{minipage}
	\label{fig:Condition_a_violation}
\end{figure}
\vspace{-0.4in}

We refer to this arrival rate structure shown in Figure 1 as a uni-modal structure. To formalize this concept, we provide the following definition.

\begin{definition} \label{def.uni.modal}
An arrival rate policy $\vec{\lambda}\in [0,\Lambda]^N\times \{0\}$ has a uni-modal structure when there exists an integer $\hat{n}\in \{0,\ldots,N-1 \}$ where $\lambda_n \leq \lambda_{n+1}$ for  $0\leq n < \hat{n}$ and $\lambda_{n}\geq \lambda_{n+1}$ for $\hat{n} \leq n < N$. 
\end{definition}

Notice that when $\hat{n}=0$, we have $\lambda_{n}\geq \lambda_{n+1}$ for any $n=0,\ldots,N-1$. Namely, monotone decreasing structure is a special case of uni-modal structure with $\hat{n}=0$. {If the optimal arrival rate has a uni-modal structure with $\hat{n}>0$, the optimal pricing policy is not monotone and the optimal quoted price first decreases and then increases as the number of customers in the system increases. Otherwise, the optimal arrival rate policy is monotone decreasing, inducing a monotone increasing quoted price.} {In Section \ref{subsec.monotonicity.main.results}, we show that the optimal arrival rate always has a uni-modal structure and provide a sufficient condition for $\hat{n}=0$. Section \ref{subsec:necessary_condition} investigates whether the sufficient condition in Section \ref{subsec.monotonicity.main.results} for monotonicity is also necessary and Section \ref{subsec:policy_itration} provides a policy iteration algorithm under uni-modal structure.}

\subsection{Major Results} \label{subsec.monotonicity.main.results}

In this section, we show that the optimal arrival rate policy is always uni-modal and is monotone decreasing ($\hat{n}=0$) under a certain condition. {Given the relation between prices and arrival rates, these results show that in general, the optimal quoted price may decrease and then increase when the number of customers in the system increases. However, under certain conditions, optimal quoted price is always increasing. The proofs for Theorems \ref{Theorem.uni_modal} and \ref{Theorem.monotonic} are in Section \ref{ec.proof.section.4}.} 

\begin{theorem} \label{Theorem.uni_modal}
The optimal arrival rate policy has a uni-modal structure.
\end{theorem}

\begin{theorem} \label{Theorem.monotonic}
{Let $\CC_s=\frac{c_h+c_s\theta_s}{\mu+\theta_s}$ and $\CC_q= \frac{c_h}{\theta_q}+c_q$.} When $\CC_s\leq \CC_q$, the optimal arrival rate policy has a monotone decreasing structure.
\end{theorem}

{Now we discuss the intuition behind these theorems and try to rationalize why the optimal policy is monotone when $\frac{c_h+c_s\theta_s}{\mu+\theta_s}\leq \frac{c_h}{\theta_q}+c_q$.} Recall that $c_h$ is the holding cost per customer per unit time and $\mu,\,\theta_s,\,\theta_q$ represent the service rate, service abandonment rate, and waiting abandonment rate. A customer arriving at a system with idle servers enters service immediately. Given exponentially distributed service and abandonment times, the expected duration until the customer leaves the system (either due to service completion or abandonment) is $\frac{1}{\mu+\theta_s}$, resulting in an expected holding cost of  $\frac{c_h}{\mu+\theta_s}$. The probability of service abandonment is $\frac{\theta_s}{\mu+\theta_s}$, resulting in an expected abandonment cost of $\frac{c_s\theta_s}{\mu+\theta_s}$. For a system with idle servers, the total expected cost per arriving customer is $\frac{c_h+c_s\theta_s}{\mu+\theta_s}$. {Thus, we can interpret $\mathcal{C}_s \coloneq \frac{c_h+c_s\theta_s}{\mu+\theta_s}$ as a measure of \textit{the expected total cost after joining service}.}

When a customer arrives at a system without idle servers, it enters into a queue to wait. {Before the customer enters service, the total time for that customer to be in the queue has an upper bound $1/\theta_q$, which implies that the expected holding cost before the customer enters service has an upper bound $c_h/\theta_q$. If the customer abandons before entering service, an additional abandonment cost $c_q$ is incurred. Therefore, we can interpret $\mathcal{C}_q\coloneq \frac{c_h}{\theta_q}+c_q$ as an upper bound of \textit{the expected total cost of a customer before joining service}.}

{When the expected total cost after joining service is lower than the upper bound of that before joining service ($\mathcal{C}_s=\frac{c_h+c_s\theta_s}{\mu+\theta_s}\leq  \frac{c_h}{\theta_q}+c_q=\mathcal{C}_q$), 
we can interpret that
the expected cost for an arriving customer who enters service immediately is lower than the potential cost for a customer who abandons prior to entering service. In this regime,} admitting customers who can be serviced immediately is desirable compated to admitting those who must wait in the queue and eventually abandon. {In other words, the system prefers the customer to eventually enter service because it can be costly for the customer to abandon without entering service.} Given this rationale, it is natural to argue that a strategy of {increasing the quoted price to reduce} the number of arriving customers (decreasing the arrival rate) as the system becomes more congested should enhance profitability.

{On the other hand, if $\mathcal{C}_s > \mathcal{C}_q$, we can interpret that customers who are not in service cost less compared to those who are in service. Then, it could be profitable for the system if customers abandon without service. Allowing more customers to join (increasing the arrival rate) as the system becomes more congested, potentially inducing more abandonments, resulting in a higher profit.}

{Expedited service in the ride-hailing setting will fall under the $\CC_s>\CC_q$ regime. In this priority system, a ``server" is a committed driver, and ``service" begins the moment a driver is assigned. Customers enter the priority system upon paying their premium and wait in a virtual pool for driver assignment. If customers abandon before driver assignment (from the queue), the system incurs minimal expected total cost before joining the service (at most $\CC_q$) because no driver's time or fuel has been committed to the request. However, if customers abandon after driver assignment (after entering service), the system incurs a high expected cost $\CC_s$ per customer. This includes the driver's time and fuel, which could have been used to profitably serve another customer. Consequently, when launching priority service, the platform should initially lower its premium to quickly assign all available drivers, then raise prices to benefit from potential abandonments.}

Additionally, when $c_h>0$ and both $\theta_s$ and $\theta_q$ are close to 0, {the expected total cost after joining service $\mathcal{C}_s$ approaches $\frac{c_h}{\mu}$ while $\mathcal{C}_q$ approaches infinity. This shows that when the abandonments from the system are negligible, then the system always has $\CC_s\leq \CC_q$}, indicating the optimal arrival rate policy is always monotone decreasing from Theorem \ref{Theorem.monotonic}. This recovers the results for queueing system without abandonments where the optimal arrival rate policy is always monotone.

\subsection{Is {$\CC_s\leq \CC_q$} a Necessary Condition for Monotonicity?} \label{subsec:necessary_condition}

\begin{figure}[H]
	\centering
	\caption{Scatter Plots for the Conditions of Monotonicity}
	\begin{minipage}{0.45\textwidth}
		\centering
		\includegraphics[width=\columnwidth]{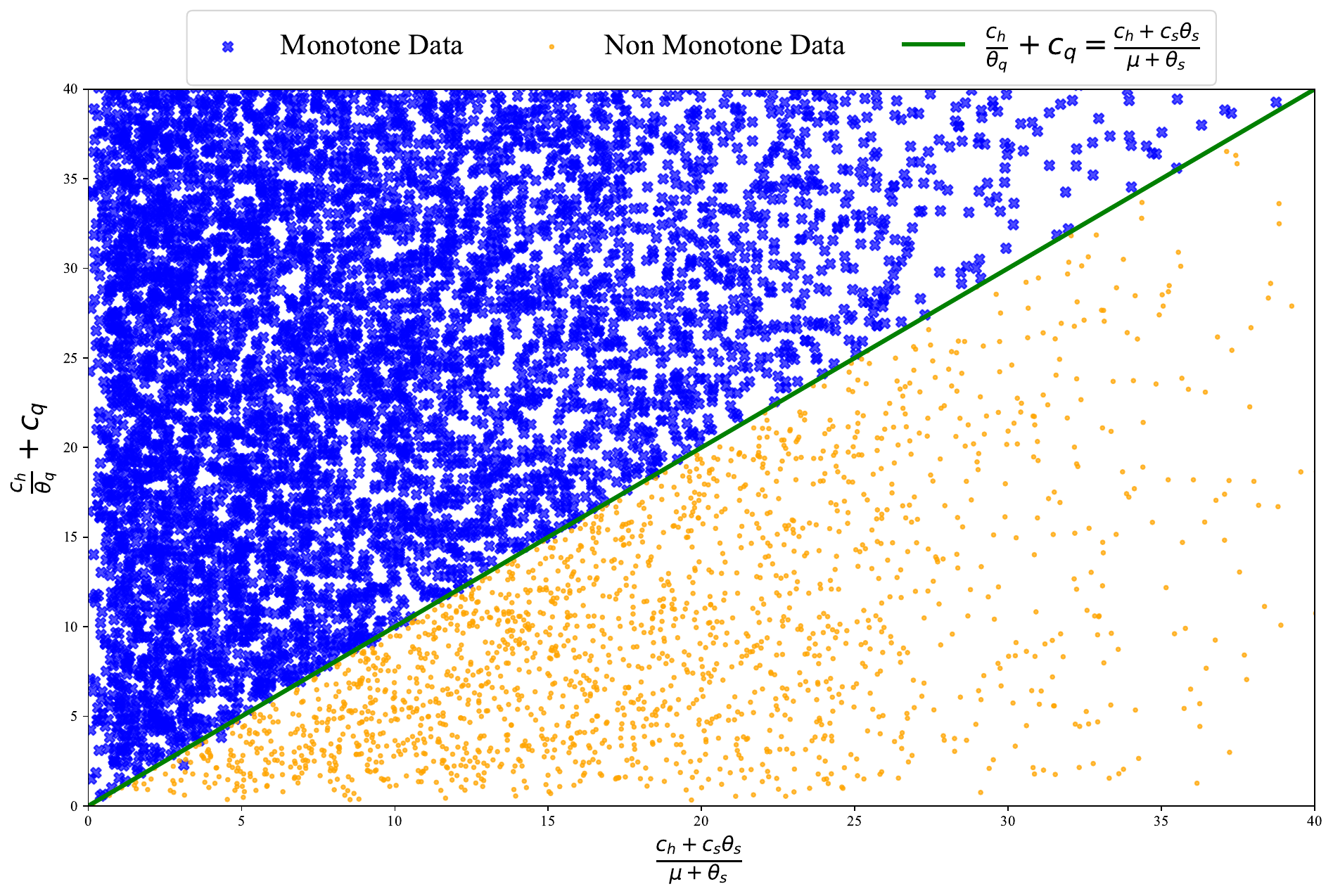}
		\subcaption{Evaluation distribution as $U(20,50)$}
        \label{Monotone_condition_unif}
	\end{minipage}
	\begin{minipage}{0.45\textwidth}
		\centering
		\includegraphics[width=\columnwidth]{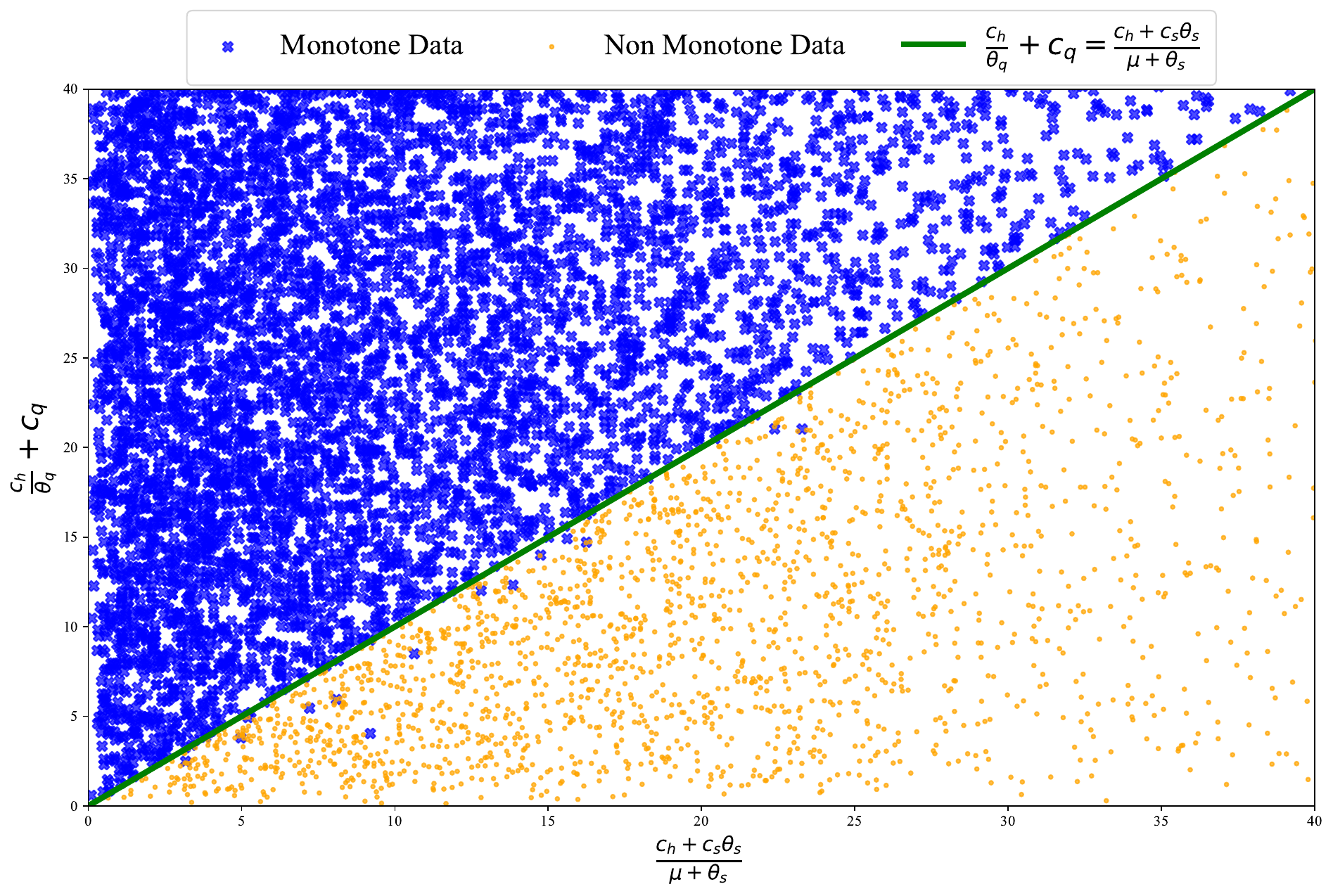}
		\subcaption{Evaluation distribution as $Exp(1/35)$}
        \label{Monotone_condition_exp}
	\end{minipage}
\end{figure}

{From Theorem \ref{Theorem.monotonic}, we know that the optimal policy is always monotone when $\CC_s\leq \CC_q$. However, the system parameters
$\Lambda= 47.08,\, \mu = 42.99,\,N=4,\,m=1,\,\theta_s= 5.24,\,\theta_q=10.7,\, c_s= 27.09,\,c_q= 1.49,\,c_h = 9.48$, and $F(\cdot)\sim U(0,31)$ yield $3.13 \approx \CC_s> \CC_q \approx 2.37$, but the optimal policy is monotone. This suggests that $\mathcal{C}_s\leq \mathcal{C}_q$ may not be a necessary condition for the optimal policy to be monotone.}  
In this section, we investigate the necessity of this condition by conducting numerical experiments.
We consider 10,000 randomly generated sets of system parameters. We select the maximal arrival rate $\Lambda$, service rate $\mu$, service abandonment rate $\theta_s$, waiting abandonment rate $\theta_q$, service abandonment cost $c_s$, waiting abandonment cost $c_q$, and holding cost $c_h$ independently from a uniform distribution with parameters 0 and 50. We choose the evaluation distribution to be either a uniform distribution with range (20,50) or an exponential distribution with mean 35. The buffer size $N$ is a random integer selected from the set of integers $\{2,3,\ldots,20\}$ with equal probability. We start by considering the number of servers to be $m=1$. The results are shown in Figures \ref{Monotone_condition_unif} and \ref{Monotone_condition_exp}.

{Figures \ref{Monotone_condition_unif} and \ref{Monotone_condition_exp} are scatter plots for the randomly generated instances, with each point representing a distinct set of system parameters. The color coding indicates the structural property of the optimal policy for that specific instance (blue cross for monotone, yellow dot for non-monotone).} 
We notice that the line $\mathcal{C}_s = \frac{c_h+ c_s\theta_s}{\mu+\theta_s} = \frac{c_h}{\theta_q} + c_q=\mathcal{C}_q$ divides each plot into two parts. For the part above the line ($\frac{c_h+ c_s\theta_s}{\mu+\theta_s} \leq \frac{c_h}{\theta_q} + c_q$), Theorem \ref{Theorem.monotonic} shows that the optimal policy is monotone. To understand if $ \frac{c_h+ c_s\theta_s}{\mu+\theta_s} \leq \frac{c_h}{\theta_q} + c_q$ is a necessary condition for the optimal policy to be monotone, we focus on the other half below the line ($ \frac{c_h+ c_s\theta_s}{\mu+\theta_s} > \frac{c_h}{\theta_q} + c_q$).

\vspace{-0.4in}

From Figure \ref{Monotone_condition_unif}, we observe that when the evaluation distribution is uniform, there are a few data points with the optimal arrival rate policy being monotone decreasing but $ \frac{c_h+ c_s\theta_s}{\mu+\theta_s} > \frac{c_h}{\theta_q} + c_q$. This indicates that although in general when $ \frac{c_h+ c_s\theta_s}{\mu+\theta_s} > \frac{c_h}{\theta_q} + c_q$ the optimal policy is not monotone, this condition is not necessary.
Figure \ref{Monotone_condition_exp} shows that when the evaluation distribution is exponential, there are also a few data points below the line $\frac{c_h+ c_s\theta_s}{\mu+\theta_s} = \frac{c_h}{\theta_q} + c_q$ where the optimal arrival rate policy is monotone decreasing, but $\frac{c_h+ c_s\theta_s}{\mu+\theta_s} > \frac{c_h}{\theta_q} + c_q$. This numerical result further confirms that our condition in Theorem \ref{Theorem.monotonic} is not necessary.

Zooming into those cases in Figure \ref{Monotone_condition_unif}-\ref{Monotone_condition_exp} with monotone decreasing optimal arrival rate policy where the condition is violated, we discover that they all have a rather small buffer size (i.e., $N-m$ less than or equal to 5). A possible explanation is that the buffer size of those systems is too small to make the arrival rate increase at the beginning. From the intuition discussed in Section 4.1, when $\frac{c_h+ c_s\theta_s}{\mu+\theta_s} > \frac{c_h}{\theta_q} + c_q$, the system could potentially make profit from customers abandoning the system before entering service. However, a small buffer size constrains the queue capacity,  thereby limiting the number of customers who might abandon the system before service. Consequently, the potential profits generated by the pre-service abandonment may not materialize. In that situation, the optimal strategy is to let fewer customers enter when the system is congested, which indicates a monotone decreasing structure for the optimal arrival rate policy. {We also consider the case when the system has multiple servers and provide the corresponding figures in Section EC.3.}

In general, the numerical experiments suggest that the condition is mostly necessary except for a few cases. Thus, these numerical experiments suggest that $\mathcal{C}_s\leq \mathcal{C}_q$ could serve as a rule of thumb for determining if the optimal arrival policy is monotone decreasing (unless the buffer size is small).

\subsection{{Policy Iteration under Uni-modal Structure}} \label{subsec:policy_itration}

\noindent In this section, we present in Algorithm 1 a policy iteration procedure {under uni-modal structure} to solve the Markov decision process defined in Section 3. 
Since the Markov decision process is unichain with a compact action space, the convergence of Algorithm 1 follows from Theorem 4.1 in \citet{Hordijk&Puterman1987}. 
In our pseudo-code, $k\in\BN$ denotes the number of iterations that the algorithm has run.
Let $\Vec{\lambda}(k)=( \lambda_0(k),\ldots,\lambda_{N}(k) )$ be the decision rule in iteration $k$ where $\lambda_n(k)$ is the arrival rate for state $n\in\{0,\ldots,N\}$. Correspondingly, we define $g^k$ as the long-run average gain and vector $\mathbf{h}^{\vec{\lambda}(k)} = (h^k(0),\ldots,h^k(N))$ as the bias function under policy $\vec{\lambda}(k)$ in iteration $k$. {In the ``for" loop, step (b) utilizes the uni-modal structure of the policy. 
Once the arrival rate begins decreasing, it continues to decrease monotonically. Thus,  the action space in state $n+1$ in iteration $k$ reduces to $[0,{\lambda}_k(n)]$ (instead of the entire action space $[0,\Lambda]$) when the arrival rate starts to decrease at $\lambda_{n-1}(k)>\lambda_n(k)$.} {When $\CC_s\leq \CC_q$, from Theorem \ref{Theorem.monotonic} we know that the optimal arrival rate policy is monotone decreasing. Then, we can modify step (b) in Algorithm 1 as ``If $n<N-1$, set $A(n+1) \leftarrow [0,{\lambda}_n(k)]$'' to utilize this structural property.}

\begin{algorithm}
	\caption{Policy Iteration Algorithm under Uni-modal Structure}
	\begin{algorithmic}[1]
		\State Input $N,m,\Lambda,\,\mu,\,\theta_s,\,\theta_q,\, F(\cdot),\,c_h,\,c_s,\,c_q$.
		\State Output: Optimal Policy $\Vec{\lambda}^*$.
		\State Choose policy $\Vec{\lambda}(0)$ where ${\lambda}_n(0) = 0 $ for all $n\in \{0,\ldots,N\}$ and set $k\leftarrow 0$.
		\Repeat
		\State Find scalar $g^k\in\BR$ and vector $\mathbf{h}^{\vec{\lambda}(k)}$ by solving $g^k\mathbf{e} +(\mathbf{I}-\mathbf{P}^{\Vec{\lambda}(k)})\mathbf{h}^{\vec{\lambda}(k)}=\mathbf{r}^{\Vec{\lambda}(k)}$ subject to $h^{k}(0)=0$, where $\mathbf{P}^{\Vec{\lambda}(k)}$ along with $\mathbf{r}^{\Vec{\lambda}(k)}$ are defined in Section 3. Here, $\mathbf{e}$ is a vector of all ones and $\mathbf{I}$ is the identity matrix.
		\State Set $A(0)\leftarrow [0,\Lambda]$ and ${\lambda}_N(k+1)=0.$
		\For{each state $n \in \{0,\ldots,N-1\}$}
		\State (a) Choose action $\lambda \in A(n)$ that maximizes the expected reward plus the bias, i.e.,
		
		\[
		{\lambda}_n(k+1) \leftarrow\inf \Bigg( \argmax_{\lambda\in A(n)}\Big\{ r_n(\lambda)+\sum_{j=0}^N \mathbf{P}(j|n,\lambda) h^k(j)  \Big\} \Bigg). 
		\]
		\State (b) {If $0<n<N-1$ and $\lambda_{n-1}(k)>\lambda_n(k)$, set $A(n+1) \leftarrow [0,{\lambda}_n(k)]$.}
        \State {$\quad \,$ If $n=0$ or both $0<n<N-1$ and $\lambda_{n-1}(k)\leq \lambda_n(k)$, set $A(n+1)\leftarrow [0,\Lambda]$.}
		\EndFor
		\State Set $k \leftarrow k + 1$.
		\Until{{policy $\Vec{\lambda}(k)=\Vec{\lambda}(k-1)$ and set $\Vec{\lambda}^*=\Vec{\lambda}(k)$}}.
		\State \textbf{return} $\Vec{\lambda}^*=\Vec{\lambda}(k)$ as the optimal policy.
	\end{algorithmic}\end{algorithm}

{Since the action space in our problem $[0,\Lambda]$ is not finite, the number of iterations required to achieve exact convergence (satisfying the equality in line 13 of Algorithm 1) can be large. Therefore, when implementing this algorithm, we can either discretize the action space or modify the stopping criterion as $||\Vec{\lambda}(k)-\Vec{\lambda}(k-1)||_2\leq \epsilon$, where $\epsilon$ is the chosen precision for the stopping criterion.
}

{To demonstrate the computational efficiency of Algorithm 1, we compare the running time of finding the optimal policy using Algorithm 1 ($\tau_1$) to that of regular policy iteration ($\tau_0$). 
We calculate the relative percentage improvement (RPI), defined as $(1-\tau_1/\tau_0)\times 100\%$. In our experiments, we use the randomly generated data set in Section \ref{subsec:necessary_condition} and for each set of randomly generated system parameters, we record the running time for finding the optimal policy with Algorithm 1 ($\tau_1$) and the running time for finding the optimal policy with regular policy iteration ($\tau_0$). The resulting RPI from this numerical experiment can be as large as 54.44\%.
We also use Algorithm 1 to find the optimal dynamic policy and the corresponding optimal gain in Section 6. In this case, Algorithm 1 (incorporating the uni-modal structure) performs significantly better than the regular policy iteration algorithm and reduces the computational time 41.05\% on average.}

\section{Heuristics} \label{sec.Heuristics}
Although dynamic pricing always provides the best long-run average gain, the algorithm to find the optimal dynamic pricing policy could be time-consuming even with the monotonicity result, especially when the state space is large. Furthermore, employing a dynamic policy can be impractical. {High-frequency price changes that appear arbitrary to customers may trigger confusion and dissatisfaction, which results in negative publicity when those customers vocalize their complaints.}
Hence, the dynamic pricing policy may be difficult to implement and against public interest. Therefore, it is crucial to find a pricing strategy that is easy to deploy and adapts to real life scenarios with a near-optimal performance.

In this section, we provide two heuristics, the cutoff-static policy and the two-price policy, to fulfill these goals. We first introduce the cutoff-static policy in Section 5.1 and prove properties that help to find the best cutoff-static policy efficiently. We then derive theoretical guarantees on the performance for the best cutoff-static policy. {In Section 5.2, we develop the two-price policy as an improvement of the cutoff-static policy and provide concrete theoretical results on the optimality gap between the best two-price policy and the optimal dynamic policy.}

\subsection{Best Cutoff-static Policy} \label{subsec.cutoff.static}
In this section, we propose a heuristic which is easier to implement than the optimal dynamic policy but yields near optimal performance. The general idea is to combine static pricing and admission control, in which we charge {the same price $p$} to the incoming customers until the number of customers in the system reaches a certain threshold $K$. No customer is allowed to enter the system beyond that threshold. We refer to this heuristic as the cutoff-static policy.  
{Since the price and arrival rate of the system have a one-to-one relationship given by the function $\bar{F}^{-1}$, a fixed price $p$ means a fixed arrival rate $\Lambda\bar{F}(p)=\delta\in[0,\Lambda]$ until the number of customers in the system reaches $K$. }
Here, we do not use $\lambda$ as the arrival rate because we want to distinguish the fixed arrival rate policy in this section from the dynamic arrival rate policy in the previous section.

\begin{definition} \label{def.cut.off.static}
	A cutoff-static policy with a threshold $K\in\{0,\ldots,N\}$ and an arrival rate $\delta \in [0,\Lambda]$, denoted $\pi_{K}^\delta \in \Pi$, is defined as follows:
    {\begingroup\renewcommand{\baselinestretch}{0.8}\normalsize
	$$\pi_{K}^\delta(n) = 
	\begin{cases}
	\delta \quad \text{ for } 0\leq n<K,\\
	0 \quad \text{ for } K\leq n< N. 
	\end{cases}
	$$ 
    \endgroup}
	Correspondingly, we denote the steady-state probability distribution under policy $\pi_{K}^\delta$ as $P_n(\pi_K^\delta)$. To express $P_n(\pi_K^\delta)$ in closed form, we have $a_0(\delta)=1$ for any $\delta\in[0,\Lambda]$ and 
    {\begingroup\renewcommand{\baselinestretch}{0.8}\normalsize
	\[
	P_n(\pi_K^\delta)=
	\begin{cases}
		P_0(\pi_K^\delta)\times a_n(\delta) &\text{for } 0\leq n\leq K,
		\vspace{0.1 in}\\
		 \dfrac{1}{\sum_{n=0}^N a_n(\delta)} &\text{for } n=0,
	\end{cases}
	\text{ where }
	a_n(\delta)  =
	\begin{cases}
		\prod_{i=1}^n \frac{\delta}{\gamma_i} & \text{for } 0\leq n\leq K,\\
		0 & otherwise.
	\end{cases}
	\]
    \endgroup}
Following the formulation in Section 3, we then have the long-run average gain under the cutoff-static policy $\pi_K^\delta$ as 
$
g_C^\delta(K) =r_K(0)\times P_K(\pi_K^\delta) + \sum_{n=0}^{K-1} r_n(\delta)\times P_n(\pi_K^\delta) .
$
\end{definition}

With this definition, we notice that all cutoff-static policies are monotone decreasing since either $\pi_K^\delta(n)=\pi_{K}^\delta(n+1)$ for $n\not = K-1$ or $\pi_{K}^\delta(K-1)=\delta\geq 0=\pi_{K}^\delta(K)$. 
We then define the best cutoff-static policy.

\begin{definition} \label{def.best.cutoff}
\begin{enumerate}
    \item The best cutoff-static policy with threshold $K$, denoted $\pi_C^*(K)$ with arrival rate $\delta_K$, and its gain, denoted $g_C^*(K)$, are defined as 
    $$
   \pi_C^*(K) = \pi_K^{\delta_K} \text{ where } \delta_K = \inf\left\{\argmax_{\delta\in [0,\Lambda]} g_C^\delta(K)\right\} \text{ and } g_C^*(K)=\max_{\delta\in[0,\Lambda]} g_C^\delta(K).
    $$
    If we choose $K=N$, the best cutoff-static policy {becomes the best static policy for the original problem}. We denote this gain as the best static gain $g_S^*$ where $g_S^*=g_C^*(N) = \max_{\delta\in [0,\Lambda]} g_C^\delta (N).$
    \item The best cutoff-static policy achieves the largest gain $g_C^*$ among all cutoff-static policies
    \[
    g_{C}^* = \max_{K=0,\ldots,N} \Big\{ \max_{\delta\in[0,\Lambda]} g_C^\delta(K) \Big\}.
    \]
     We refer to the threshold of the best cutoff-static policy $K_C$ as the optimal threshold in which $K_C = \min\left\{ \argmax_{K=0,\ldots,N} g_{C}^*(K) \right\}$.

Correspondingly, we denote the best cutoff-static policy as $\pi_{C}^*$.
\end{enumerate}

\end{definition}

{Now we discuss how to find the optimal threshold $K_C$ efficiently. Instead of enumerating all possible thresholds from $0$ to $N$ and solving for the best cutoff-static policy with that threshold, we provide lower bounds and upper bounds for the optimal threshold $K_C$ under some conditions to speed up this procedure. In the following theorem, we show that $K_C\geq m$ except for the case in which the best cutoff static policy has a zero gain. }

\begin{theorem} \label{Theorem.skip.server}
When $g^*_C(1)>0$, then $g^*_C(K+1)\geq g^*_C(K)$ for all $K=1,\ldots,m-1$. When $g^*_C(1)= 0$, then $g^*_C({K})=0$ for all $K=1,\ldots,m$.
\end{theorem}

From Theorem \ref{Theorem.skip.server}, when $g_C^*(1)>0$, we have $g^*_C(K+1)\geq g^*_C(K)$ for $K=1,\ldots,m-1$, indicating that $g^*_C(m)=\max_{K=0,\ldots,m}g^*_C(K)$ which means $K_C\geq m$ { and $m$ is a lower bound of $K_C$. When $g^*_C(1)=0$, we have $g^*_C(K)=0$ for all $K=1,\ldots,m$, indicating that $g^*_C(K)$ may be positive only when $K> m$. If $g_C^*(K)$ is never positive when $K>m$, we have $g_C^*=0$ and $K_C=0$. Otherwise, we have $K_C>m$. This implies that $K_C\geq m$ except for the case when the best cutoff static policy has a zero gain. Therefore, Theorem \ref{Theorem.skip.server} implies that $K_C$ can never take values between $1$ to $m-1$.}

Furthermore, Theorem \ref{Theorem.skip.server} also is intuitive. Namely, we either use all the servers available or idle all of them. When $g_C^*(1)>0$, there exists a pricing strategy that can make the system profitable by utilizing one server. Then, logically, if the system has more than one server, it makes sense to assume that a similar pricing strategy will make the system profitable when all servers are busy. In other words, there should not be any forced server idleness if operating only one server is already profitable. In comparison, if $g_C^*(1)=0$, it means the system can never be profitable when operating only one server. Then, it is reasonable to conclude that idling all servers will not change the fact that the system is not profitable. {We prove Theorem \ref{Theorem.skip.server} in detail in Section \ref{ec.proof.section.5.1}.}

{Now we provide an upper bound for the optimal threshold $K_C$. In Proposition \ref{prop.upper.0}, we identify conditions under which $K_C=0$. Proposition \ref{prop.upper.0}  also shows that the system is never profitable when the maximal willingness to pay of a customer does not exceed  $\min \left\{ \mathcal{C}_s=\frac{c_h+c_s\theta_s}{\mu+\theta_s}, \mathcal{C}_q= \frac{c_h}{\theta_q}+c_q \right\}$.}

\begin{proposition} \label{prop.upper.0}
If $\Bar{F}\left(\min \left\{ \mathcal{C}_s,\mathcal{C}_q \right\}  \right)=0$, then $K_C=0$ and $g^*= 0$.
\end{proposition}

{ Proposition \ref{prop.upper.0} provides a sufficient condition for $K_C=0$ and $g^*=0$ (i.e., the system can never be profitable). We provide an intuitive explanation for this condition. 
Recall from Section 4.1 that $\CC_s$ can be interpreted as the expected total cost after joining service and $\CC_q$ as an upper bound of the expected total cost of a customer before joining service. When $\Bar{F}\left(\min \left\{ \mathcal{C}_s,\mathcal{C}_q \right\}  \right)=0$, the maximum willingness to pay of a customer does not exceed $\CC_s$ and $\CC_q$.
This suggests that the profit generated from serving a customer will not exceed the expected cost of serving that customer.
Thus, it is never profitable to accept customers and $K_C=0$.
The following proposition further refines the upper bound $K_C$ when $\CC_s\leq \CC_q$.

\begin{proposition} \label{prop.upper.refine}

When $\CC_s\leq \CC_q$, if there exists $b\geq m$ such that $\Bar{F}\left( \frac{\gamma_b}{\gamma_{b+1}} \mathcal{C}_s +  \frac{\gamma_{b+1}-\gamma_m}{\gamma_{b+1}} \mathcal{C}_q \right)=0$, then $K_C \leq b$.
\end{proposition}
}

{We prove Propositions \ref{prop.upper.0} and \ref{prop.upper.refine} in Section \ref{ec.proof.section.5.1}. Using Theorem \ref{Theorem.skip.server}, and Propositions \ref{prop.upper.0} and \ref{prop.upper.refine}, we can efficiently find the best cutoff-static policy. We consider all values of $K$ from $m$ to the upper bounds derived in Propositions \ref{prop.upper.0} and \ref{prop.upper.refine}, and conduct a one-variable optimization problem $g_C^*(K) = \max_{\delta\in[0,\Lambda]} g_C^{\delta}(K)$ for each threshold $K$ to find the best cutoff-static gain $g^*_C(K)$ and the corresponding arrival rate $\delta_K$. We then pick the largest $g^*_C(K)$ to achieve the best cutoff-static policy and the corresponding $\delta_K$ as the optimal arrival rate. If $g^*_C(K)\leq 0$ for all $K\geq m$, it means that system is not profitable and $K_C=0$.}

The rest of this section focuses on showing when the best cutoff-static policy is near optimal. First, we introduce the concept of a regular distribution, which is a type of probability distributions used to model customer demand in revenue management (also referred to as Myerson's regularity). 

\begin{definition} \label{def.regular.distribution}
An evaluation distribution with a continuous cumulative distribution function $F(\cdot)$ is a regular distribution if $\lambda\Bar{F}^{-1}(\lambda/\Lambda)$ is a continuous concave function for $\lambda\in [0,\Lambda].$ 
\end{definition}

Regular distributions include many commonly used distributions for revenue management, such as the Gaussian, uniform, and exponential distributions. Therefore, it is reasonable for us to assume our evaluation distribution is regular. 

We also define new notations for the long-run average revenues and costs for our Markov decision process. 
For a given policy $\Vec{\lambda}\in\Pi$, define the long-run average revenue as
$$
\mathbf{R}(\Vec{\lambda}) = \sum_{n=0}^N \lambda_n\Bar{F}^{-1}(\lambda_n/\Lambda) \times P_n(\Vec{\lambda})
$$
and the long-run average total cost as 
\[
\mathbf{C}(\Vec{\lambda}) = \sum_{n=1}^N (nc_h + \min\{n,m\}c_s\theta_s + \max\{n-m,0\}c_q\theta_q)\times P_n(\Vec{\lambda}). \addrefnum\label{total.exp.cost}
\]
Consequently, we have $g^{\Vec{\lambda}} = \mathbf{R}(\Vec{\lambda})-\mathbf{C}(\Vec{\lambda})$ for any $\Vec{\lambda}\in \Pi$. Then, $\mathbf{R}(\Vec{\lambda}^*)$ and $\mathbf{C}(\Vec{\lambda}^*)$ are the long-run average revenue and cost under the optimal dynamic policy $\vec{\lambda}^*$ and $\mathbf{R}(\pi_K^\delta)$ and $\mathbf{C}(\pi_K^\delta)$ are those under a cutoff-static policy $\pi_K^\delta$. The following results link $\mathbf{R}(\Vec{\lambda}^*)$ and $\mathbf{C}(\Vec{\lambda}^*)$ with $\mathbf{R}(\pi_K^{\hat{\delta}})$ and $\mathbf{C}(\pi_K^{\hat{\delta}})$ where $\hat{\delta}$ is defined to be the effective (long-run average) arrival rate of the system under the optimal dynamic policy $\Vec{\lambda}^*$:
\[
\hat{\delta} = \sum_{n=0}^N \lambda_n^*\times P_n(\Vec{\lambda}^*)  \addrefnum \label{effective.arrival.rate}.
\] 

\begin{theorem} \label{Theorem.performance.reve}
When the evaluation distribution is a regular distribution, for any $m\leq K\leq N$, the policy $\pi_K^{\hat{\delta}}$ satisfies 
\[
\mathbf{R}(\pi_K^{\hat{\delta}})  \geq \mathbf{R}(\Vec{\lambda}^*) \times \left(1-\dfrac{a_K(\Lambda)}{\sum_{n=0}^K a_n(\Lambda)}\right).
\]
\end{theorem}

\begin{theorem} \label{Theorem.performance.cost}
	For any $m\leq K\leq N$,   
	\begin{enumerate}
		\item  if $c_q\theta_q\geq c_s\theta_s>0$, the policy $\pi_K^{\hat{\delta}}$ satisfies
		\[
		\mathbf{C}(\pi_K^{\hat{\delta}})\leq \mathbf{C}(\Vec{\lambda^*})\times \max_{\delta\in [0,\Lambda]} \left( \dfrac{\sum_{n=1}^K n  P_n(\pi_K^{{\delta}})}{{\delta}/\max\{\mu+\theta_s, \theta_q\}} + \dfrac{c_q\theta_q-c_s\theta_s}{c_h+c_s\theta_s}\times \dfrac{\sum_{n=m+1}^K n  P_n(\pi_K^{{\delta}})}{{\delta}/\max\{\mu+\theta_s, \theta_q\}}  \right).
		\]
		\item if $c_s\theta_s\geq c_q\theta_q>0$, the policy $\pi_K^{\hat{\delta}}$ satisfies
		\[
		\mathbf{C}(\pi_K^{\hat{\delta}})\leq  \mathbf{C}(\Vec{\lambda^*})\times \max_{\delta\in [0,\Lambda]} \left( \dfrac{\sum_{n=1}^K n  P_n(\pi_K^{{\delta}})}{{\delta}/\max\{\mu+\theta_s, \theta_q\}} + \dfrac{c_s\theta_s-c_q\theta_q}{c_h+c_q\theta_q}\times \dfrac{\sum_{n=1}^K \min\{n,m\}  P_n(\pi_K^{{\delta}})}{{\delta}/\max\{\mu+\theta_s, \theta_q\}}  \right).
		\]
	\end{enumerate}
\end{theorem}

Since we have a lower bound on the revenue and an upper bound on the cost, Theorems \ref{Theorem.performance.reve} and \ref{Theorem.performance.cost} provide a lower bound on the long-run average gain under the cutoff-static policy $\pi_K^{\hat{\delta}}$ with threshold $K$ and arrival rate $\hat{\delta}$ defined in (\ref{effective.arrival.rate}).
{We prove Theorems \ref{Theorem.performance.reve} and \ref{Theorem.performance.cost} in Section \ref{ec.proof.section.5.1} and we discuss further insights regarding these theorems in the following remark.
}

\begin{remark} \label{rem.intuition.Theorem.performance}

To illustrate Theorems \ref{Theorem.performance.reve} and \ref{Theorem.performance.cost}, consider the case when $K=m=1$ and $0<c_s\theta_s\leq c_q\theta_q$:
\[
\begin{cases}
	\mathbf{R}(\pi_1^{\hat{\delta}})  \geq \mathbf{R}(\Vec{\lambda}^*) \times \left(1-\dfrac{\frac{\Lambda}{\mu+\theta_s}}{1+\frac{\Lambda}{\mu+\theta_s}}\right) = \mathbf{R}(\Vec{\lambda}^*) \times \dfrac{\mu+\theta_s}{\Lambda+\mu+\theta_s}, 
	\vspace{0.1in}\\
	\mathbf{C}(\pi_1^{\hat{\delta}})\leq  \mathbf{C}(\Vec{\lambda^*})\times \max_{\delta\in [0,\Lambda]} \dfrac{({\frac{\delta}{\mu+\theta_s}})/({1+\frac{\delta}{\mu+\theta_s}})}{\delta/\max\{\mu+\theta_s, \theta_q\}}  = \mathbf{C}(\Vec{\lambda^*})\times \dfrac{\max\{\mu+\theta_s, \theta_q\}}{\mu+\theta_s}.
\end{cases}
\]
We can now bound the gain of the best cutoff-static policy in terms of the revenue and the total cost of the optimal dynamic policy as follows:
\begin{align*}
    g^*_C =& \max_{K=0,\ldots,N} \Big\{ \max_{\delta\in[0,\Lambda]} g_C^\delta (K) \Big\} \geq g_C^{\hat{\delta}}(1) = \mathbf{R}(\pi_1^{\hat{\delta}}) - \mathbf{C}(\pi_1^{\hat{\delta}})\\
    \geq & \mathbf{R}(\Vec{\lambda}^*) \times \dfrac{\mu+\theta_s}{\Lambda+\mu+\theta_s} - \mathbf{C}(\Vec{\lambda^*})\times \dfrac{\max\{\mu+\theta_s, \theta_q\}}{\mu+\theta_s}.   \addrefnum \label{tightness.analysis}
\end{align*}
{ If we further assume $\mu+\theta_s\geq \theta_q$, then (\ref{tightness.analysis}) can be reduced to $g^*_C \geq \mathbf{R}(\Vec{\lambda}^*) \times \frac{\mu+\theta_s}{\Lambda+\mu+\theta_s} - \mathbf{C}(\Vec{\lambda^*})$ with $K=m=1$. If we let service rate $\mu$ increases, indicating a growth in the processing power of the system, then $\frac{\mu+\theta_s}{\Lambda+\mu+\theta_s}$ will be close to 1 
and $g_C^*\approx \mathbf{R}(\vec{\lambda}^*) -  \mathbf{C}(\vec{\lambda}^*) = g^*$. This observation suggests that when a single-server system has  $0<c_s\theta_s\leq c_q\theta_q$ and $\mu+\theta_s\geq \theta_q$ and the service rate approaches to infinity, picking the threshold $K=1$ leads to a near optimal heuristic that is easy to implement.}
\end{remark}

{

Indeed, numerical experiments confirm our intuitive argument in Remark \ref{rem.intuition.Theorem.performance} and the tightness of Theorems \ref{Theorem.performance.reve} and \ref{Theorem.performance.cost}. Under a numerical example of a single-server system with parameters $\Lambda = 10$, $N=6$, $\theta_s=10$, $\theta_q=10$, $c_s=3$, $c_q=27$, $c_h=11$, and an evaluation distribution of $Exp(35)$, we increase the service rate $\mu$ from $10$ to $10^4$. Given these system parameters satisfy $0<c_s\theta_s\leq c_q\theta_q$ and $\mu+\theta_s\geq \theta_q$, if we pick the threshold $K=m=1$, we can simplify the statement in Theorems \ref{Theorem.performance.reve} as $\mathbf{R}(\pi_1^{\hat{\delta}})  \geq  \mathbf{R}(\Vec{\lambda}^*) \times \frac{\mu+\theta_s}{\Lambda+\mu+\theta_s}$ and that in Theorem \ref{Theorem.performance.cost} as $\mathbf{C}(\pi_1^{\hat{\delta}})\leq \mathbf{C}(\Vec{\lambda^*})\times \frac{\max\{\mu+\theta_s, \theta_q\}}{\mu+\theta_s} = \mathbf{C}(\Vec{\lambda^*})$. Figures \ref{fig:tightness_revenue} and \ref{fig:tightness_cost} demonstrate the relationship between $\mathbf{R}(\pi_1^{\hat{\delta}})$ and $\mathbf{R}(\Vec{\lambda}^*) \times \frac{\mu+\theta_s}{\Lambda+\mu+\theta_s}$, and between $\mathbf{C}(\pi_1^{\hat{\delta}})$ and $\mathbf{C}(\Vec{\lambda^*})$ as a function of $\mu$ (in log scale). 

\begin{figure}[h]
	\centering
	\caption{Tightness of Theorems \ref{Theorem.performance.reve} and \ref{Theorem.performance.cost} as the service rate $\mu$ grows}
	\begin{minipage}{0.45\textwidth}
		\centering
		\includegraphics[width=\columnwidth]{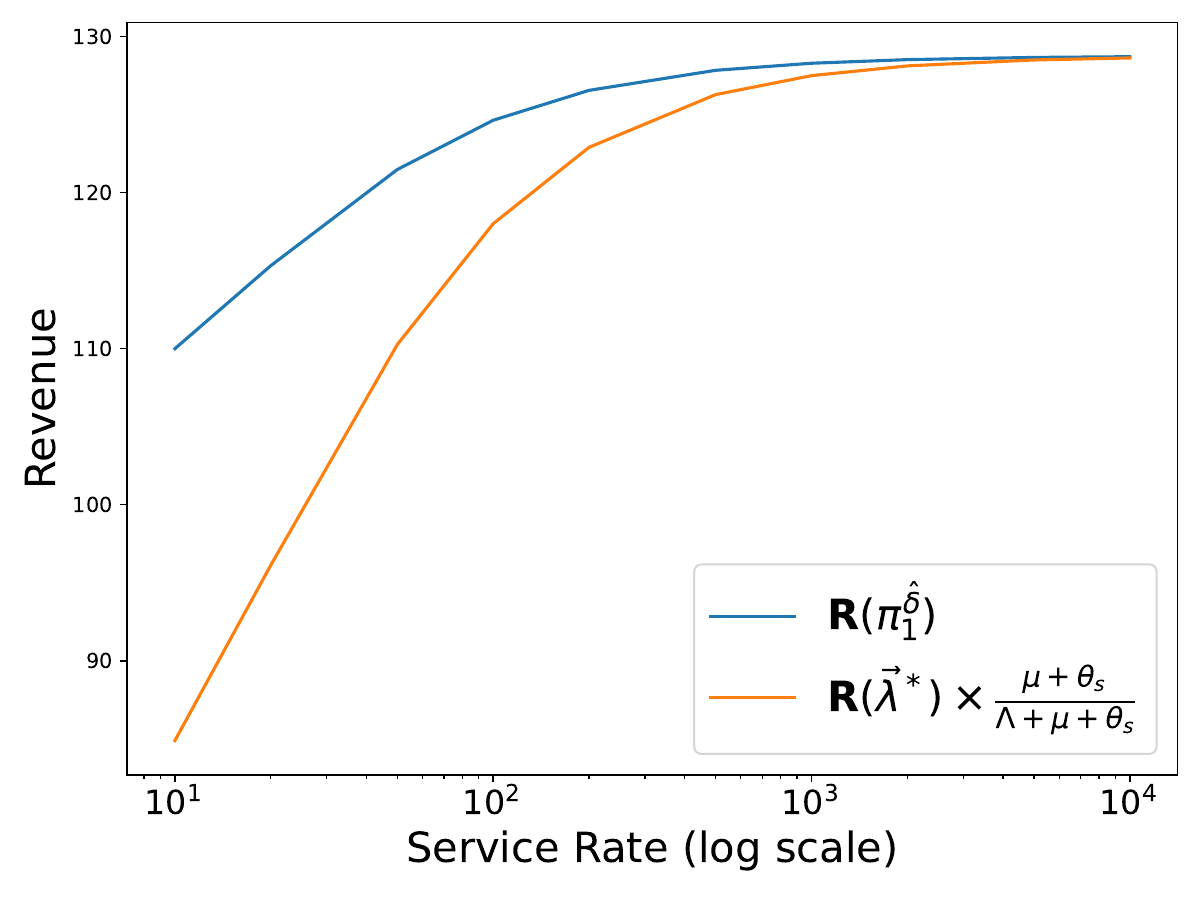}
		\subcaption{Tightness for revenue}
        \label{fig:tightness_revenue}
	\end{minipage}
	\begin{minipage}{0.45\textwidth}
		\centering
		\includegraphics[width=\columnwidth]{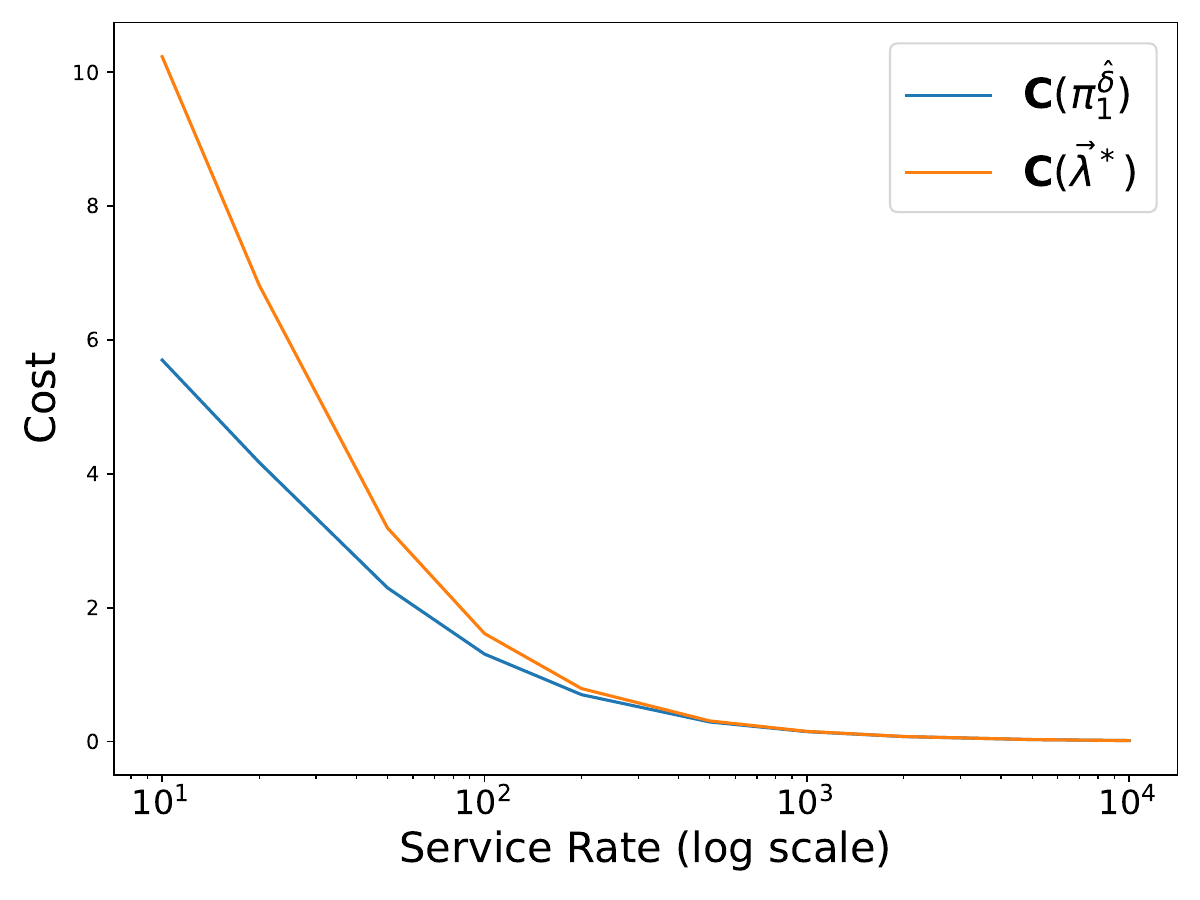}
		\subcaption{Tightness for cost}
        \label{fig:tightness_cost}
	\end{minipage}
\end{figure}

In Figure \ref{fig:tightness_revenue}, we observe that $\mathbf{R}(\pi_1^{\hat{\delta}})$ and $\mathbf{R}(\Vec{\lambda}^*) \times \frac{\mu+\theta_s}{\Lambda+\mu+\theta_s}$ converge as the service rate $\mu$ approaches infinity, suggesting the tightness of Theorem \ref{Theorem.performance.reve}. In Figure \ref{fig:tightness_cost}, we observe that $\mathbf{C}(\pi_1^{\hat{\delta}})$, and $\mathbf{C}(\Vec{\lambda^*})$ converge as the service rate $\mu$ goes to infinity, suggesting the tightness of Theorem \ref{Theorem.performance.cost}.

Furthermore, we show that the best cutoff-static gain achieves at least $\frac{15}{19}\approx78.9\%$ of the optimal dynamic gain under the conditions in Proposition \ref{prop.better.bound}, which is proven in Section \ref{ec.proof.section.5.1}.

\begin{proposition}  \label{prop.better.bound}
When $\bar{F}(\CC_s)>0$ and $\CC_s\leq \CC_q$, if $\Bar{F}\left( \frac{\gamma_m}{\gamma_{m+1}} \CC_s +  \frac{\gamma_{m+1}-\gamma_m}{\gamma_{m+1}}\CC_q \right)=0$, then $\frac{g_C^*}{g^*}\geq \frac{15}{19}$.
\end{proposition}
}

\subsection{Best Two-price Policy} \label{subsec.two.price}

One of the drawbacks of the best cutoff-static policy is that it may yield a zero gain when the optimal dynamic gain is positive. {For instance when $\Lambda = 13.65,\,N=8,\,\mu=1.55,\,\theta_s=4.68,\,\theta_q=46.05,\,c_s=26.92,\,c_q=4.47,$ and 	$c_h=16.59$ with the evaluation distribution $U(20,50)$, the best cutoff-static policy gives a zero gain while the gain of the optimal dynamic policy is 3.49.} In order to provide a better heuristic for those cases, we introduce a new heuristic: the two-price policy. {For a two-price policy, the service provider charges the incoming customers one price $p_s$ when there are idling servers in the system and another price $p_q$ when the customers need to wait in the queue. We also do not accept arrivals when the number of customers reaches a certain threshold, which we know from 
Theorem \ref{Theorem.skip.server} is at least $m$ except for the case when the best cutoff static policy has a zero gain. We first present a general definition for the two-price policy in Definition \ref{def.two.price} using the arrival rates $\delta_s$ and $\delta_q$ corresponding to the prices $p_s$ and $p_q$, respectively.}

\begin{definition} \label{def.two.price}
A two-price policy with a threshold $K\in\{m,\ldots,N\}$ and arrival rates $\delta_s,\,\delta_q  \in [0,\Lambda]$, denoted $\pi_{K}^{\delta_s,\delta_q} \in \Pi$, is defined as follows:
{\begingroup\renewcommand{\baselinestretch}{0.8}\normalsize
\[
\pi_{K}^{\delta_s,\delta_q}(n) =
\begin{cases}
\delta_s \qquad \text{for } 0\leq n < m,\\
 \delta_q \qquad \text{for } m\leq n < K,\\
 0\, \qquad \text{ for } K\leq n\leq N,
\end{cases}
\] 
\endgroup}
where we refer to $\delta_s$ as the no-wait arrival rate (the arrival rate for the customers who can immediately enter the service) and $\delta_q$ as the wait arrival rate (the arrival rate for the customers who need to wait in the queue). We denote the gain of policy $\pi_K^{\delta_s,\delta_q}$ as $g^{\delta_s,\delta_q}_T(K)$. When $\delta_s=\delta_q$, the two-price policy becomes the best cutoff-static policy with only one price.
\end{definition}

For any $\delta_s,\delta_q\in[0,\Lambda]$, we denote the steady-state probabilities under policy $\pi_{K}^{\delta_s,\delta_q}$ as
{\begingroup\renewcommand{\baselinestretch}{0.8}\normalsize
\[
P_n(\pi_{K}^{\delta_s,\delta_q}) = 
\begin{cases}
P_0(\pi_{K}^{\delta_s,\delta_q}) \times a_n(\delta_s,\delta_q) & \text{for } 0\leq n\leq K,
\vspace{0.1 in}\\
\dfrac{1}{\sum_{n=0}^N a_n(\delta_s,\delta_q)} & \text{for } n=0,
\end{cases}
\]
where $a_0(\delta_s,\delta_q) = 1$ and

\[
a_n(\delta_s,\delta_q) =
\begin{cases}
	\prod_{i=1}^n \frac{\delta_s}{\gamma_i} & \text{for } 0\leq n\leq m, \vspace{0.1 in} \\
	a_m(\delta_s,\delta_q) \prod_{i=m+1}^K \frac{\delta_q}{\gamma_i} & \text{for } m<n\leq K, \vspace{0.1 in} \\
	0 & \text{otherwise}.
\end{cases}
\]
\endgroup}
The long-run average gain under this policy is 
\[
g_T^{\delta_s,\delta_q}(K) = r_K(0) \times P_K(\pi_{K}^{\delta_s,\delta_q}) + \sum_{n=0}^{m-1} r_n(\delta_s)\times P_n(\pi_{K}^{\delta_s,\delta_q}) + \sum_{n=m}^{K-1} r_n(\delta_q)\times P_n(\pi_{K}^{\delta_s,\delta_q}) .
\]

{The flexibility of charging two prices allows this policy to mimic all the possible structures that the optimal pricing policy can have. When $\delta_s\geq \delta_q$, according to Definition \ref{def.mono}, the arrival rate policy $\pi_{K}^{\delta_s,\delta_q}$ has a monotone decreasing structure and the corresponding pricing policy has a monotone increasing structure.} When $\delta_s<\delta_q$, the two-price policy has the uni-modal structure (Definition \ref{def.uni.modal}) with $\hat{n}=m$. From Theorems \ref{Theorem.uni_modal} and \ref{Theorem.monotonic}, we know the optimal dynamic policy always has uni-modal structure and it is monotone decreasing {when $\CC_s\leq \CC_q$}. Therefore, by changing the relation of $\delta_s$ and $\delta_q$, the two-price policy can depict all the possible structures of the optimal dynamic policy. {Now, we define the best two-price policy.}

\begin{definition}
	\begin{enumerate}
		\item The best two-price policy with threshold $K$ and its gain $g_T^*(K)$ are defined as $g_T^*(K)=\max_{(\delta_s,\delta_q)\in [0,\Lambda]^2} g_T^{\delta_s,\delta_q}(K)$ and 
		$$
		\pi_{T}^*(K) = \pi_K^{\delta_s(K),\delta_q(K)} \text{ where } \left(\delta_s(K),\delta_q(K)\right) \in \argmax_{(\delta_s,\delta_q)\in [0,\Lambda]^2} g_T^{\delta_s,\delta_q}(K) 
		$$
		\item The best two-price policy achieves the largest gain $g_T^*$ among all two-price policies
		\[
		g_T^* = \max_{K=m,\ldots,N} \Big\{ \max_{(\delta_s,\delta_q)\in [0,\Lambda]^2} g_T^{\delta_s,\delta_q}(K) \Big\}.
		\]
		{We call the threshold of the best two-price policy as the optimal two-price threshold in which $K_T = \min\left\{ \argmax_{K=m,\ldots,N} g_T^{\delta_s,\delta_q}(K)\right\}$.
		We denote the best two-price policy as $\pi_T^*$ with optimal threshold $K_T$ and arrival rates $\delta_s^*,\,\delta_q^*$.}
	\end{enumerate}
	
\end{definition}

{Now, we focus on how to find the best two-price policy. 
We enumerate all the possible $K$ values from $m$ to $N$ and solve a two-variable optimization problem to find the best two-price gain $g_T^*(K)$ for each $K$. Finally, we select the maximum among all $g_T^*(K)$ to get $g_T^*$, the threshold $K_T$ with the maximum two-price gain, and the prices $p_s$ and $p_q$ corresponding to $\delta_s$ and $\delta_q$ are the best two prices. {If $g_T^*(K)\leq 0$ for all $K\geq m$, we have $g_T^*=0$ and let $K_T=0$.} Using this procedure, we compute an optimal two-price policy for the system with $\Lambda = 13.65,\,N=8,\,\mu=1.55,\,\theta_s=4.68,\,\theta_q=46.05,\,c_s=26.92,\,c_q=4.47,$ and 	$c_h=16.59$ mentioned at the beginning of this section. The gain of the best two-price policy is 3.48 with $\delta_s=1.47,\,\delta_q=4.91$ and $K_T=8$, which is a significant improvement over the best cutoff-static policy with gain 0 and very close to the optimal gain 3.49. This improvement is not limited to just a specific instance. From the numerical experiments in Section 6, we observe that the difference in the relative performance of the best two-price and cutoff-static policies is significant.} 

{  {As a direct generalization of the cutoff static policy, the performance analysis in Section \ref{subsec.cutoff.static} about the best cutoff static policy carries on to the best two-price policy.}
To further demonstrate the strength of the two-price policy by providing theoretical guarantees, we present the following result that gives an upper bound on the optimality gap between the two-price policy and the optimal policy. In order to formally present the result, we define the following short-hand notations. When the evaluation distribution $F(\cdot)$ is a regular distribution (namely $\lambda\bar{F}^{-1}(\lambda/\Lambda)$ is concave with respect to $\lambda$), let $\tilde{\CR}_s$ and $\tilde{\delta_s}$ be the myopic reward and rate for service, respectively. Then, \[
\tilde{\mathcal{R}_s} = \max_{\lambda\in[0,\Lambda]} \left\{\lambda\Bar{F}^{-1}(\lambda/\Lambda) - \lambda \mathcal{C}_s   \right\} \text{ and } \tilde{\delta_s} = \inf \left(\argmax_{\lambda\in[0,\Lambda]} \left\{\lambda\Bar{F}^{-1}(\lambda/\Lambda) - \lambda \mathcal{C}_s   \right\} \right).
\]
Similarly, define the myopic reward $\tilde{\CR}_q$ and rate $\tilde{\delta_q}$ for abandoning as \[
\tilde{\mathcal{R}_q} = \max_{\lambda\in[0,\Lambda]} \left\{\lambda\Bar{F}^{-1}(\lambda/\Lambda) - \lambda \mathcal{C}_q   \right\} \text{ and } \tilde{\delta_q} = \inf \left(\argmax_{\lambda\in[0,\Lambda]} \left\{\lambda\Bar{F}^{-1}(\lambda/\Lambda) - \lambda \mathcal{C}_q   \right\} \right).
\]
Notice that if $\lambda\bar{F}^{-1}(\lambda/\Lambda)$ is concave, when $\CC_s\leq \CC_q$, we have $\tilde{\CR}_s\geq \tilde{\CR}_q$, and when $\CC_s\geq \CC_q$, $\tilde{\CR}_s \leq \tilde{\CR}_q$. 
We also use $P_s^T$, $P_q^T$, and $P_N^T$ as the probability that a server is idle, the probability that all servers are busy when the system is not full, and the probability of balking, respectively, under the two-price policy $\pi_N^{\tilde{\delta_s},\tilde{\delta_q}}$ with no-wait arrival rate $\tilde{\delta_s}$, wait arrival rate $\tilde{\delta_q}$ and a threshold $N$:
\[
P_s^T = \sum_{n=0}^{m-1} P_n\left(\pi_N^{\tilde{\delta_s}, \tilde{\delta_q}}\right) ,\, P_q^T = \sum_{n=m}^{N-1} P_n\left(\pi_N^{\tilde{\delta_s}, \tilde{\delta_q}}\right), \text{ and } P_N^T = P_N\left(\pi_N^{\tilde{\delta_s}, \tilde{\delta_q}}\right).
\]
Finally, let $P_s^\Lambda$ be the probability that a server is idle and $P_q^\Lambda$ be the probability that all servers are busy under the static policy with the maximal arrival rate $\Lambda$ when $N$ tends to $ \infty$. Then, 
\[
P^{\Lambda}_s =  \dfrac{\sum_{n=0}^{m-1} a_n(\Lambda)}{\sum_{n=0}^{\infty} a_n(\Lambda)} \text{ and } P_q^\Lambda = \dfrac{\sum_{n=m}^{\infty} a_n(\Lambda)}{\sum_{n=0}^{\infty} a_n(\Lambda)}. \addrefnum \label{infinite buffer steady state}
\]
The existence of $P_s^\Lambda$ and $P_q^\Lambda$ is guaranteed and the closed form expressions are provided in Appendix B of \citet{Garnett.et.al.2002}. 
We next state our result that provides upper bounds for the optimality gap between the best two-price policy and the optimal dynamic policy.
\begin{theorem} \label{Theorem.tp.perform}
If the evaluation distribution $F(\cdot)$ is a regular distribution, we have 
\[
g^*-g_T^* \leq
\begin{cases}
\left(\tilde{\CR}_s - \tilde{\CR}_q\right) P_q^T + \tilde{\CR}_s P_N^T  &\text{when } \mathcal{C}_s\leq \mathcal{C}_q,  \\
\left(\tilde{\CR_q} - \tilde{\CR_s}\right) \left(P_q^\Lambda - P_q^T\right) + \tilde{\delta_q} (\mathcal{C}_s-\mathcal{C}_q) P_q^T + \tilde{\CR}_s P_N^T &\text{when } \mathcal{C}_s \geq \mathcal{C}_q.
\end{cases}
\]
\end{theorem}
In Theorem \ref{Theorem.tp.perform}, $\left(\tilde{\CR}_s - \tilde{\CR}_q\right) P_q^T$ when $\CC_s\leq \CC_q$ and $\left(\tilde{\CR_q} - \tilde{\CR_s}\right) \left(P_q^\Lambda - P_q^T\right) + \tilde{\delta_q} (\mathcal{C}_s-\mathcal{C}_q) P_q^T$ when $\CC_s\geq \CC_q$ represent the approximation errors induced by using the best two-price policy  and $\tilde{\CR}_s P_N^T$ can be viewed as an error induced by having a finite buffer. 
The detailed proof of Theorem \ref{Theorem.tp.perform} is in Section \ref{ec.proof.section.5.2}.

To further demonstrate the strength of Theorem \ref{Theorem.tp.perform}, in Proposition \ref{Prop.tp.optimal} we show that when $\CC_s=\CC_q$, the bounds in this result are tight asymptotically as $N\to \infty$ and the best two-price policy is optimal. The proof is included in Section \ref{ec.proof.section.5.2}.
\begin{proposition} \label{Prop.tp.optimal}
When $\mathcal{C}_s=\mathcal{C}_q$ and $N\to \infty$, we have $g^*-g_T^* \to 0$, and, hence, the gain of the best two-price policy converges to the gain of the optimal dynamic policy.
\end{proposition}
 
}

\section{Numerical Experiments} \label{sec.numerical.example}

In this section, we present numerical experiments documenting the performance of the best cutoff-static policy and the best two-price policy. 
{We employ the best static policy with gain $g_S^*$ as a benchmark. From Definition 5.2, the best static policy is the best cutoff-static policy with threshold $N$. To find the best static policy, we solve the one-variable optimization problem $\max_{\delta\in[0,\Lambda]} g_C^{\delta}(N)$. To evaluate the performance of our heuristics, we use a performance ratio defined as follows. }

{
\begin{definition}
Let $g^*$ denote the optimal long-run average profit. The performance ratio for a heuristic policy $H \in \{S, C, T\}$ (representing the Static, Cutoff-static, and Two-price policies, respectively) is defined as $R_H = g_H^*/g^*$, where $g_H^*$ represents the long-run average profit achieved by the best policy of type $H$.
\end{definition}
 }

In the following numerical experiments, we randomly generate 10,000 sets of system parameter in the same manner as in Section \ref{subsec:necessary_condition}. In particular, we randomly select the maximal arrival rate $\Lambda$, service rate $\mu$, service abandonment rate $\theta_s$, waiting abandonment rate $\theta_q$, service abandonment cost $c_s$, waiting abandonment cost $c_q$ and holding cost $c_h$ independently from a uniform distribution with parameters 0 and 50. The buffer size $N$ is a random integer selected from the set of integer $\{2,\ldots,20\}$ with equal probability and we choose the evaluation distribution to be either a uniform distribution with range (20,50) or an exponential distribution with mean 35. We fix the number of servers $m$ to be 1. For each set of system parameters, we compute the optimal dynamic gain and the gain under the heuristics to estimate the performance ratios $R_S,\,R_C,$ and $R_T$. Table \ref{Histogram_1_server} shows the frequency of performance ratios in various ranges when the evaluation distribution is uniform.

{\begingroup\renewcommand{\baselinestretch}{0.6}\normalsize

\begin{table}[h]   
	\centering
	\caption{Performance ratios of heuristics with evaluation distribution as $U(20,50)$}
	\begin{minipage}{0.5\textwidth} 
		\centering
		\renewcommand{\arraystretch}{1.5} 
		\begin{tabular}{|>{}c|>{}c|>{}c|>{}c|}  
			\hline
			{\tiny Performance ratio} & {\tiny Best static policy} & {\tiny Best cutoff-static policy} & {\tiny Best two-price policy}\\ \hline
			{\small $0\sim0.1$} & 10 & 10 & 0\\ \hline
			{\small $0.1\sim0.2$} & 1 & 1 & 0 \\ \hline
            {\small $0.2\sim0.3$} & 1 & 1 & 0 \\ \hline
            {\small $0.3\sim0.4$} & 4 & 1 & 0\\ \hline
            {\small $0.4\sim0.5$} & 14 & 1 & 0 \\ \hline
            {\small $0.5\sim0.6$} & 25 & 5 & 0 \\ \hline
            {\small $0.6\sim0.7$} & 63 & 10 & 0 \\ \hline
            {\small $0.7\sim0.8$} & 114 & 19 & 0 \\ \hline
            {\small $0.8\sim0.9$} & 415 & 47 & 0\\ \hline
            {\small $0.9\sim1$} & 9,353 & 9,905 & 10,000 \\ \hline
		\end{tabular}
        \vspace{0.1in}
		\subcaption{Full scale}
		\label{Histogram_1_server.png}
	\end{minipage}%
	\hfill
	\begin{minipage}{0.5\textwidth} 
		\centering
		\renewcommand{\arraystretch}{1.5}
		\begin{tabular}{|>{}c|>{}c|>{}c|>{}c|}  
			\hline
			{\tiny Performance ratio} & {\tiny Best static policy} & {\tiny Best cutoff-static policy} & {\tiny Best two-price policy}\\ \hline
			{\small $0.90\sim0.91$} & 93 & 3 & 1\\ \hline
			{\small $0.91\sim0.92$} & 93 & 10 & 0 \\ \hline
            {\small $0.92\sim0.93$} & 114 & 28 & 0 \\ \hline
            {\small $0.93\sim0.94$} & 145 & 47 & 0\\ \hline
            {\small $0.94\sim0.95$} & 211 & 94 & 0 \\ \hline
            {\small $0.95\sim0.96$} & 282 & 197 & 0 \\ \hline
            {\small $0.96\sim0.97$} & 386 & 382 & 1 \\ \hline
            {\small $0.97\sim0.98$} & 604 & 654 & 3 \\ \hline
            {\small $0.98\sim0.99$} & 1,162 & 1,311 & 23\\ \hline
            {\small $0.99\sim1$} & 6,263 & 7,179 & 9,972 \\ \hline
		\end{tabular}
        \vspace{0.1in}
		\subcaption{Zoom in}
		\label{Histogram_1_server_90.png}
	\end{minipage}
	\label{Histogram_1_server}
\end{table}
\endgroup}
\vspace{-0.3in}

From Table \ref{Histogram_1_server.png}, we observe that the best two-price policy has the best performance with all the samples giving performance ratios of at least 0.9. The best cutoff-static policy has 9,905 out of 10,000 sample points with a performance ratio above 0.9 and the best static policy only has 9,353 out of 10,000 sample points with a performance ration above 0.9. This shows that the best cutoff-static policy can also be considered near optimal. However, there are 10 cases where the best cutoff-static policy gives a  performance ratio smaller than 0.1, which demonstrates the potential drawbacks of this heuristic under certain combinations of system parameters.  
This indicates that the best cutoff-static policy can be arbitrarily bad in a relative sense when approximating the optimal dynamic gain. We also include the best static policy as a benchmark to show that both of our heuristics are significant improvements compared to even the best static policy, in which only 9,353 out of 10,000 data points have performance ratios above 0.9. To further demonstrate the advantages of the two-price policy, we zoom into the last row of Table \ref{Histogram_1_server.png} in Table \ref{Histogram_1_server_90.png} to focus on the data where the performance ratio is over 0.9.

From Table \ref{Histogram_1_server_90.png}, we notice that the best two-price policy has 9,972 out of 10,000 data points with performance ratio above 0.99. In comparison, the best cutoff-static policy only has 7,179 out of 10,000 data points with performance ratio above 0.99. The best static policy has the worst performance among all three heuristics with only 6,263 out of 10,000 data points with performance ratio above 0.99. This shows that the best two-price policy is a significant improvement over the best cutoff-static policy, which is in turn considerably better than the best static policy.

\begin{figure}[h]   
    \centering
    \includegraphics[width=\linewidth]{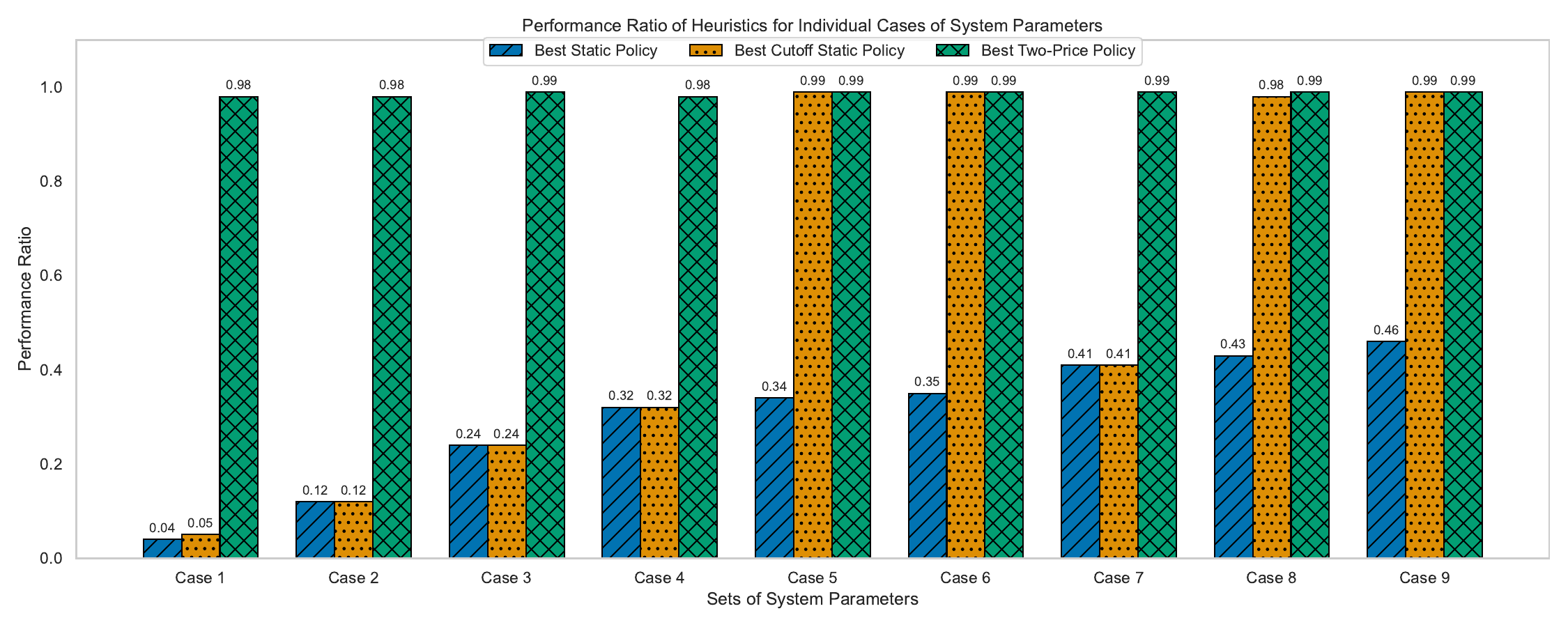}
    \caption{Performance ratio of heuristics for individual cases of system parameters with evaluation distribution $U(20,50)$}
    \label{heuristic_performance_plot}
\end{figure}
\vspace{-0.3in}

{To further demonstrate the improvement of the best two-price policy over the best static policy, we randomly select 9 cases from the dataset used in Table \ref{Histogram_1_server} in which the best static policy performs poorly as a heuristic (with the performance ratio $R_S$ of the best static policy being less than 0.5) and display the performance ratios for the three heuristics in each case in Figure \ref{heuristic_performance_plot}. 
We observe that our proposed heuristics significantly improved the performance compared to the best static policy. In Cases 5, 6, 8, 9, the best cutoff-static policy achieves a performance ratio close to 1, while the best static policy only yields a performance ratio around 0.4. In Cases 1, 2, 3, 4, 7, although the best cutoff-static policy does not provide much improvement over the best static policy, the best two-price policy enhances the performance ratio to close to 1. }
{Additional numerical experiments with an exponential evaluation distribution are provided in Section EC.5. Although the exponential distribution has an unbounded support, which does not align with our assumption of finite support of the evaluation distribution, our numerical results show that our algorithms and heuristics are robust and perform well in this setting.}

We now focus on understanding why the best two-price policy outperforms the cutoff-static policy significantly in some cases. We believe the better performance of the best two-price policy $\pi_T^*$ is due to its ability to mimic the structure of the optimal dynamic policy $\vec{\lambda}^*$, especially when the optimal dynamic policy is not monotone decreasing. Recall the discussion in Section 5.3 that all two-price policies have uni-modal structure. Furthermore, when $\delta_s \geq \delta_q$, the two-price policy $\pi_T^*$ has a monotone decreasing structure. Then, for each set of system parameters, there are four possible situations for the structure of the optimal dynamic policy and the best two-price policy $\pi_T^*$: (1) both $\pi_T^*$ and $\vec{\lambda}^*$ are monotone decreasing; (2) neither $\pi_T^*$ nor $\vec{\lambda}^*$ is monotone decreasing; (3) $\vec{\lambda}^*$ is monotone decreasing but $\pi_T^*$ is not; (4) $\pi_T^*$ is monotone decreasing but $\vec{\lambda}^*$ is not. Table 1 records the frequency out of 10,000 randomly generated data for each of the four possibilities above. Table \ref{dynamic_tp_compare_uni} uses {the randomly generated datasets} that are used in Figures \ref{Histogram_1_server.png} and \ref{Histogram_1_server_90.png}. Table \ref{dynamic_tp_compare_exp} uses {the randomly generated datasets} that are used in Table \ref{Histogram_exp}. 

{\begingroup\renewcommand{\baselinestretch}{0.6}\normalsize
\begin{table}[H] 
	\centering
	\caption{Monotonicity of optimal dynamic policy vs.\ two-price policy}
	\begin{minipage}{0.5\textwidth} 
		\centering
		\renewcommand{\arraystretch}{1.5} 
		\setlength{\tabcolsep}{10pt} 
		\begin{tabular}{|c|c|c|}
			\hline
			{\tiny $F\sim U(20,50)$} & {\tiny$ \vec{\lambda}^*$ Monotone} & {\tiny $\vec{\lambda}^*$ Non Monotone} \\ \hline
			{\tiny $\pi_T^*$ Monotone} & 7762 & 10 \\ \hline
			{\tiny $\pi_T^*$ Non Monotone} & 1 & 2227 \\ \hline
		\end{tabular}
        \vspace{0.1in}
		\subcaption{\small Uniform evaluation distribution}
		\label{dynamic_tp_compare_uni}
	\end{minipage}%
	\hfill
	\begin{minipage}{0.5\textwidth} 
		\centering
		\renewcommand{\arraystretch}{1.5}
		\setlength{\tabcolsep}{10pt}
		\begin{tabular}{|c|c|c|}
			\hline
			{\tiny $F\sim Exp (1/35)$} & {\tiny $\vec{\lambda}^*$ Monotone} & {\tiny $\vec{\lambda}^*$ Non Monotone} \\ \hline
			{\tiny $\pi_T^*$ Monotone} & 7721 & 2 \\ \hline
			{\tiny $\pi_T^*$ Non Monotone} & 0 & 2277 \\ \hline
		\end{tabular}
        \vspace{0.1in}
		\subcaption{\small Exponential evaluation distribution}
		\label{dynamic_tp_compare_exp}
	\end{minipage}

\end{table}
\endgroup}
\vspace{-0.3in}

From Table \ref{dynamic_tp_compare_uni}, we observe that the structure of the best two-price policy $\pi_T^*$ aligns with the structure of the optimal dynamic policy $\vec{\lambda}^*$ in most cases. There is only one case where $\vec{\lambda}^*$ is monotone decreasing but $\pi_T^*$ is not monotone ($\delta_s < \delta_q$) and there are 10 cases where $\vec{\lambda}^*$ policy is not monotone decreasing but $\pi_T^*$ is (i.e., $\delta_s\geq \delta_q$). In total, there are 11 out of 10,000 data points in which the structure of $\pi_T^*$ is different from that of $\vec{\lambda}^*$. Table \ref{dynamic_tp_compare_exp} gives a similar result for an exponential evaluation distribution. In particular, $\vec{\lambda}^*$ and $\pi_T^*$ have structure differences in only 2 out of 10,000 cases. 
In general, the numerical experiments suggest the best two-price policy can capture the monotone and non-monotone structure of the optimal dynamic policy. In contrast, the best cutoff-static policy $\pi_C^*$ always has a monotone decreasing structure.

As we mentioned earlier, the best cutoff-static policy sometimes yields arbitrarily bad relative performance, namely it yields zero gain while the gain of the optimal dynamic policy is positive. This motivated the development of the best two-price policy. 
To further illustrate the strength of the best two-price policy, we give two randomly selected examples where the gain of the best cutoff-static policy is non-zero and the best two-price policy significantly outperforms the best cutoff-static policy.

\begin{figure}[h]
	\centering
	\caption{Comparison of arrival rates for three pricing policies at each state}
	\begin{minipage}{0.45\textwidth}
		\centering
		\includegraphics[width=\columnwidth]{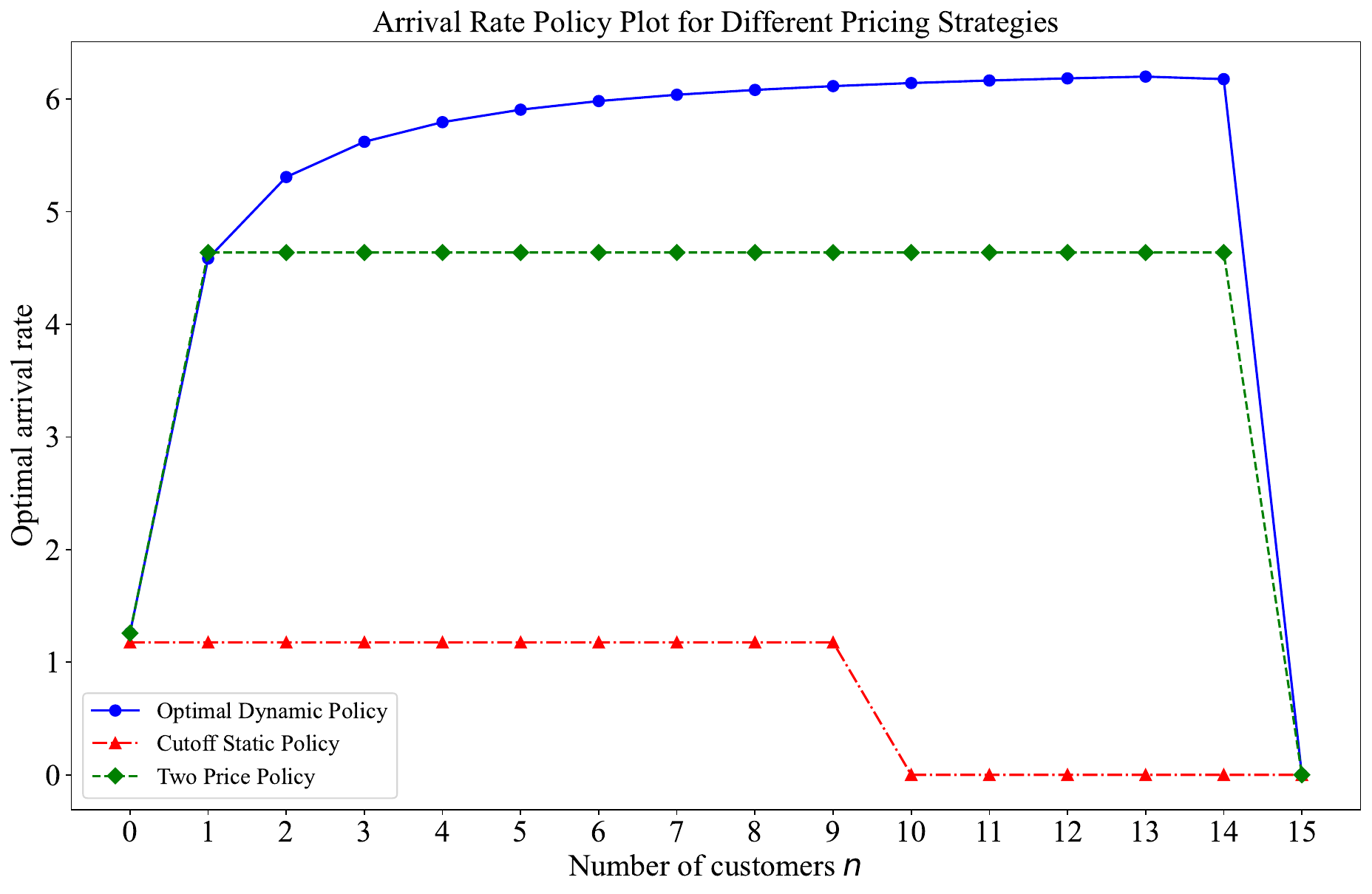}
		\subcaption{Single-server system}
        \label{arrival_rate_comparison_1_server}
	\end{minipage}
	\begin{minipage}{0.45\textwidth}
		\centering
		\includegraphics[width=\columnwidth]{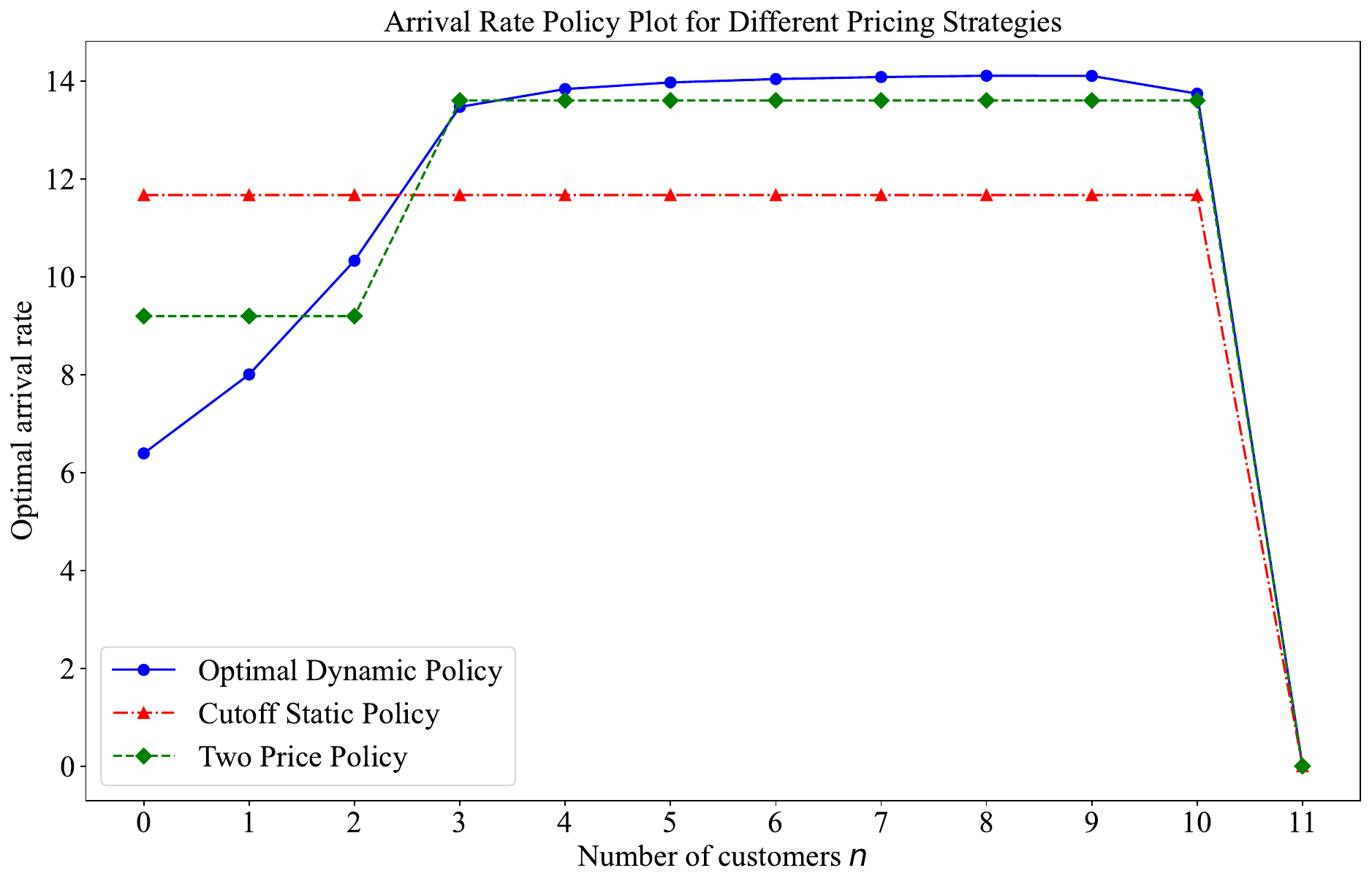}
		\subcaption{Multi-server system}
        \label{arrival_rate_comparison_multi_server}
	\end{minipage}
\end{figure}

Figure \ref{arrival_rate_comparison_1_server} shows the optimal dynamic policy, the best cutoff-static policy, and the best two-price policy when the system parameters are $N= 15,\,\Lambda = 16.517,\,\mu=1.55,\,m=1,\,\theta_s=19.32,\,\theta_q=36.24,\,c_s=30.767,\,c_q=6.85,\,c_h=15.117$, and the evaluation distribution follows $U(20,50)$. Clearly, the best cutoff-static policy only captures the arrival rate at state 0 (with an idling server) and state 15 (when the system is full) but misses the rest of the states (when there is congestion in the system). In comparison, the best two-price policy captures the arrival rate at state 0,1,15 and approximates it better for the other states when there is waiting in the system. Given that the best two-price policy resembles the optimal dynamic policy much more closely than the best cutoff-static policy, we expect the gain of the best two-price policy to be significantly higher than that of the best cutoff-static policy. {Indeed, in this case, the performance ratio of the best two-price policy is 0.998 while that of the best cutoff-static policy is 0.593, which is a significant increase.}

\vspace{-0.4in}

We next look at an example with multiple servers. Figure \ref{arrival_rate_comparison_multi_server} shows the structures of optimal dynamic policy, best cutoff-static policy, and best two-price policy when the system parameters are $N= 11,\,\Lambda = 50,\,\mu=1,\,m=3,\,\theta_s=2,\,\theta_q=30,\,c_s=45,\,c_q=20,\,c_h=45$ and the evaluation distribution follows $U(20,50)$. The best cutoff-static policy approximates the structure of the optimal dynamic policy better compared to the case in Figure \ref{arrival_rate_comparison_1_server}. But the best two-price policy mimics the structure of the optimal dynamic policy even better than the best cutoff-static policy. {Under this set of system parameters, the performance ratio of the best cutoff-static policy is 0.864 and that of the best two-price policy is 0.971.} This indicates that the best two-price policy also approximates the structure of optimal dynamic policy better than the best cutoff-static policy and significantly improves the performance of the best cutoff-static policy for a multi-server system.

\section{Conclusion}\label{sec.conclustion}

This paper investigates the optimal pricing policy in a multiple-server queueing system where customers can abandon during service or while waiting. We formulate the optimal pricing problem as a Markov decision process. By assuming the customers follow the same evaluation distribution independent of each other, we show that finding the optimal pricing policy is equivalent to finding the optimal arrival rate policy in all states. 

We show that the presence of customer abandonments complicates the pricing problem in that the structure of the optimal arrival rate policy is not always monotone. We prove that the optimal arrival rate policy always follows a uni-modal structure. 
{Namely, the corresponding quoted price may first decrease and then increase as the number of customers in the system grows.}
Also, under reasonable assumptions on the costs and abandonment rates, we prove that the optimal arrival rate policy is monotone decreasing, {which implies that the optimal quoted price increases as the system gets more congested. Based on these results, we develop a policy iteration algorithm that exploits these structures to efficiently solve for the optimal arrival rate policy.}

Furthermore, we introduce two heuristics to simplify the optimal dynamic arrival rate policy. The cutoff-static policy charges the incoming customers the same price until the number of customers in the system reaches a certain threshold.  
The two-price policy charges one price if there is an idling server in the system and another price if all servers are busy until a certain threshold. {We describe how to identify the best versions of these heuristics and derive theoretical results that provide upper bounds on the profit loss for both the best cutoff static policy and the best two-price policy.}

Finally, we conduct numerical experiments to demonstrate the performance of our heuristics under different system parameters and evaluation distributions. The results show that the best cutoff-static policy is near optimal in general but it cannot approximate the optimal dynamic policy under some system parameters. In comparison, the best two-price policy is always near optimal.  

 {Future research may extend our framework to consider other reimbursement mechanisms when customers abandon the system. Incorporating reimbursement that depends on the price paid upon entry can be challenging because it requires the system to keep track of the quoted price for each customer in the system, which can significantly complicate the state space for the Markov decision process. However, space reduction in some asymptotic regimes, similar to \citet{Kim&Randhawa2017}, \citet{Varma&etal2023}, and \citet{wenxin2025},  may provide valuable insights. }

\newpage

\ACKNOWLEDGMENT{This work was supported by NSF under grant 2127778. The second author was also supported by NSF under grant 2348409.}

\bibliographystyle{my-plainnat}
\bibliography{references}

\ECSwitch
\vspace{-0.15in}
\ECHead{Electronic Companion to ``Pricing in Queues with Abandonments: Optimal Policies and Practical Heuristics''}

{This Electronic Companion provides the detailed proofs of the theoretical results in Sections \ref{sec.monotonicity} and \ref{sec.Heuristics}  in Sections \ref{ec.proof.section.4} and \ref{ec.proof.section.5}, respectively. Moreover, Sections EC.3 and EC.4 contain additional numerical experiments for Sections 4.2 and 6, respectively.}

\section{Proofs of Theorem \ref{Theorem.uni_modal} and \ref{Theorem.monotonic}}  \label{ec.proof.section.4}
\noindent In this section, we prove the major theorems in Section 4.
We start by defining the difference operator $\Delta h:\BR\to \BR$ for the bias function $h(n)\in\BR$ as follows:
\[
\Delta h(n) = h(n+1) - h(n) \text{ for } n=0,\ldots, N-1. 
\]
Similarly, we define the double difference of the bias function $\Delta^2 h(n)$ as 
\begin{align*}
\Delta^2 h(n) &= \Delta(\Delta h(n)) = \Delta h(n+1) - \Delta h(n) \addrefnum \label{Delta2h} \\
& =[h(n+2) - h(n+1)] - [h(n+1) - h(n)] 
\end{align*} for $n=0,\ldots,N-2$. 
When $\Delta h(n)\leq 0$ for all $n=0,\ldots,N-1$, we say the bias function $h(n)$ is \textit{decreasing}; when $\Delta^2 h(n) \leq 0$ for all $n=0,\ldots, N-2$, we say $h(n)$ is \textit{concave}.

Using the definition of $\Delta h(n)$, we define 
\[
\Delta w(n) = 
\begin{cases}
	(\mu+\theta_s) \Delta h(n) + c_h + c_s\theta_s & \text{ for } 0\leq n <m,\\
	\theta_q \Delta h(n) + c_h + c_q\theta_q & \text{ for } m\leq n < N.
\end{cases}
\]
Furthermore, for $n=0,\ldots,N-1$, consider the difference between the reward functions of states $n$ and $n+1$ when the arrival rate for both states is $\lambda\in[0,\Lambda]$:
\begin{align*}
r_{n}(\lambda) - r_{n+1}(\lambda) =& \lambda \bar{F}^{-1}(\lambda/\Lambda) - nc_h - \min\{n,m\}c_s\theta_s -\max\{n-m,0\}c_q\theta_q \\
&- (\lambda \bar{F}^{-1}(\lambda/\Lambda) - (n+1)c_h - \min\{n+1,m\}c_s\theta_s -\max\{n+1-m,0\}c_q\theta_q)\\
= &  c_h + c_s\theta_s \mathbf{1}_{n< m} + c_q\theta_q \mathbf{1}_{n\geq m},  \addrefnum \label{reve.diff}
\end{align*}
which does not depend on $\lambda$. Here $\mathbf{1}_{n<m}$ and $\mathbf{1}_{n\geq m}$ are indicator functions.

{The rest of the section is organized as follows. We first state the lemmas used to prove Theorems \ref{Theorem.uni_modal} and \ref{Theorem.monotonic} and prove Theorems \ref{Theorem.uni_modal} and \ref{Theorem.monotonic} in Section \ref{ec.proof.tm41&tm42}. Then, we prove the lemmas stated in Section \ref{ec.proof.tm41&tm42} in Section \ref{ec.proof.lemma.in.1.1}.}

\subsection{Lemmas for Proofs of Theorems \ref{Theorem.uni_modal} and \ref{Theorem.monotonic}} \label{ec.proof.tm41&tm42}

{In this section, we prove Theorems \ref{Theorem.uni_modal} and \ref{Theorem.monotonic} by stating several related lemmas. We first establish the relationship between the optimal arrival rate $\lambda_n^*$ at each state $n$ and the double difference of the bias function $\Delta^2h(n)$. }

\begin{lemma} \label{lemma.dece}
For a bias function $h(n)$ satisfying the Bellman equation in (\ref{Bellman.equation}) and optimal arrival rate policy $\vec{\lambda}^*$ defined in (\ref{def.opt.rate}), we have
\begin{enumerate}
\item $\lambda_0^* \Delta^2 h(0)\leq \Delta w(0) \leq \lambda_1^* \Delta^2 h(0)$.
\item $\lambda_n^* \Delta^2 h(n)\leq \Delta w(n) +\gamma_n \Delta^2 h(n-1) \leq \lambda_{n+1}^* \Delta^2 h(n)$ for $1\leq n\leq N-2$.
\item $\Delta w(N-1) +\gamma_{N-1}\Delta^2 h(N-2)\leq 0.$
\end{enumerate}
Here $\gamma_n$ is the rate for state $n$ defined in Section 3.
\end{lemma}

{
From Items 1 and 2 in Lemma \ref{lemma.dece}, we have that $\lambda_n^*\Delta^2h(n) \leq \lambda_{n+1}^*\Delta^2h(n)$ for all $0\leq n\leq N-2$. Then, the sign of the double difference will determine whether the optimal arrival rate is increasing or decreasing in that state when $\lambda_n^*\not =0$. Namely, when $\Delta^2 h(n)<0$, we have $\lambda_n^*\geq \lambda_{n+1}^*$, implying that the optimal arrival rate policy is decreasing in that state. Item 3 in Lemma \ref{lemma.dece} is a critical boundary condition that will be used while proving that the optimal arrival rate policy is decreasing under certain conditions. We defer the proof of this lemma to the next section.
}

\begin{remark} \label{remark_all_posi_rate}
In the remaining proofs of this section, we assume $\lambda_n^*>0$ for all $0\leq n\leq N-1$. We can impose this assumption because it will not change the results in our theorems. If there exists $\lambda_n^*=0$ for some $0\leq n<N$, then $P_i(\vec{\lambda}^*)=0$ for all $n< i<N$ by its definition in Section 3. Thus, any arrival rate policy $\vec{\lambda}$ with $\lambda_n=\lambda_n^*$ for $0\leq i\leq n$ and $\lambda_i$ being any value in $[0,\Lambda]$ for $n<i<N$ also gives the optimal long-run average gain, which means that $\vec{\lambda}$ satisfies the Bellman equation in (\ref{Bellman.equation}). Therefore, by the definition of the optimal arrival rate in (\ref{def.opt.rate}), we have $\lambda^*_i=0$ for all $n\leq i<N$, which is monotone decreasing for $i\geq n$. Thus, our results will follow by replacing $N$ in the proofs by $\min\{n\in\{0,\ldots,N\}:\lambda^*_n=0\}$, ensuring the optimal arrival rates at all positive recurrent states to be positive.
\end{remark}

{The next lemma shows that when $\CC_s\leq \CC_q$, the double differences of the bias function are less than or equal to 0.}

\begin{lemma} \label{lemma.mono}
If $\lambda_n^*>0$ for all $n=0,\ldots,N-1$, then $\CC_s=\frac{c_h+c_s\theta_s}{\mu+\theta_s}\leq \frac{c_h}{\theta_q}+c_q=\CC_q$ implies that $\Delta^2 h(n)\leq 0$ for $0\leq n\leq N-2$.
\end{lemma}

We defer the proof of this lemma to the next section and we are ready to prove Theorem \ref{Theorem.monotonic}.

\begin{proof}{Proof of Theorem \ref{Theorem.monotonic}}
From Remark \ref{remark_all_posi_rate}, we can assume $\lambda_n^* >0$ for all $0\leq n<N$. Then, 
from Lemma \ref{lemma.mono}, we have $\Delta^2 h(n)\leq 0$ for $0\leq n\leq N-2$. 
Lemma \ref{lemma.dece} gives us 
\[
\lambda_n^* \Delta^2 h(n)  \leq \lambda_{n+1}^* \Delta^2 h(n) \text{ for }0\leq n\leq N-2.  \addrefnum \label{rate.leq}
\]
If $\Delta^2 h(n)<0$, we have $\lambda_n^*\geq \lambda_{n+1}^*$ for $0\leq n\leq N-2$. 
If $\Delta^2 h(n)=0$ for some $0\leq n\leq N-2$, then (\ref{def.opt.rate}) yields 
\begin{align*}
r_{n}(\lambda_{n}^*)+\sum_{j=0}^N \mathbf{P}(j|n,\lambda_n^*) h(j) &\geq r_n(\lambda_{n+1}^*)+\sum_{j=0}^N \mathbf{P}(j|n,\lambda_{n+1}^*) h(j), \addrefnum \label{sub.n.with.n+1} \\
r_{n+1}(\lambda_{n+1}^*)+\sum_{j=0}^N \mathbf{P}(j|n+1,\lambda_{n+1}^*) h(j) &\geq r_{n+1}(\lambda_{n}^*)+\sum_{j=0}^N \mathbf{P}(j|n+1,\lambda_n^*) h(j). \addrefnum \label{sub.n+1.with.n}
\end{align*}
Summing up (\ref{sub.n.with.n+1}) and (\ref{sub.n+1.with.n}) and rearranging the terms, we have 
\begin{align*}
	& r_{n+1}(\lambda_{n+1}^*)+\sum_{j=0}^N \mathbf{P}(j|n+1,\lambda_{n+1}^*) h(j) - r_n(\lambda_{n+1}^*)-\sum_{j=0}^N \mathbf{P}(j|n,\lambda_{n+1}^*) h(j) \\
	\geq&  r_{n+1}(\lambda_{n}^*)+\sum_{j=0}^N \mathbf{P}(j|n+1,\lambda_n^*) h(j) - r_{n}(\lambda_{n}^*) - \sum_{j=0}^N \mathbf{P}(j|n,\lambda_n^*) h(j). \addrefnum \label{sum_of_n&n+1}
\end{align*}
From (\ref{bias.diff.0.leq}) and (\ref{bias.diff.0.geq}), (\ref{sum_of_n&n+1}) can be reduced to
\begin{align*}
\lambda_{1}^* \Delta^2 h(0) + \Delta h(0) -\Delta w(0) \geq \lambda_{0}^* \Delta^2 h(0) + \Delta h(0) -\Delta w(0).
\end{align*}
Similarly, for $1\leq n\leq N-2$, from (\ref{bias.diff.<m.leq})-(\ref{bias.diff.>m.geq}), (\ref{sum_of_n&n+1}) can be reduced to
\[
\lambda_{n+1}^* \Delta^2 h(n) + \Delta h(n) -\Delta w(n) -\gamma_n \Delta^2 h(n-1) \geq \lambda_{n}^* \Delta^2 h(n) + \Delta h(n) -\Delta w(n) -\gamma_n \Delta^2 h(n-1).
\]
Thus, we conclude that if $\Delta^2 h(n)=0$ for some $0\leq n\leq N-2$, then (\ref{sum_of_n&n+1}) achieves equality. Since the sum of (\ref{sub.n.with.n+1}) and (\ref{sub.n+1.with.n}) achieves equality, then it follows that the equality also holds in (\ref{sub.n.with.n+1}) and (\ref{sub.n+1.with.n}), 
which implies that $\lambda_n^*$ and
$\lambda_{n+1}^*$ are optimal arrival rates for both states $n$ and $n+1$. It now follows from (\ref{def.opt.rate}) that $\lambda_n^* \leq \lambda_{n+1}^*$ and similarly that
\[
\lambda_{n+1}^* = \inf \Bigg( \arg \max_{\lambda\in A_{n+1}}\Big\{ r_{n+1}(\lambda)+\sum_{j=0}^N \mathbf{P}(j|n+1,\lambda) h(j)  \Big\} \Bigg)\leq \lambda_n^*.
\]
We have shown that $\lambda_n^*\geq \lambda_{n+1}^*$ for $0\leq n\leq N-2$, both when $\Delta^2 h(n)<0$ and $\Delta^2 h(n)=0$. The result now follows from Lemma \ref{lemma.mono} and the fact that $\lambda_N^*=0$.
\end{proof}

Moving on to prove Theorem \ref{Theorem.uni_modal}, we first define $\tilde{n}$ as:
\[
\tilde{n} = 
\begin{cases}
N-1, \text{ if } \Delta^2 h(n) > 0 \text{ for all } 0\leq n\leq N-2, \\
\min \{ n: \Delta^2 h(n) \leq 0 \text{ for } 0\leq n\leq N-2 \}, \text{ otherwise}.
\end{cases} \addrefnum \label{def.tilde_n}
\]
Thus, $\tilde{n}$ is the first state where $\Delta^2 h(n)$ is non-positive. We will show that the optimal arrival rate policy $\vec{\lambda}^*$ satisfies Definition \ref{def.uni.modal} with $\hat{n} = \tilde{n}$, indicating the optimal arrival rate policy is uni-modal. {The following lemma establishes that all double differences after state $\tilde{n}$ are non-positive.}

\begin{lemma} \label{lemma.uni_modal}
If $\lambda_n^*>0$ for all $n=0,\ldots,N-1$ and $\tilde{n}\leq N-2$, then $\Delta^2 h(n)\leq 0$ for all $\tilde{n} \leq n\leq N-2$. Also, either $\tilde{n}\geq m$ or $\tilde{n}=0$.
\end{lemma}

\begin{remark} \label{remark.tilde.n=N-1}
In Lemma \ref{lemma.uni_modal}, we show that when $\tilde{n}\leq N-2$, we have either $\tilde{n}\geq m$ or $\tilde{n}=0$. Now we consider the case when $\tilde{n}=N-1$. If $N\geq m+1$, then $\tilde{n}=N-1\geq m$ follows immediately. If $N=m$, we show that $\tilde{n}$ is never equal to $N-1$. Suppose $\tilde{n} = N-1 = m-1$. From the definition of $\tilde{n}$, we have $\Delta^2 h(m-2)>0$ and Lemma \ref{lemma.sign.n0} yields $\Delta w(m-2)>0$. Then, from (\ref{Delta.w.induction.leq.m.step}), we have $\Delta w(m-1) = \Delta w(m-2) + (\mu+\theta_s) \Delta^2 h(m-2)>0$, which implies $\Delta w(m-1)+ \gamma_{m-1} \Delta^2 h(m-2)>0$, contradicting  Part 3 of Lemma \ref{lemma.dece} for $N=m$. Therefore, when $N=m$, $\tilde{n}\not=N-1$ always, which implies $\tilde{n}\leq N-2$. Then, from Lemma \ref{lemma.uni_modal}, we have $\tilde{n}=0$ when $N=m$.
\end{remark}

{We defer the proof of Lemma \ref{lemma.uni_modal} to the next section and prove Theorem \ref{Theorem.uni_modal}.}

\begin{proof}{Proof of Theorem \ref{Theorem.uni_modal}.}
From Remark \ref{remark_all_posi_rate}, we can assume $\lambda_n^* >0$ for all $0\leq n<N$. 
Now define $\tilde{n}$ as shown in (\ref{def.tilde_n}). In order to prove the theorem, we intend to show that $\lambda_n^*\leq \lambda_{n+1}^*$ for $n<\tilde{n}$ and $\lambda_n^*\geq \lambda_{n+1}^*$ for $n\geq\tilde{n}$. 
If $\tilde{n}=N-1$, then $\Delta^2 h(n)>0$ for all $0\leq n\leq N-2$. From Lemma \ref{lemma.dece} and $\Delta^2 h(n)>0$ for $0\leq n\leq N-2$, then
\[
\lambda_n^* \Delta^2 h(n) \leq \lambda_{n+1}^* \Delta^2 h(n) \Rightarrow \lambda_{n}^* \leq \lambda_{n+1}^* \text{ for } 0\leq n\leq N-2.
\]
Given $\lambda_N^*=0$, we have $\lambda_{N-1}^* \geq 0 = \lambda_N^*$, and hence $\vec{\lambda}^*$ has a uni-modal structure with $\hat{n} = \tilde{n}$. 
If $\tilde{n}\leq N-2$, by Lemma \ref{lemma.uni_modal}, we have $\Delta^2 h(n)>0$ for $n< \tilde{n}$ and $\Delta^2 h(n)\leq 0$ for $n\geq \tilde{n}$.
When $n<\tilde{n}$, the above argument yields $\lambda_n^*\leq \lambda_{n+1}^*$ for $n<\tilde{n}$. 
Then, when $n\geq \tilde{n}$, from the proof of Theorem \ref{Theorem.monotonic}, we have $\lambda_n^*\geq \lambda_{n+1}^*$. \Halmos
\end{proof}

\begin{remark}
It follows from Lemma \ref{lemma.uni_modal} and Remark \ref{remark.tilde.n=N-1} and the proof of Theorem \ref{Theorem.uni_modal} that it suffices to consider policies satisfying Definition \ref{def.uni.modal} with $\hat{n}=0$ (monotone decreasing) or $\hat{n}\geq m$ (uni-modal but not monotone decreasing). This is consistent with the intuition provided in the last paragraph of Section 4.1.
\end{remark}

\subsection{Proof of Lemmas in Section EC 1.1}  \label{ec.proof.lemma.in.1.1}

We first prove Lemma \ref{lemma.dece}. We derive this lemma by using the Bellman equations and moving substituting the optimal arrival rate for different states.
\begin{proof}{Proof of Lemma \ref{lemma.dece}}
Since $h(n)$ satisfies the Bellman equation (\ref{Bellman.equation}) and $\lambda_n^*$ is the optimal arrival rate at state $n$ where $0\leq n\leq N-1$, we have
\begin{align*}
	g^*+h(n) =& \max_{\lambda\in A_n}\Bigg\{ r_n(\lambda)+\sum_{j=0}^N \mathbf{P}(j|n,\lambda) h(j)  \Bigg\} 
	\geq r_n(\lambda_{n+1}^*)+\sum_{j=0}^N \mathbf{P}(j|n,\lambda_{n+1}^*) h(j).  \addrefnum \label{bias_diff_geq_n}
\end{align*}
Similarly, we have
\begin{align*}
	g^* + h(n+1) =& \max_{\lambda\in A_{n+1}}\Bigg\{ r_{n+1}(\lambda)+\sum_{j=0}^N \mathbf{P}(j|n+1,\lambda) h(j)  \Bigg\} 
	 \geq  r_{n+1}(\lambda_{n}^*)+\sum_{j=0}^N \mathbf{P}(j|n+1,\lambda_n^*) h(j).
	\addrefnum \label{bias_diff_geq_n+1}
\end{align*}
When $n=0$, from (\ref{reve.diff}), (\ref{bias_diff_geq_n}), and $\gamma_1=\mu+\theta_s$:
\begin{align*}
	\Delta h(0) =& h(1) - h(0) = g^* + h(1) - (g^*+ h(0)) \\
	\leq & r_{1}(\lambda_{1}^*) + \sum_{j=0}^N \mathbf{P}(j|1,\lambda_{1}^*) h(j) - r_0(\lambda_{1}^*)-\sum_{j=0}^N \mathbf{P}(j|0,\lambda_{1}^*) h(j)\\
	=& r_{1}(\lambda_{1}^*) + \lambda_{1}^* h(2) + \gamma_{1} h(0) + (1-\lambda_{1}^*-\gamma_{1}) h(1)  - [r_{0}(\lambda_{1}^*) + \lambda_{1}^* h(1) + (1-\lambda_{1}^*) h(0)]\\	
	=& r_{1}(\lambda_{1}^*) - r_0(\lambda_{1}^*) + \lambda_{1}^* (h(2) - h(1)) + (1-\lambda_{1}^*-(\mu+\theta_s))(h(1)-h(0))\\
	=& -(c_h+c_s\theta_s) + \lambda_{1}^* [h(2) - h(1) - (h(1)-h(0))] + (1-(\mu+\theta_s)) (h(1) - h(0))  \\
	=& \lambda_{1}^* \Delta^2 h(0) + \Delta h(0) -(c_h+c_s\theta_s) -(\mu+\theta_s) \Delta h(0) \\
	=& \lambda_{1}^* \Delta^2 h(0) + \Delta h(0) -\Delta w(0). \addrefnum \label{bias.diff.0.leq}
\end{align*}
Hence, $\Delta h(0)\leq \lambda_{1}^* \Delta^2 h(0) + \Delta h(0) -\Delta w(0)$ and we immediately have
$
\Delta w(0) \leq \lambda_1^* \Delta^2 h(0).
$
Similarly, from (\ref{reve.diff}) and (\ref{bias_diff_geq_n+1}), we have
\begin{align*}
	\Delta h(0) =& h(1) - h(0) = g^* + h(1) - (g^*+ h(0)) \\
	\geq & r_{1}(\lambda_{0}^*) + \sum_{j=0}^N \mathbf{P}(j|1,\lambda_{0}^*) h(j) - r_0(\lambda_{0}^*)-\sum_{j=0}^N \mathbf{P}(j|0,\lambda_{0}^*) h(j)\\
	=& r_{1}(\lambda_{0}^*) + \lambda_{0}^* h(2) + \gamma_{1} h(0) + (1-\lambda_{0}^*-\gamma_{1}) h(1)  - [r_{0}(\lambda_{0}^*) + \lambda_{0}^* h(1) + (1-\lambda_{0}^*) h(0)]\\
	=& r_{1}(\lambda_{0}^*) - r_0(\lambda_{0}^*) + \lambda_{0}^* (h(2) - h(1)) + (1-\lambda_{0}^*-(\mu+\theta_s))(h(1)-h(0))\\
	=& -(c_h+c_s\theta_s) + \lambda_{0}^* [h(2) - h(1) - (h(1)-h(0))] + (1-(\mu+\theta_s)) (h(1) - h(0))  \\
	=& \lambda_{0}^* \Delta^2 h(0) + \Delta h(0) -(c_h+c_s\theta_s) -(\mu+\theta_s) \Delta h(0) \\
	=& \lambda_{0}^* \Delta^2 h(0) + \Delta h(0) -\Delta w(0), \addrefnum \label{bias.diff.0.geq}
\end{align*}
and consequently 
$
\Delta w(0) \geq \lambda_0^* \Delta^2 h(0).
$ Thus, we have $\lambda_0^* \Delta^2 h(0)\leq \Delta w(0) \leq \lambda_1^* \Delta^2 h(0),$ proving the first part of the Lemma.

When $n\leq N-2$ and $0<n\leq m-1$, from (\ref{reve.diff}), (\ref{bias_diff_geq_n}) and $\gamma_{n+1} =\gamma_n+(\mu+\theta_s)$, we have:
\begin{align*}
\Delta h(n) =& h(n+1) - h(n) = g^* + h(n+1) - (g^*+ h(n)) \\
\leq & r_{n+1}(\lambda_{n+1}^*) + \sum_{j=0}^N \mathbf{P}(j|n+1,\lambda_{n+1}^*) h(j) - r_n(\lambda_{n+1}^*)-\sum_{j=0}^N \mathbf{P}(j|n,\lambda_{n+1}^*) h(j)\\
=& r_{n+1}(\lambda_{n+1}^*) - r_n(\lambda_{n+1}^*) + \lambda_{n+1}^* (h(n+2) - h(n+1))\\
&+ \gamma_n (h(n) - h(n-1)) + (1-\lambda_{n+1}^*-\gamma_n-(\mu+\theta_s))(h(n+1)-h(n))\\
=& -(c_h+c_s\theta_s) + \lambda_{n+1}^* [h(n+2) - h(n+1) - (h(n+1)-h(n))]\\
&+ (1-(\mu+\theta_s)) (h(n+1) - h(n)) - \gamma_n [h(n+1)-h(n)-(h(n) - h(n-1))]  \\
=& \lambda_{n+1}^* \Delta^2 h(n) + \Delta h(n) -(c_h+c_s\theta_s) -(\mu+\theta_s) \Delta h(n)   -\gamma_n \Delta^2 h(n-1)\\
=& \lambda_{n+1}^* \Delta^2 h(n) + \Delta h(n) -\Delta w(n) -\gamma_n \Delta^2 h(n-1). \addrefnum \label{bias.diff.<m.leq}
\end{align*}
Hence, $\Delta h(n)\leq \lambda_{n+1}^* \Delta^2 h(n) + \Delta h(n) -\Delta w(n) -\gamma_n \Delta^2 h(n-1)$
and we immediately have
\begin{align*}
\Delta w(n) + \gamma_n\Delta^2 h(n-1) \leq \lambda_{n+1}^* \Delta^2 h(n).
\end{align*}
Similarly, from (\ref{reve.diff}) and (\ref{bias_diff_geq_n+1}), we have
\begin{align*}
\Delta h(n) =& h(n+1) - h(n) = g^* + h(n+1) - (g^*+ h(n)) \\
 \geq & r_{n+1}(\lambda_{n}^*) + \sum_{j=0}^N \mathbf{P}(j|n+1,\lambda_{n}^*) h(j) - r_n(\lambda_{n}^*)-\sum_{j=0}^N \mathbf{P}(j|n,\lambda_{n}^*) h(j)\\
=& r_{n+1}(\lambda_{n}^*) - r_n(\lambda_{n}^*) + \lambda_{n}^* (h(n+2) - h(n+1))\\
&+ \gamma_n (h(n) - h(n-1)) + (1-\lambda_{n}^*-\gamma_n-(\mu+\theta_s))(h(n+1)-h(n))\\
=& -(c_h+c_s\theta_s) + \lambda_{n}^* \Delta^2 h(n) + (1-(\mu+\theta_s)) \Delta h(n) -\gamma_n \Delta^2 h(n-1)\\
=& \lambda_{n}^* \Delta^2 h(n) + \Delta h(n) -\Delta w(n) -\gamma_n \Delta^2 h(n-1), \addrefnum \label{bias.diff.<m.geq}
\end{align*}
and consequently 
$
\Delta w(n) + \gamma_n\Delta^2 h(n-1) \geq \lambda_{n}^* \Delta^2 h(n)
$. We have proved the second part of the Lemma for $n\leq N-2$ and $0<n\leq m-1$.

When $m\leq n \leq  N-2$, so that $\gamma_{n+1}=\gamma_n+\theta_q$, from (\ref{reve.diff}) and  (\ref{bias_diff_geq_n}), we have
\begin{align*}
\Delta h(n) =& h(n+1) - h(n) = g^* + h(n+1) - (g^*+ h(n)) \\
\leq & r_{n+1}(\lambda_{n+1}^*) + \sum_{j=0}^N \mathbf{P}(j|n+1,\lambda_{n+1}^*) h(j) - r_n(\lambda_{n+1}^*)-\sum_{j=0}^N \mathbf{P}(j|n,\lambda_{n+1}^*) h(j)\\
=& r_{n+1}(\lambda_{n+1}^*) - r_n(\lambda_{n+1}^*) + \lambda_{n+1}^* (h(n+2) - h(n+1))\\
&+ \gamma_n (h(n) - h(n-1)) + (1-\lambda_{n+1}^*-\gamma_n-\theta_q)(h(n+1)-h(n))\\
=& -(c_h+c_q\theta_q) + \lambda_{n+1}^* \Delta^2 h(n) + (1-\theta_q) \Delta h(n) -\gamma_n \Delta^2 h(n-1)\\
=& \lambda_{n+1}^* \Delta^2 h(n) + \Delta h(n) -\Delta w(n) -\gamma_n \Delta^2 h(n-1), \addrefnum \label{bias.diff.>m.leq}
\end{align*}
which implies that $\Delta w(n) + \gamma_n\Delta^2 h(n-1) \leq \lambda_{n+1}^* \Delta^2 h(n).$ 
Similarly, from (\ref{reve.diff}) and (\ref{bias_diff_geq_n+1}):
\begin{align*}
	\Delta h(n) =& h(n+1) - h(n) = g^* + h(n+1) - (g^*+ h(n)) \\
	\geq & r_{n+1}(\lambda_{n}^*) + \sum_{j=0}^N \mathbf{P}(j|n+1,\lambda_{n}^*) h(j) - r_n(\lambda_{n}^*)-\sum_{j=0}^N \mathbf{P}(j|n,\lambda_{n}^*) h(j)\\
	=& r_{n+1}(\lambda_{n}^*) - r_n(\lambda_{n}^*) + \lambda_{n}^* (h(n+2) - h(n+1))\\
	&+ \gamma_n (h(n) - h(n-1)) + (1-\lambda_{n}^*-\gamma_n-\theta_q)(h(n+1)-h(n))\\
	=& -(c_h+c_q\theta_q) + \lambda_{n}^* \Delta^2 h(n) + (1-\theta_q) \Delta h(n) -\gamma_n \Delta^2 h(n-1)\\
	=& \lambda_{n}^* \Delta^2 h(n) + \Delta h(n) -\Delta w(n) -\gamma_n \Delta^2 h(n-1), \addrefnum \label{bias.diff.>m.geq}
\end{align*}
and consequently $
\Delta w(n) + \gamma_n\Delta^2 h(n-1) \geq \lambda_{n}^* \Delta^2 h(n)
$, proving the second part of the Lemma.

When $n=N-1$, from the definition of $\Delta w(n)$, we have 
\[
\Delta w(N-1) = [(\mu+\theta_s) \mathbf{1}_{N= m} + \theta_q \mathbf{1}_{N> m}] \Delta h(N-1)  + [c_h + c_s\theta_s \mathbf{1}_{N= m} + c_q\theta_q \mathbf{1}_{N> m}].
\]
Then $\lambda_N^*=0 \in A_{N-1}$, $\gamma_N = \gamma_{N-1} + (\mu+\theta_s) \mathbf{1}_{N= m} + \theta_q \mathbf{1}_{N> m}$, and (\ref{bias_diff_geq_n}) yield 
\begin{align*}
	\Delta h(N-1) =& h(N) - h(N-1) = g^* + h(N) - (g^*+ h(N-1)) \\
	\leq & r_{N}(\lambda_{N}^*) + \sum_{j=0}^N \mathbf{P}(j|N,\lambda_{N}^*) h(j) - r_{N-1}(\lambda_N^*)-\sum_{j=0}^N \mathbf{P}(j|N-1,\lambda_N^*) h(j)\\
	= & -  [c_h + c_s\theta_s \mathbf{1}_{N= m} + c_q\theta_q \mathbf{1}_{N> m}] + \gamma_N h(N-1) + (1-\gamma_N) h(N) \\
	& - \gamma_{N-1} h(N-2) - (1 - \gamma_{N-1}) h(N-1)  \\
	= & -  [c_h + c_s\theta_s \mathbf{1}_{N= m} + c_q\theta_q \mathbf{1}_{N> m}] + h(N) - h(N-1)\\
	& - (\gamma_N- \gamma_{N-1})[h(N)-h(N-1)] \\
	& - \gamma_{N-1}[(h(N)-h(N-1))-(h(N-1)-h(N-2))]\\
	= & \Delta h(N-1) -\gamma_{N-1} \Delta^2 h(N-2) -  [c_h + c_s\theta_s \mathbf{1}_{N= m} + c_q\theta_q \mathbf{1}_{N> m} ]\\
	& - [(\mu+\theta_s) \mathbf{1}_{N= m} + \theta_q \mathbf{1}_{N> m}]\Delta h(N-1)  \\
	= & \Delta h(N-1) - \Delta w(N-1) -\gamma_{N-1}\Delta^2 h(N-2). 
\end{align*}
Hence, $\Delta h(N-1)\leq  \Delta h(N-1) - \Delta w(N-1) -\gamma_{N-1}\Delta^2 h(N-2)$
and we immediately have
$\Delta w(N-1) +\gamma_{N-1}\Delta^2 h(N-2)\leq 0$, which completes the proof. \Halmos
\end{proof}

{We then derive the required lemmas to prove Lemmas \ref{lemma.mono} and \ref{lemma.uni_modal}. To simplify the notation, recall the sign function }
\[
\mathrm{sgn}(x) = \begin{cases}
	1 &\text{ when } x>0,\\
	0 &\text{ when } x=0,\\
	-1 &\text{ when } x<0.
\end{cases}
\]
Also, when $\mathrm{sgn}(x)=\mathrm{sgn}(y)$ and $a\geq 0$, we have $\mathrm{sgn}(x+ay)=\mathrm{sgn}(x)=\mathrm{sgn}(y)$.

{The following lemma shows that the quantities in Lemma \ref{lemma.dece} should have the same sign.}
\begin{lemma} \label{lemma.sign}
If $\lambda_n^*>0$ for all $n=0,\ldots,N-1$, then
\begin{enumerate}
\item $\mathrm{sgn}(\Delta^2 h(n)) = \mathrm{sgn}(\Delta w(n) + \gamma_n \Delta^2 h(n-1))$ when $n=1,\ldots,N-2$.
\item $\mathrm{sgn}(\Delta^2 h(0)) = \mathrm{sgn}(\Delta w(0)). $
\end{enumerate}
\end{lemma}
\begin{proof}{Proof of Lemma \ref{lemma.sign}}
From Lemma \ref{lemma.dece}, if $1\leq n\leq N-2$ and $\Delta^2 h(n)>0$, then $\Delta w(n) + \gamma_n\Delta^2 h(n-1) \geq \lambda_n^* \Delta^2 h(n) >0 $. Conversely, when $\Delta w(n) + \Delta^2 h(n-1) > 0$, we have $\lambda_{n+1}^*\Delta^2 h(n)\geq \Delta w(n) + \gamma_n\Delta^2 h(n-1) > 0$. 
Since $\lambda_{n+1}^*>0$, we conclude $\Delta^2 h(n)>0$. Similarly, we get the same result when $\Delta^2h(n)<0$. When $\Delta^2 h(n)=0$, we have 
\[
0=\lambda_n^* \Delta^2 h(n) \leq \Delta w(n) +\gamma_n \Delta^2 h(n-1) \leq \lambda_{n+1}^* \Delta^2 h(n)=0 
\Rightarrow \Delta w(n) +\gamma_n \Delta^2 h(n-1) = 0. 
\]
Conversely, when $\Delta w(n) +\gamma_n \Delta^2 h(n-1)=0$, we have $\lambda_n^* \Delta^2 h(n) \leq 0 \leq \lambda_{n+1}^* \Delta^2 h(n)$. 
Since $\lambda_n^*,\lambda_{n+1}^*>0$, we conclude $\Delta^2 h(n)=0$. This proves the first part of the lemma and part 2 follows similarly. \Halmos
\end{proof}

\noindent {The next lemma shows that the signs of $\Delta w(n)$ and $\Delta^2 h(n) $ are the same for all $0\leq n<m$.}
\begin{lemma} \label{lemma.sign.n0}
If $\lambda_n^*>0$ for all $n=0,\ldots,N-1$, then for $n\leq N-2$ and $0\leq n< m$, we have $\mathrm{sgn}(\Delta^2 h(n)) = \mathrm{sgn}(\Delta w(n))= \mathrm{sgn}(\Delta w(0))$.
\end{lemma}

\begin{proof}{Proof of Lemma \ref{lemma.sign.n0}}
We use induction. For $n=0$, Lemma \ref{lemma.sign} implies $\mathrm{sgn}(\Delta^2 h(0)) = \mathrm{sgn}(\Delta w(0))$. 
Now assume $\mathrm{sgn}(\Delta^2 h(k)) = \mathrm{sgn}(\Delta w(k))= \mathrm{sgn}(\Delta w(0))$ for $k=0,\ldots,n-1$ and $n<m-1$. We will show the statement holds for $k=n$. 
For $1\leq n < m$, from (\ref{Delta2h}), we have 
\begin{align*}
\Delta w(n) &=  (\mu+\theta_s) \Delta h(n) + c_h + c_s\theta_s = (\mu+\theta_s) \Delta h(n-1) + (\mu+\theta_s)\Delta^2h(n-1) + c_h + c_s\theta_s\\
&=\Delta w(n-1) + (\mu+\theta_s)\Delta^2h(n-1). \addrefnum \label{Delta.w.induction.leq.m.step}
\end{align*}
From the fact that $\mu+\theta_s>0$ and the induction hypothesis, we have
\[
\mathrm{sgn} (\Delta w(n)) = \mathrm{sgn} (\Delta w(n-1) + (\mu+\theta_s)\Delta^2h(n-1)) =\mathrm{sgn}(\Delta w(0)).
\]
Since $\gamma_n>0$, from Lemma \ref{lemma.sign} and the induction hypothesis, we now have $\mathrm{sgn}(\Delta^2 h(n)) = \mathrm{sgn} (\Delta w(n) + \gamma_n \Delta^2 h(n-1)) =  \mathrm{sgn}(\Delta w(0))$. \Halmos
\end{proof}

{The next lemma establishes conditions for $\Delta w(n)$ and $\Delta^2 h(n-1)$ to have the same sign when $n \leq m$.}

\begin{lemma} \label{lemma.sign.nm}
Suppose that $\lambda_n^*>0$ for all $n=0,\ldots,N-1$. If $\mathrm{sgn}(\Delta w(n)) = \mathrm{sgn}(\Delta^2h(n-1))$ for $m\leq n\leq N-1$, then we have 
$\mathrm{sgn}(\Delta w(i)) = \mathrm{sgn}(\Delta^2h(i-1))= \mathrm{sgn}(\Delta^2h(n-1))$ for all $n\leq i\leq N-1$.
\end{lemma}

\begin{proof}{Proof of Lemma \ref{lemma.sign.nm}}
We use induction. The result immediately holds for $i=n$. Now assume $\mathrm{sgn}(\Delta w(k)) = \mathrm{sgn}(\Delta^2h(k-1))$ for $k=n,\ldots,i \leq N-2$ and we will show that the result holds for $k=i+1$. Since $\gamma_i>0$, we have
$
\mathrm{sgn} (\Delta w(i) + \gamma_i \Delta^2 h(i-1)) = \mathrm{sgn}(\Delta w(i)) = \mathrm{sgn}(\Delta^2h(i-1))
$.  
Since $m\leq i\leq N-2$, from (\ref{Delta2h}), we derive the relation between $\Delta w(i+1)$ and $\Delta w(i)$ as
\begin{align*}
\Delta w(i+1) &=  \theta_q \Delta h(i+1) + c_h + c_q\theta_q = \theta_q \Delta h(i) + \theta_q \Delta^2 h(i) + c_h + c_q\theta_q =\Delta w(i) + \theta_q \Delta^2h(i). \addrefnum \label{Delta.w.induction.step}
\end{align*}
From Lemma \ref{lemma.sign}, we have $\mathrm{sgn}(\Delta^2 h(i)) = \mathrm{sgn} (\Delta w(i) + \gamma_i\Delta^2 h(i-1)) =\mathrm{sgn}(\Delta w(i)) = \mathrm{sgn}(\Delta^2h(i-1))$. Since $\theta_q>0$, we have the sign of $\Delta w(i+1)$ is the same as that of $\Delta w(i)$ and $\Delta^2 h(i)$ by
\[
\mathrm{sgn} (\Delta w(i+1)) = \mathrm{sgn} (\Delta w(i) + \theta_q \Delta^2h(i)) = \mathrm{sgn} (\Delta^2 h(i)).
\]
Therefore, we show $\mathrm{sgn}(\Delta w(i)) = \mathrm{sgn}(\Delta^2h(i-1))= \mathrm{sgn}(\Delta^2h(n-1))$ for $n\leq i\leq N-1$. \Halmos
\end{proof}

Now we are ready to prove Lemma \ref{lemma.mono}.
\begin{proof}{Proof of Lemma \ref{lemma.mono}}
Assume that there exists some $n$ such that $\Delta^2 h(n)>0$ where $0\leq n <m$. 
By Lemma \ref{lemma.sign.n0}, we have $\Delta w(m-1)>0$ and $\Delta^2 h(m-1)>0$. From (\ref{Delta2h}), we then have
\begin{align*}
\Delta w(m) =& \theta_q \Delta h(m) + c_h +  c_q\theta_q  
= \theta_q \Big( \Delta h(m-1) + \frac{c_h}{\theta_q} +c_q \Big) + \theta_q \Delta^2 h(m-1) \addrefnum \label{Delta.w.m.induction.step} \\ 
\geq & \theta_q \Big( \Delta h(m-1) + \frac{c_h+c_s\theta_s}{\mu+\theta_s} \Big) + \theta_q \Delta^2 h(m-1) \text{ since } \frac{c_h+c_s\theta_s}{\mu+\theta_s}\leq \frac{c_h}{\theta_q}+c_q\\
=&\frac{\theta_q}{\mu+\theta_s} \Big( (\mu+\theta_s) \Delta h(m-1) + c_h+ c_s\theta_s \Big) + \theta_q \Delta^2 h(m-1)\\
=& \frac{\theta_q}{\mu+\theta_s} \Delta w(m-1) + \theta_q \Delta^2 h(m-1) > 0.
\end{align*}
Given that $\Delta w(m) >0$ and $\Delta^2 h(m-1)>0$, from Lemma \ref{lemma.sign.nm}, we have $\Delta w(N-1) >0$ and $\Delta^2 h(N-2)>0$, which contradicts Part 3 of Lemma \ref{lemma.dece}. This indicates that $\Delta^2 h(n)\leq 0$ for $0\leq n<m$ and, from Lemma \ref{lemma.sign.n0}, $\Delta w(n)\leq 0$ for $0\leq n <m$.

Now suppose there exists $m\leq n\leq N-2$ such that $\Delta^2 h(n)>0$ and let $n$ be the smallest such $n$. Then $\Delta w(n) + \gamma_n \Delta^2 h(n-1)>0$ from Lemma \ref{lemma.sign}. Since $\Delta h^2(n-1)\leq 0$ by the definition of $n$, then $\Delta w(n)>0$. From (\ref{Delta.w.induction.step}), we have $\Delta w(n+1)=\Delta w(n)+\theta_q \Delta^2 h(n)$. Thus, we have $w(n+1)>0$ and $\Delta^2 h(n)>0$. From Lemma \ref{lemma.sign.nm}, we conclude $\Delta w(N-1) >0$ and $\Delta^2 h(N-2)>0$, which contradicts Part 3 of Lemma \ref{lemma.dece}. This indicates that $\Delta^2 h(n)\leq 0$ for all $m\leq n\leq N-2$. \Halmos
\end{proof}

Now we move on to the proof of Lemma \ref{lemma.uni_modal}.

\begin{proof}{Proof of Lemma \ref{lemma.uni_modal}}
When $\frac{c_h+c_s\theta_s}{\mu+\theta_s}\leq \frac{c_h}{\theta_q}+c_q$, the results follow from Lemma \ref{lemma.mono} with $\tilde{n} = 0$. 
Then, we only need to prove the statements when $\frac{c_h+c_s\theta_s}{\mu+\theta_s}> \frac{c_h}{\theta_q}+c_q$. First, we assume $N\geq m+2$ and consider the cases $\tilde{n} = N-2$, $m\leq \tilde{n} \leq N-3$, and $0\leq \tilde{n}\leq m-1$ separately.

When $\tilde{n} = N-2$, by the definition of $\tilde{n}$, $\Delta^2(n)\leq 0$ for all $\tilde{n}\leq n\leq N-2$ follows immediately. Since $N\geq m+2$, then $\tilde{n}\geq m$.

When $m\leq \tilde{n}\leq N-3$, the second statement follows immediately and we have $\Delta^2 h(\tilde{n})\leq 0$ and $\Delta^2 h (\tilde{n}-1)>0$. From Lemma \ref{lemma.sign}, we have
\[
\Delta^2 h(\tilde{n})\leq 0 \Rightarrow \Delta w(\tilde{n}) + \gamma_{\tilde{n}} \Delta^2 h(\tilde{n}-1)\leq 0 \Rightarrow \Delta w(\tilde{n})<0 
\text{ since } \Delta^2 h (\tilde{n}-1)>0.
\]
Since $m\leq \tilde{n}\leq N-3$, then (\ref{Delta.w.induction.step}) implies $\Delta w(\tilde{n}+1) =\Delta w(\tilde{n}) + \theta_q \Delta^2h(\tilde{n})< 0$. 
If $\Delta^2 h(\tilde{n})<0$, then from Lemma \ref{lemma.sign.nm}, we have $\Delta^2 h(n)< 0$ for all $n\geq \tilde{n}$. If $\Delta^2 h(\tilde{n}) = 0$, then we have
$\Delta w(\tilde{n}+1) + \gamma_{\tilde{n}+1} \Delta^2 h(\tilde{n})< 0 $, and Lemma \ref{lemma.sign} implies $\Delta^2 h(\tilde{n}+1)<0$. For $m\leq \tilde{n}\leq N-3$, (\ref{Delta.w.induction.step}) implies $\Delta w(\tilde{n}+2) =\Delta w(\tilde{n}+1) + \theta_q \Delta^2h(\tilde{n}+1)< 0$. It now follows from Lemma \ref{lemma.sign.nm} that $\Delta^2 h(n)<0$ for all $n\geq \tilde{n}+1$. Combining all cases above, we have $\Delta^2 h(\tilde{n})\leq 0$ for all $n\geq \tilde{n}$.

When $0\leq \tilde{n}\leq m-1$, Lemma \ref{lemma.sign.n0} yields $\Delta w(n)\leq 0$ and $\Delta^2 h(n)\leq 0$ for $0\leq n\leq m-1$. Then, we have $\Delta^2 h(0)\leq 0$, which means $\tilde{n}=0$, proving the second statement. Also, we have $\Delta w(m-1)\leq 0$ and $\Delta^2 h(m-1)\leq 0$. Then, from (\ref{Delta.w.m.induction.step}), we deduce
\begin{align*}
	\Delta w(m) =& \theta_q \Big( \Delta h(m-1) + \frac{c_h}{\theta_q} +c_q \Big) + \theta_q \Delta^2 h(m-1)\\
	< & \theta_q \Big( \Delta h(m-1) + \frac{c_h+c_s\theta_s}{\mu+\theta_s} \Big) + \theta_q \Delta^2 h(m-1) \text{ since } \frac{c_h+c_s\theta_s}{\mu+\theta_s}> \frac{c_h}{\theta_q}+c_q\\
	=& \frac{\theta_q}{\mu+\theta_s} \Delta w(m-1) + \theta_q \Delta^2 h(m-1) \leq 0 \text{ which implies } \Delta w(m)<0.
\end{align*}
Given $\Delta w(m)< 0$ and $\Delta^2 h(m-1)\leq 0$, if $\Delta^2 h(m-1)<0$, then from Lemma \ref{lemma.sign.nm}, $\Delta^2 h(n)< 0$ for all $n\geq m$. If $\Delta^2 h(m-1)=0$, then we have $\Delta w(m) + \gamma_{m} \Delta^2 h(m-1)< 0 $ and Lemma \ref{lemma.sign} implies that $\Delta^2 h(m)<0$. Then,  (\ref{Delta.w.induction.step}) implies $\Delta w(m+1) =\Delta w(m) + \theta_q \Delta^2h(m)< 0$. It now follows from Lemma \ref{lemma.sign.nm} that $\Delta^2 h(n)<0$ for all $n\geq m+1$. 

Now consider $N\leq m+1$. Since $\tilde{n}\leq N-2\leq m-1$, then using the argument for the case $0\leq \tilde{n}\leq m-1$, we have $\tilde{n}=0$ and $\Delta^2 h(n)\leq 0$ for all $n\geq \tilde{n}$.

Combining all cases above, we show that $\Delta^2 h(n)\leq 0$ for all $n\geq \tilde{n}$ and either $\tilde{n}\geq m$ or $\tilde{n}=0$. \Halmos
\end{proof}

\section{Proofs of Results in Section \ref{sec.Heuristics}} \label{ec.proof.section.5}

In this section, we prove the major results in Section \ref{sec.Heuristics}. We prove the properties of the optimal threshold and the performance of the cutoff-static policy in Section \ref{ec.proof.section.5.1} and show the results related to the two-price policy in Section \ref{ec.proof.section.5.2}.

\subsection{Proofs of Results in Section \ref{subsec.cutoff.static}} \label{ec.proof.section.5.1}

We first prove Theorem \ref{Theorem.skip.server} by showing that the increment between $g_C(K+1)$ and $g_C(K)$ is non-negative for $K=0,\ldots,m-1$.

\begin{proof}{Proof of Theorem \ref{Theorem.skip.server}}
	Recall the definition of $a_n(\delta)$ from the Section 3. Given $n\leq m$,
	$$
	n\times a_n(\delta) = n\times \prod_{i=1}^n\dfrac{\delta}{i(\mu+\theta_s)}=\dfrac{\delta}{\mu+\theta_s}\times\prod_{i=1}^{n-1}\dfrac{\delta}{i(\mu+\theta_s)}=\dfrac{\delta}{\mu+\theta_s}\times a_{n-1}(\delta)=a_1(\delta)a_{n-1}(\delta).
	$$
	We consider the cases $g_C^*(1)>0$ and  $g_C^*(1)=0$ separately.
	
	If $g_C^*(1)>0$, we have 
	\[
	g_C^*(1)=\dfrac{\delta_1\Bar{F}^{-1}(\delta_1/\Lambda)-a_1(\delta_1)(c_h+c_s\theta_s)}{1+a_1(\delta_1)}>0 \Rightarrow \delta_1\Bar{F}^{-1}(\delta_1/\Lambda)-a_1(\delta_1)(c_h+c_s\theta_s) > 0.
	\]
	Now consider $K=2,\ldots,m$, for any $\delta\in[0,\Lambda]$, we have
	\begin{align*}
		g_C^*(K) &\geq g_C^{\delta}(K) = \dfrac{\sum_{n=0}^{K-1} a_n(\delta)\times[\delta\Bar{F}^{-1}(\delta/\Lambda)]-\sum_{n=1}^K a_n(\delta)\times[n(c_h+c_s\theta_s)]}{\sum_{n=0}^K a_n(\delta)}\\
		&=\dfrac{\sum_{n=0}^{K-1} a_n(\delta)\times[\delta\Bar{F}^{-1}(\delta/\Lambda)]-\sum_{n=1}^K a_{n-1}(\delta)a_1(\delta)\times(c_h+c_s\theta_s)}{\sum_{n=0}^K a_n(\delta)}\\
		&=\dfrac{\sum_{n=0}^{K-1} a_n(\delta)\times[\delta\Bar{F}^{-1}(\delta/\Lambda)-a_1(\delta)(c_h+c_s\theta_s)]}{\sum_{n=0}^K a_n(\delta)},
	\end{align*}
	where the inequality holds from the definition of $g_C^*(K)$. Plugging in $\delta=\delta_1$, we have \[
	g_C^*(K)\geq \dfrac{\sum_{n=0}^{K-1} a_n(\delta_1)\times[\delta_1\Bar{F}^{-1}(\delta_1/\Lambda)-a_1(\delta_1)(c_h+c_s\theta_s)]}{\sum_{n=0}^K a_n(\delta_1)}>0.
	\]
	Then, given $g_C^*(K)>0$, we have
	\begin{align*}
		g_C^*(K)=\dfrac{\sum_{n=0}^{K-1} a_n(\delta_K)\times[\delta_K\Bar{F}^{-1}(\delta_K/\Lambda)-a_1(\delta_K)(c_h+c_s\theta_s)]}{\sum_{n=0}^K a_n(\delta_K)}>0,
	\end{align*}
	which implies that $ \delta_K\Bar{F}^{-1}(\delta_K/\Lambda)-a_1(\delta_K)(c_h+c_s\theta_s)>0.$
	
	Now consider the difference between $g_C^*(K+1)$ and $g_C^*(K)$ for $K=1,\ldots,m-1$:
	\begin{align*}
		&g_C^*(K+1)-g_C^*(K) \geq  g_C^{\delta_{K}}(K+1) - g_C^{\delta_{K}}(K)\\
		 = & \dfrac{\sum_{n=0}^{K} a_n(\delta_K) \times[\delta_K\Bar{F}^{-1}(\delta_K/\Lambda)-a_1(\delta_K)(c_h+c_s\theta_s)]}{\sum_{n=0}^{K+1} a_n(\delta_K)}\\
		&- \dfrac{\sum_{n=0}^{K-1} a_n(\delta_K)\times[\delta_K\Bar{F}^{-1}(\delta_K/\Lambda)-a_1(\delta_K)(c_h+c_s\theta_s)]}{\sum_{n=0}^K a_n(\delta_K)}\\
		=& [\delta_K\Bar{F}^{-1}(\delta_K/\Lambda)-a_1(\delta_K)(c_h+c_s\theta_s)] \Big( \dfrac{\sum_{n=0}^{K} a_n(\delta_K)}{\sum_{n=0}^{K+1} a_n(\delta_K)} - \dfrac{\sum_{n=0}^{K-1} a_n(\delta_K)}{\sum_{n=0}^{K} a_n(\delta_K)} \Big)\\
		=& [\delta_K\Bar{F}^{-1}(\delta_K/\Lambda)-a_1(\delta_K)(c_h+c_s\theta_s)] \Big( \dfrac{a_K(\delta_{K})}{\sum_{n=0}^{K} a_n(\delta_K)} - \dfrac{a_{K+1}(\delta_{K})}{\sum_{n=0}^{K+1} a_n(\delta_K)}  \Big)\\
		=& \dfrac{[\delta_K\Bar{F}^{-1}(\delta_K/\Lambda)-a_1(\delta_K)(c_h+c_s\theta_s)]a_K(\delta_{K})}{\sum_{n=0}^K a_n(\delta_K)\times \sum_{n=0}^{K+1} a_n(\delta_K)}\Big(\sum_{n=0}^{K+1} a_n(\delta_K)-\frac{\delta_K}{(K+1)(\mu+\theta_s)}\sum_{n=0}^{K} a_n(\delta_K)\Big). \addrefnum\label{diff.g_c^k+1,k}
	\end{align*}
	Moreover, 
	\begin{align*}
		&\sum_{n=0}^{K+1} a_n(\delta_K)-\frac{\delta_K}{(K+1)(\mu+\theta_s)}\sum_{n=0}^{K} a_n(\delta_K)\\
		=& 1+\sum_{n=1}^{K+1} \frac{\delta_K}{n(\mu+\theta_s)}a_{n-1}(\delta_K) -\frac{\delta_K}{(K+1)(\mu+\theta_s)}\sum_{n=0}^{K} a_n(\delta_K)\\
		= &1+ \sum_{n=0}^{K} \frac{\delta_{K}}{(n+1)(\mu+\theta_s)}a_n(\delta_K)-\frac{\delta_K}{(K+1)(\mu+\theta_s)}\sum_{n=0}^{K} a_n(\delta_K)\\
		= & 1+ \sum_{n=0}^{K} a_n(\delta_K)\dfrac{\delta_K}{\mu+\theta_s} \Big( \frac{1}{n+1}-\frac{1}{K+1} \Big) \geq 0 \text{ since } n\leq K.
	\end{align*} 
	Since $ \delta_K\Bar{F}^{-1}(\delta_K/\Lambda)-a_1(\delta_K)(c_h+c_s\theta_s)>0$, it now follows from (\ref{diff.g_c^k+1,k}) that $g_C^*(K+1)\geq g_C^*(K)$.
	 
	On the other hand, if $g_C^*(1)=0$, then we have $g_C^\delta(1)\leq 0$ for any $\delta\in [0,\Lambda]$,  which implies \[
	\delta\Bar{F}^{-1}(\delta/\Lambda)-a_1(\delta)(c_h+c_s\theta_s)\leq 0
	\] for all $\delta\in[0,\Lambda]$. Thus, for $K=2,\ldots,m$, we have 
	$$
	g_C^*(K)=\dfrac{\sum_{n=0}^{K-1} a_n(\delta_K)\times[\delta_K\Bar{F}^{-1}(\delta_K/\Lambda)-a_1(\delta_K)(c_h+c_s\theta_s)]}{\sum_{n=0}^K a_n(\delta_K)}\leq 0.
	$$
	Since $g_C^{\delta}(K)=0$ when $\delta=0$, then we must have $g_C^*(K)=0$. \Halmos
\end{proof}

{Before proving Propositions \ref{prop.upper.0} and \ref{prop.upper.refine}, define $\CR_n(\lambda)$ as follows: 
\[
\mathcal{R}_n(\lambda) = 
\begin{cases}
\lambda\Bar{F}^{-1}(\lambda/\Lambda) - \lambda \mathcal{C}_s  & \text{for } n<m,\\
\lambda\Bar{F}^{-1}(\lambda/\Lambda) - \lambda \left[ \frac{\gamma_m}{\gamma_{n+1}} \mathcal{C}_s+ \frac{\gamma_{n+1}-\gamma_{m}}{\gamma_{n+1}} \mathcal{C}_q \right] & \text{for } m\leq n<N,
\end{cases}
\]
where $\mathcal{C}_s = \frac{c_h+c_s\theta_s}{\mu+\theta_s}$ and $\mathcal{C}_q = \frac{c_h}{\theta_q}+c_q$. Since $\lambda_N=0$ from the definition in Section \ref{sec.model_description}, we have $\CR_N(\lambda_N)=\CR_N(0)=0$. We can rewrite the long-run average gain $g^{\vec{\lambda}}$ under policy $\vec{\lambda}$ as follows:
\begin{align*}
    g^{\Vec{\lambda}} =&\sum_{n=0}^{N-1} \lambda_n\bar{F}^{-1}(\lambda_n/\Lambda)P_n(\vec{\lambda}) - \sum_{n=1}^N\left[nc_h + \min\{n,m\}c_s\theta_s +\max\{n-m,0\}c_q\theta_q\right] P_n(\vec{\lambda})\\
    =&\sum_{n=0}^{N-1} \lambda_n\bar{F}^{-1}(\lambda_n/\Lambda)P_n(\vec{\lambda}) - \sum_{n=1}^N\left[nc_h + \min\{n,m\}c_s\theta_s +\max\{n-m,0\}c_q\theta_q\right] \frac{\lambda_{n-1}}{\gamma_n}P_{n-1}(\vec{\lambda})\\
    =& \sum_{n=0}^{m-1} \lambda_n\bar{F}^{-1}(\lambda_n/\Lambda)P_n(\vec{\lambda}) - \sum_{n=1}^{m} n(c_h+c_s\theta_s) \frac{\lambda_{n-1}}{n(\mu+\theta_s)}P_{n-1}(\vec{\lambda}) + \sum_{n=m}^{N-1} \lambda_n\bar{F}^{-1}(\lambda_n/\Lambda)P_n(\vec{\lambda}) \\
    &- \sum_{n=m+1}^{N} \left[m(\mu+\theta_s)\frac{c_h+c_s\theta_s}{\mu+\theta_s} + (n-m)\theta_q \left(\frac{c_h}{\theta_q}+c_q\right)\right] \frac{\lambda_{n-1}}{m(\mu+\theta_s) + (n-m)\theta_q}P_{n-1}(\vec{\lambda})\\
    =& \sum_{n=0}^{m-1} \left[ \lambda_n\Bar{F}^{-1}(\lambda_n/\Lambda) - \lambda_n\CC_s \right] P_n(\Vec{\lambda}) 
    + \sum_{n=m}^{N-1} \left[ \lambda_n\Bar{F}^{-1}(\lambda_n/\Lambda) - \lambda_n \left(  \frac{\gamma_m}{\gamma_{n+1}}\CC_s + \frac{\gamma_{n+1}-\gamma_m}{\gamma_{n+1}}\CC_q\right) \right] P_n(\Vec{\lambda})\\
    =& \sum_{n=0}^{N-1} \mathcal{R}_n(\lambda_n) P_n(\vec{\lambda}) = \sum_{n=0}^{N} \mathcal{R}_n(\lambda_n) P_n(\vec{\lambda}). \addrefnum \label{gain.reformulation}
\end{align*}
The expression of the long-run average gain in (\ref{gain.reformulation}) is used in the proofs on Propositions \ref{prop.upper.0} and \ref{prop.upper.refine}. To prove the upper bound on the optimal threshold $K_C$, we will show that there exists an integer $n$ such that any cutoff-static policy with non-zero arrival rate beyond $n$ will be sub-optimal.}
{
\begin{proof}{Proof of Proposition \ref{prop.upper.0}.}
Since $\Bar{F}(\min\{\CC_s,\CC_q\})=0$, we have $\Bar{F}^{-1}(\lambda/\Lambda)\leq \min  \{\CC_s,\CC_q\}$ for all $\lambda\in [0,\Lambda]$. When $0\leq n<m$, we have $\CR_n(\lambda)=\lambda\Bar{F}^{-1}(\lambda/\Lambda) - \lambda \mathcal{C}_s \leq \lambda\Bar{F}^{-1}(\lambda/\Lambda) - \lambda\min  \{\CC_s,\CC_q\}\leq 0$ for all $\lambda\in [0,\Lambda]$. For $m\leq n<N$, we have $\frac{\gamma_m}{\gamma_{n+1}} \mathcal{C}_s+ \frac{\gamma_{n+1}-\gamma_{m}}{\gamma_{n+1}} \mathcal{C}_q \geq \min\{\CC_s,\CC_q\}$, which implies 
\[
\CR_n(\lambda) = \lambda\Bar{F}^{-1}(\lambda/\Lambda) - \lambda\left[ \frac{\gamma_m}{\gamma_{n+1}} \mathcal{C}_s+ \frac{\gamma_{n+1}-\gamma_{m}}{\gamma_{n+1}} \mathcal{C}_q \right]  \leq \lambda\Bar{F}^{-1}(\lambda/\Lambda) - \lambda\min             \{\CC_s,\CC_q\}\leq 0.
\]
This indicates that $\CR_n(\lambda)\leq 0$ for all $\lambda\in[0,\Lambda]$. Thus, for any policy $\vec{\lambda}$, we have $\CR_n(\lambda_n)\leq 0$ for all $n=0,\ldots,N-1$ and $g^{\vec{\lambda}}=\sum_{n=0}^{N-1} \CR_n(\lambda_n) P_n(\vec{\lambda})\leq 0$ from (\ref{gain.reformulation}). Since $\lambda_n=0$ for all $n=0,\ldots,N-1$ yields a long-run average gain of $0$, we have $g^*=0$ and $\lambda_n=0$ for all $n=0,\ldots,N-1$ as the optimal policy. Since the argument above holds for all policies, it holds for the cutoff-static policy, which implies $K_C=0$. \Halmos
\end{proof}
}

{
Before proving Proposition \ref{prop.upper.refine}, we show that $\CR_n(\lambda)$ is a non-decreasing function with respect to $0\leq n<N$ for a fixed $\lambda$ when $\CC_s\leq \CC_q$. We provide this useful observation in the following remark.
\begin{remark} \label{remark.non-decreasing.CR}
Observe that $\frac{\gamma_{n+1}-\gamma_{m}}{\gamma_{n+1}} =\frac{(n+1-m)\theta_q}{m(\mu+\theta_s)+(n+1-m)\theta_q}$ is an increasing function with respect to $n\geq m-1$. When $\CC_s\leq \CC_q$, then $\frac{\gamma_m}{\gamma_{n+1}} \mathcal{C}_s+ \frac{\gamma_{n+1}-\gamma_{m}}{\gamma_{n+1}} \mathcal{C}_q$ is an increasing function of $n\geq m-1$, and $\frac{\gamma_m}{\gamma_{n+1}} \mathcal{C}_s+ \frac{\gamma_{n+1}-\gamma_{m}}{\gamma_{n+1}} \mathcal{C}_q=\CC_s$ when $n=m-1$. Thus, $\CR_n(\lambda)$, which is as a constant when $n<m$, is a non-decreasing function of $n$ for $0\leq n<N$ for any fixed $\lambda\in[0,\Lambda]$.
\end{remark}
}
\begin{proof}{Proof of Proposition \ref{prop.upper.refine}.}
{
 To show that $K_C\leq b$, it suffices to show that for all $\ell \geq b$, there exists some $k\leq b$ such that $g_C^*(\ell) \leq g_C^{\delta_\ell}(k)$ where $\delta_\ell$ is the optimal arrival rate with threshold $\ell$ (as defined in Definition \ref{def.best.cutoff}). If $\delta_\ell =0$, we have $g_C^*(\ell)=0$. Then, pick any $k\leq b$ to achieve $g_C^{\delta_\ell}(k)=0$. 
}

{
 Now we assume $\delta_{\ell}>0$. If $\CR_0(\delta_{\ell})\leq 0$, then $\CR_n(\delta_\ell)\leq 0$ for all $0\leq n<\ell$ from Remark \ref{remark.non-decreasing.CR}, which implies $g_C^*(\ell)=0$ and $\delta_\ell=0$. This reaches a contradiction. Thus, we conclude $\CR_0(\delta_\ell)>0$.
 }

{
 Since $\bar{F}(\frac{\gamma_m}{\gamma_{b+1}} \mathcal{C}_s+ \frac{\gamma_{b+1}-\gamma_{m}}{\gamma_{b+1}} \mathcal{C}_q)=0$, then $\bar{F}^{-1}(\delta_{\ell}/\Lambda)\leq \frac{\gamma_m}{\gamma_{b+1}} \mathcal{C}_s+ \frac{\gamma_{b+1}-\gamma_{m}}{\gamma_{b+1}} \mathcal{C}_q$, which implies 
\[
 \CR_b(\delta_{\ell}) = \delta_\ell\Bar{F}^{-1}(\delta_\ell/\Lambda) - \delta_\ell\left[ \frac{\gamma_m}{\gamma_{b+1}} \mathcal{C}_s+ \frac{\gamma_{b+1}-\gamma_{m}}{\gamma_{b+1}} \mathcal{C}_q \right] \leq 0,
 \]
and $R_n(\delta_\ell)\leq 0$ for all $n\geq b$ from Remark \ref{remark.non-decreasing.CR}. Given $\CR_0(\delta_\ell)>0$, $\CR_b(\delta_\ell)\leq 0$ and $\CR_n(\delta_\ell)$ is non-increasing, there exists $0<k\leq b$ where $\CR_{k-1}(\delta_\ell)>0$ and $\CR_{k}(\delta_\ell)\leq 0$. This also implies that $\CR_n(\delta_\ell) >0$ for all $n<k$ and $\CR_n(\delta_\ell)\leq 0$ for all $n\geq k$.
}

{
Recall the definition of the steady-state probability $P_n(\pi_\ell^{\delta_\ell})$ of state $n$ under the cutoff-static policy $\pi_\ell^{\delta_\ell}$ in Section \ref{subsec.cutoff.static}. Since $\sum_{n=0}^k a_n(\delta_\ell)\leq \sum_{n=0}^\ell a_n(\delta_\ell)$ from $k\leq b\leq \ell$, then we have $P_n(\pi_k^{\delta_\ell})\geq P_n(\pi_\ell^{\delta_\ell})$ for all $0\leq n\leq k$. Therefore, employing (\ref{gain.reformulation}), we conclude that
\begin{align*}
    g_C^*(k)\geq & g_C^{\delta_\ell}(k) = \sum_{n=0}^{k-1} \CR_n(\delta_\ell) P_n(\pi_k^{\delta_\ell}) \underset{}{\geq} \sum_{n=0}^{k-1} \CR_n(\delta_\ell) P_n(\pi_\ell^{\delta_\ell}) \\
    \underset{}{\geq } &\sum_{n=0}^{k-1} \CR_n(\delta_\ell) P_n(\pi_\ell^{\delta_\ell}) + \sum_{n=k}^{\ell-1} \CR_n(\delta_\ell) P_n(\pi_\ell^{\delta_\ell})
    = \sum_{n=0}^{\ell-1} \CR_n(\delta_\ell) P_n(\pi_\ell^{\delta_\ell}) = g_C^*(\ell).
\end{align*}
The first inequality holds due to $\CR_n(\delta_\ell)>0$ and $P_n(\pi_k^{\delta_\ell})\geq P_n(\pi_\ell^{\delta_\ell})$ for all $0\leq n\leq k$. The second inequality holds because $\CR_n(\delta_\ell)\leq 0$ and $P_n(\pi_\ell^{\delta_\ell})> 0$ for all $n\geq k$. \Halmos
}
\end{proof}

{The rest of this section focuses on proving the results related to the performance guarantees of the best cutoff static policy.}
To prove Theorem \ref{Theorem.performance.reve}, we show $P_K(\pi_K^\delta)$ is an increasing function of $\delta$.
\begin{lemma} \label{lemma.a_n.derivative}
The steady-state probability $P_K(\pi_K^\delta)$ at state $K$ under policy $\pi_K^\delta$ is a non-decreasing function of $\delta$ for $\delta\in[0,\Lambda]$.
\end{lemma}

\begin{proof}{Proof.}
From Definition \ref{def.cut.off.static}, we compute the derivative of $a_n(\delta)$ for $n\geq 0$:
\begin{align*}
a'_n(\delta) = \dfrac{d}{d\delta} a_n(\delta) =  \dfrac{d}{d\delta} \Bigg( \dfrac{\delta^n}{\prod_{i=1}^n \gamma_i} \Bigg) =\dfrac{n\delta^{n-1}}{\prod_{i=1}^n \gamma_i}  = \dfrac{n}{\delta} \dfrac{\delta^n}{\prod_{i=1}^n \gamma_i} = \dfrac{n}{\delta}\times  a_n(\delta). \addrefnum \label{a_n_derivative}
\end{align*}
From Definition \ref{def.cut.off.static} and (\ref{a_n_derivative}), we have 
\begin{align*}
\dfrac{d}{d\delta} P_K(\pi_K^\delta) =& \dfrac{d}{d\delta} \Bigg( \dfrac{a_K(\delta)}{\sum_{n=0}^K a_n(\delta)} \Bigg) = \dfrac{1}{[\sum_{n=0}^K a_n(\delta)]^2} \Big[ {a'_K(\delta) \times \sum_{n=0}^K a_n(\delta) -  a_K(\delta) \times \sum_{n=0}^K a'_n(\delta)} \Big]\\
=& \dfrac{1}{[\sum_{n=0}^K a_n(\delta)]^2} \Big[ \dfrac{K}{\delta}\times a_K(\delta) \times \sum_{n=0}^K a_n(\delta) -  a_K(\delta) \times \sum_{n=0}^K \dfrac{n}{\delta}\times a_n(\delta) \Big]\\
=& \dfrac{a_K(\delta) /\delta}{[\sum_{n=0}^K a_n(\delta)]^2} \Big[ \sum_{n=0}^K (K-n) a_n(\delta) \Big] \geq 0,
\end{align*}
where the last inequality holds because $a_n(\delta)\geq 0$ for all $n\geq 0$ and $\delta\in[0,\Lambda]$. \Halmos
\end{proof}

\begin{proof}{Proof of Theorem \ref{Theorem.performance.reve}.}
We consider the long-run average revenue under policy $\pi_K^{\hat{\delta}}$:
$$\CR(\pi_K^{\hat{\delta}}) = \sum_{n=0}^{K-1} \hat{\delta}\Bar{F}^{-1}(\hat{\delta}/\Lambda) \times  P_n(\pi_K^{\hat{\delta}}) = \hat{\delta}\Bar{F}^{-1}(\hat{\delta}/\Lambda) \times [1-P_K(\pi_K^{\hat{\delta}})].$$
Since $F(\cdot)$ is a regular distribution and (\ref{effective.arrival.rate}) shows that $\hat{\delta}$ is a convex combination of $\lambda_n^*$ for $n=0,\ldots,N$, then, from Jensen Inequality, 
$$
\hat{\delta}\Bar{F}^{-1}(\hat{\delta}/\Lambda) \geq \sum_{n=0}^N \lambda_n^* \Bar{F}^{-1}(\lambda_n^*/\Lambda) \times P_n(\vec{\lambda}^*) =\CR(\Vec{\lambda}^*).
$$
Furthermore, from Lemma \ref{lemma.a_n.derivative}, we know that $P_K(\pi_K^\delta)$
is a non-decreasing function of $\delta\in[0,\Lambda]$. Thus, given the maximal possible arrival rate as $\Lambda$, we have
$$
1-P_K(\pi_K^{\hat{\delta}})\geq 1-\dfrac{a_K(\Lambda)}{\sum_{n=0}^K a_n(\Lambda)}.
$$
Therefore, we have $$
\CR(\pi_K^{\hat{\delta}}) = \hat{\delta}\Bar{F}^{-1}(\hat{\delta}/\Lambda) \times [1- P_K(\pi_K^{\hat{\delta}})] \geq \CR(\Vec{\lambda}^*) \times \Big(1-\dfrac{a_K(\Lambda)}{\sum_{n=0}^K a_n(\Lambda)}\Big). \Halmos
$$
\end{proof}

Before proving Theorem \ref{Theorem.performance.cost}, we define $L^*$ as the long-run expected number of customers in the system and $T^*$ as the long-run expected total time of a customer in the system under the optimal dynamic policy $\vec{\lambda}^*$. For customers who join the system with idling servers, they enter service immediately and may abandon during service or finish the service. For those who join the system without idling servers, they join the queue and either abandon during waiting or enter service. {The average time a customer who enters service spends in the system is
1/($\mu+\theta_s$). On the other hand, if the customer eventually abandons, 
the average time that customer spends in the system is no smaller than 1/$\theta_q$. Hence, $T^*$ must exceed the minimum of $1/(\mu+\theta_s)$ and $1/\theta_q$, which indicates $T^*\geq 1/ \max \{\mu+\theta_s,\theta_q\}$.}
From \citet{Dai2013}, we know Little's Law applies to our queueing system with abandonments under the optimal dynamic policy $\vec{\lambda}^*$. Hence, we have
\[
L^* = \hat{\delta} \times T^* \geq \dfrac{\hat{\delta}}{\max \{\mu+\theta_s,\theta_q\}} \addrefnum\label{Little.law}.
\] 
\begin{proof}{Proof of Theorem \ref{Theorem.performance.cost}.}
When $c_q\theta_q\geq c_s\theta_s>0$, consider the long-run average total cost (\ref{total.exp.cost}) under the optimal dynamic policy $\vec{\lambda}^*$ and apply (\ref{Little.law}):
\begin{align*}
	\CC(\Vec{\lambda}^*) &\geq   \sum_{n=0}^N (nc_h + nc_s\theta_s)\times P_n(\Vec{\lambda}^*) =  (c_h + c_s\theta_s) \times \sum_{n=0}^N n P_n(\Vec{\lambda}^*)\\
	& = (c_h+c_s\theta_s) \times L^* \geq (c_h+c_s\theta_s) \times \dfrac{\hat{\delta}}{\max \{\mu+\theta_s,\theta_q\}},
\end{align*}
which implies 
\[
c_h+c_s\theta_s  \leq \dfrac{\CC(\Vec{\lambda}^*)}{\hat{\delta}/(\max \{\mu+\theta_s,\theta_q\})}. \addrefnum \label{Little.law.c_s.theta_s}
\]
Now consider the long-run average cost under a cutoff-static policy $\pi_K^{\hat{\delta}}$ and apply (\ref{Little.law.c_s.theta_s}): 
\begin{align*}
\CC(\pi_K^{\hat{\delta}}) &= \sum_{n=1}^m n(c_h+c_s\theta_s) \times P_n(\pi_K^{\hat{\delta}}) + \sum_{n=m+1}^K (n c_h+ mc_s\theta_s + (n-m)c_q\theta_q)\times P_n(\pi_{K}^{\hat{\delta}})\\
& = \sum_{n=1}^K n(c_h+c_s\theta_s) \times P_n(\pi_K^{\hat{\delta}}) + \sum_{n=m+1}^K (n-m)(c_q\theta_q-c_s\theta_s)\times P_n(\pi_{K}^{\hat{\delta}})\\
& \leq \CC(\Vec{\lambda^*})\times \Big( \dfrac{\sum_{n=1}^m n  P_n(\pi_K^{\hat{\delta}})}{\hat{\delta}/(\max \{\mu+\theta_s,\theta_q\})} + \dfrac{c_q\theta_q-c_s\theta_s}{c_h+c_s\theta_s} \times \dfrac{\sum_{n=m+1}^K (n-m) P_n(\pi_K^{\hat{\delta}})}{{\hat{\delta}}/(\max \{\mu+\theta_s,\theta_q\})}  \Big)\\
& \leq  \CC(\Vec{\lambda^*})\times \max_{\delta\in [0,\Lambda]} \Big( \dfrac{\sum_{n=1}^m n  P_n(\pi_K^{{\delta}})}{{\delta}/(\max \{\mu+\theta_s,\theta_q\})} + \dfrac{c_q\theta_q-c_s\theta_s}{c_h+c_s\theta_s} \times \dfrac{\sum_{n=m+1}^K (n-m) P_n(\pi_K^{{\delta}})}{{\delta}/(\max \{\mu+\theta_s,\theta_q\})}  \Big).
\end{align*}

When $c_s\theta_s\geq c_q\theta_q>0$, an similar argument to that in (\ref{Little.law.c_s.theta_s}) and  (\ref{Little.law}) yields
\begin{align*}
\CC(\Vec{\lambda}^*) &\geq \sum_{n=1}^N (nc_h + nc_q\theta_q) \times P_n(\Vec{\lambda}^*) = (c_h+c_q\theta_q) \times L^*\geq (c_h+c_q\theta_q) \times \dfrac{\hat{\delta}}{\max \{\mu+\theta_s,\theta_q\}},
\end{align*}
which implies 
\[
c_h+c_q\theta_q \leq  \dfrac{\CC(\Vec{\lambda}^*)}{\hat{\delta}/(\max \{\mu+\theta_s,\theta_q\})} \addrefnum \label{Little.law.c_q.theta_q}.
\]
Now consider the long-run average cost under a cutoff-static policy $\pi_K^{\hat{\delta}}$ and apply (\ref{Little.law.c_q.theta_q}):
\begin{align*}
	\CC(\pi_K^{\hat{\delta}}) &= \sum_{n=1}^m n(c_h+c_s\theta_s) \times P_n(\pi_K^{\hat{\delta}}) + \sum_{n=m+1}^K (n c_h+ mc_s\theta_q + (n-m)c_q\theta_q)\times P_n(\pi_{K}^{\hat{\delta}})\\
	& = \sum_{n=1}^K n(c_h+c_q\theta_q) \times P_n(\pi_K^{\hat{\delta}}) + \sum_{n=1}^m n(c_s\theta_s-c_q\theta_q)\times P_n(\pi_{K}^{\hat{\delta}}) + \sum_{n=m+1}^K m(c_s\theta_s-c_q\theta_q)\times P_n(\pi_{K}^{\hat{\delta}})\\
	& \leq \CC(\Vec{\lambda^*})\times \Big( \dfrac{\sum_{n=1}^K n  P_n(\pi_K^{\hat{\delta}})}{\hat{\delta}/(\max \{\mu+\theta_s,\theta_q\})} + \dfrac{c_s\theta_s-c_q\theta_q}{c_h+c_q\theta_q} \times \dfrac{\sum_{n=1}^K \min\{n,m\} P_n(\pi_K^{\hat{\delta}})}{\hat{\delta}/(\max \{\mu+\theta_s,\theta_q\})}  \Big)\\
	& \leq  \CC(\Vec{\lambda^*})\times \max_{\delta\in [0,\Lambda]} \Big( \dfrac{\sum_{n=1}^K n  P_n(\pi_K^{{\delta}})}{{\delta}/(\max \{\mu+\theta_s,\theta_q\})} + \dfrac{c_s\theta_s-c_q\theta_q}{c_h+c_q\theta_q} \times \dfrac{\sum_{n=1}^K \min\{n,m\} P_n(\pi_K^{{\delta}})}{{\delta}/(\max \{\mu+\theta_s,\theta_q\})}  \Big).
\end{align*}
Combining all cases above, we prove the statement. \Halmos
\end{proof}

{

\begin{proof}{Proof of Proposition \ref{prop.better.bound}}
Given $\bar{F}(\CC_s)>0$, there exists $\delta\in[0,\Lambda]$ such that $\bar{F}^{-1}(\delta/\Lambda)-\CC_s > 0$. Then, from (\ref{gain.reformulation}), we have $g_C^*(1) \geq g_C^{\delta}(1) =\delta\bar{F}^{-1}(\delta/\Lambda) -\delta \CC_s>0$. Thus, from Theorem \ref{Theorem.skip.server}, Item 1, we have $K_C\geq m$. 
Since $\Bar{F}\left( \frac{\gamma_m}{\gamma_{m+1}} \CC_s +  \frac{\gamma_{m+1}-\gamma_m}{\gamma_{m+1}} \CC_q \right)=0$ and $\CC_s\leq \CC_q$, from Proposition \ref{prop.upper.refine}, we have $K_C\leq m$, which implies $K_C=m$. Then, the system becomes an Erlang Loss system with a single class of customers. Thus, according to \citet{adam&jiaqi2025}, we have $g_C^*/g^* \geq 15/19$.
\end{proof}

}

{ 
\subsection{Proofs of Results in Section \ref{subsec.two.price}} \label{ec.proof.section.5.2}

In this section, we prove the results related to the optimality gap between the best two-price policy and the optimal dynamic policy. 
To prove Theorem \ref{Theorem.tp.perform}, we derive upper bounds for the optimal dynamic gain $g^*$ and lower bounds for the best two-price policy $g_T^*$. We first state the following lemma that gives a lower bound on the server-idling probability for any arrival rate policy $\vec{\lambda}$.
\begin{lemma} \label{lemma.idle.probability.bound}
For any arrival rate policy $\vec{\lambda}$, we have 
\[
 \frac{\sum_{n=0}^{m-1} a_n(\vec{\lambda})}{\sum_{n=0}^{N} a_n(\vec{\lambda})} = \sum_{n=0}^{m-1} P_n(\vec{\lambda}) \geq \sum_{n=0}^{m-1} P_n(\pi_N^\Lambda) = \frac{\sum_{n=0}^{m-1} a_n(\Lambda)}{\sum_{n=0}^{N} a_n({\Lambda})} \geq P_s^\Lambda,
\]
where  $\pi_N^\Lambda$ and $P_s^\Lambda$ are provided in Definition \ref{def.cut.off.static} and (\ref{infinite buffer steady state}), respectively.
\end{lemma}
\begin{proof}{Proof of Lemma \ref{lemma.idle.probability.bound}.}
Suppose $0\leq i\leq m-1$ and $m\leq j\leq N$. For $i+1\leq \ell\leq j$, since $\lambda_{\ell-1}\leq \Lambda$, $\frac{\Lambda}{\gamma_{\ell}}\geq \frac{\lambda_{\ell-1}}{\gamma_\ell}$, which implies $\prod_{\ell =i+1}^j \frac{\Lambda}{\gamma_{\ell}} - \prod_{\ell =i+1}^j \frac{\lambda_{\ell-1}}{\gamma_{\ell}}\geq 0$. Thus, we have
    \begin{align*}
        a_i(\vec{\lambda})  a_j({\Lambda}) -  a_i(\Lambda)  a_j(\vec{\lambda}) =&a_i(\vec{\lambda})  a_i({\Lambda}) \prod_{\ell =i+1}^j \frac{\Lambda}{\gamma_{\ell}} -a_i(\Lambda)  a_i(\vec{\lambda}) \prod_{\ell =i+1}^j \frac{\lambda_{\ell-1}}{\gamma_{\ell}}\\
        =& a_i(\vec{\lambda})   a_i({\Lambda}) \left[\prod_{\ell =i+1}^j \frac{\Lambda}{\gamma_{\ell}} - \prod_{\ell =i+1}^j \frac{\lambda_{\ell-1}}{\gamma_{\ell}} \right]\geq 0.
    \end{align*}
Then, we have
\begin{align*}
    &\sum_{i=0}^{m-1} a_i(\vec{\lambda})\sum_{j=0}^{N} a_j({\Lambda}) - \sum_{j=0}^{m-1} a_j(\Lambda)\sum_{i=0}^{N} a_i(\vec{\lambda}) \\
    =& \sum_{i=0}^{m-1} a_i(\vec{\lambda})\left[\sum_{j=0}^{m-1} a_j({\Lambda}) + \sum_{j=m}^{N} a_j({\Lambda})\right] - \sum_{j=0}^{m-1} a_j(\Lambda)\left[\sum_{i=0}^{m-1} a_i(\vec{\lambda}) + \sum_{i=m}^{N} a_i(\vec{\lambda})\right]\\
    =& \sum_{i=0}^{m-1} a_i(\vec{\lambda}) \sum_{j=m}^{N} a_j({\Lambda}) - \sum_{j=0}^{m-1} a_j(\Lambda) \sum_{i=m}^{N} a_i(\vec{\lambda}) = \sum_{i=0}^{m-1} \sum_{j=m}^N \left[a_i(\vec{\lambda})  a_j({\Lambda}) -  a_i(\Lambda)  a_j(\vec{\lambda}) \right] \geq 0.
\end{align*}
Finally, we deduce
    \begin{align*}
        \sum_{n=0}^{m-1} P_n(\vec{\lambda}) - \sum_{n=0}^{m-1} P_n(\pi_N^\Lambda) =&  \frac{\sum_{i=0}^{m-1} a_i(\vec{\lambda})}{\sum_{i=0}^{N} a_i(\vec{\lambda})} - \frac{\sum_{j=0}^{m-1} a_j(\Lambda)}{\sum_{j=0}^{N} a_j({\Lambda})} \\
        =& \frac{\sum_{i=0}^{m-1} a_i(\vec{\lambda})\sum_{j=0}^{N} a_j({\Lambda}) - \sum_{j=0}^{m-1} a_j(\Lambda)\sum_{i=0}^{N} a_i(\vec{\lambda})}{\sum_{i=0}^{N} a_i(\vec{\lambda})\sum_{j=0}^{N} a_j({\Lambda})} \geq 0,
    \end{align*}
which proves the first inequality in the statement.  Since $\sum_{n=0}^{N} a_n({\Lambda})\leq \sum_{n=0}^\infty a_n(\Lambda)$, then  we have
\[
\sum_{n=0}^{m-1} P_n(\vec{\lambda})\geq\sum_{n=0}^{m-1} P_n(\pi_N^\Lambda) = \frac{\sum_{n=0}^{m-1} a_n(\Lambda)}{\sum_{n=0}^{N} a_n({\Lambda})} \geq  \frac{\sum_{n=0}^{m-1} a_n(\Lambda)}{\sum_{n=0}^{\infty} a_n({\Lambda})} = P_s^\Lambda,
\]
which completes the proof of the statement. \Halmos
\end{proof}

Now we are ready to prove Theorem \ref{Theorem.tp.perform}.
\begin{proof}{Proof of Theorem \ref{Theorem.tp.perform}}
When $\mathcal{C}_s\leq \mathcal{C}_q$, we have $\CC_s\leq \frac{\gamma_m}{\gamma_{n+1}}\CC_s + \frac{\gamma_{n+1}-\gamma_m}{\gamma_{n+1}}\CC_q\leq \CC_q$, which implies 
\[
\lambda\bar{F}^{-1}(\lambda/\Lambda)-\lambda\mathcal{C}_s\geq \mathcal{R}_n(\lambda)\geq \lambda\bar{F}^{-1}(\lambda/\Lambda)-\lambda\mathcal{C}_q \text{ for all $\lambda$ and $n<N$.} \addrefnum \label{cc_s.leq.cc_q.cr_n.realtion}
\]
Then, from (\ref{gain.reformulation}), (\ref{cc_s.leq.cc_q.cr_n.realtion}), and $\lambda_N^*=0$, we have 
\[
g^* = \sum_{n=0}^{N} \CR_n(\lambda_n^*)P_n(\vec{\lambda}^*)\leq \sum_{n=0}^{N} [\lambda_n^*\bar{F}^{-1}(\lambda_n^*/\Lambda)-\lambda_n^*\mathcal{C}_s] P_n(\vec{\lambda^*}).
\]
Recall that $\hat{\delta} = \sum_{n=0}^N \lambda_n^* P_n(\Vec{\lambda}^*)$ from (\ref{effective.arrival.rate}) and $\lambda\bar{F}^{-1}(\lambda/\Lambda)$ is concave. Then, using Jensen's inequality, we have
\[
g^*\leq \sum_{n=0}^{N} [\lambda_n^*\bar{F}^{-1}(\lambda_n^*/\Lambda)-\lambda_n^*\mathcal{C}_s] P_n(\vec{\lambda^*})\leq \hat{\delta} \bar{F}^{-1}(\hat{\delta}/\Lambda)-\hat{\delta}\mathcal{C}_s \underset{}{\leq} \tilde{\CR_s}, \addrefnum \label{cc_s.leq.cc_q.upper}
\]
where the last inequality follows from the definition of $\tilde{\CR}_s$.
We next derive a lower bound for the best two-price policy $g_T^*$: 
\begin{align*}
  g_T^*\geq & g_T^{\tilde{\delta_s},\tilde{\delta_q}}(N) = \sum_{n=0}^{m-1} \CR_n(\tilde{\delta_s}) P_n(\pi_N^{\tilde{\delta_s}, \tilde{\delta_q}}) + \sum_{n=m}^{N-1} \CR_n(\tilde{\delta}_q) P_n(\pi_N^{\tilde{\delta_s}, \tilde{\delta_q}})  \\
  \underset{}{\geq} & \sum_{n=0}^{m-1} \CR_n(\tilde{\delta_s}) P_n(\pi_N^{\tilde{\delta_s}, \tilde{\delta_q}}) + \sum_{n=m}^{N-1} [\tilde{\delta}_q \bar{F}^{-1}(\tilde{\delta}_q/\Lambda) - \tilde{\delta_q}\CC_q] P_n(\pi_N^{\tilde{\delta_s}, \tilde{\delta_q}}) = \tilde{\CR_s} P_s^T + \tilde{\CR}_q P_q^T, \addrefnum \label{cc_s.leq.cc_q.lower}
\end{align*}
where the second inequality holds due to (\ref{cc_s.leq.cc_q.cr_n.realtion}). Thus, from (\ref{cc_s.leq.cc_q.upper}) and (\ref{cc_s.leq.cc_q.lower}), we have 
\[
g^*-g_T^* \leq \tilde{\CR_s} - [\tilde{\CR_s} P_s^T + \tilde{\CR}_q P_q^T] = (\tilde{\CR}_s - \tilde{\CR}_q) P_q^T + \tilde{\CR}_s P_N^T.
\]

If $\mathcal{C}_s\geq \mathcal{C}_q$, we have $\CC_s\geq \frac{\gamma_m}{\gamma_{n+1}}\CC_s + \frac{\gamma_{n+1}-\gamma_m}{\gamma_{n+1}}\CC_q\geq \CC_q$, which implies 
\[
\lambda\bar{F}^{-1}(\lambda/\Lambda)-\lambda\mathcal{C}_s\leq \mathcal{R}_n(\lambda)\leq \lambda\bar{F}^{-1}(\lambda/\Lambda)-\lambda\mathcal{C}_q \text{ for all $\lambda$ and $n<N$}. \addrefnum \label{cc_s.geq.cc_q.cr_n.realtion}
\]
We define the short-hand notation $\hat{\delta_s}$ and $\hat{\delta_q}$ as follows:
\[
\hat{\delta_s} = \sum_{n=0}^{m-1} \lambda_n^* \dfrac{P_n(\vec{\lambda}^*)}{\sum_{n=0}^{m-1} P_n(\vec{\lambda}^*)} \text{ and } \hat{\delta_q} = \sum_{n=m}^{N} \lambda_n^* \dfrac{P_n(\vec{\lambda}^*)}{\sum_{n=m}^{N} P_n(\vec{\lambda}^*)}.
\]
Then, we derive the upper bound of optimal gain $g^*$ using $P_s^\Lambda,\,\tilde{\CR}_s$, and $\tilde{\CR}_q$:
\begin{align*}
    g^* =& \sum_{n=0}^{N} \CR_n(\lambda_n^*) P_n(\vec{\lambda}^*) \underset{}{\leq} \sum_{n=0}^{m-1} [\lambda_n^*\bar{F}^{-1}(\lambda_n^*/\Lambda)-\lambda_n^*\CC_s] P_n(\vec{\lambda}^*) + \sum_{n=m}^{N} [\lambda_n^*\bar{F}^{-1}(\lambda_n^*/\Lambda)-\lambda_n^*\CC_q] P_n(\vec{\lambda}^*)\\
    =& \left( \sum_{n=0}^{m-1} P_n(\vec{\lambda}^*) \right) \left( \sum_{n=0}^{m-1} [\lambda_n^*\bar{F}^{-1}(\lambda_n^*/\Lambda)-\lambda_n^*\CC_s] \dfrac{P_n(\vec{\lambda}^*)}{\sum_{n=0}^{m-1} P_n(\vec{\lambda}^*)}  \right) \\
    &+ \left( \sum_{n=m}^{N} P_n(\vec{\lambda}^*) \right) \left( \sum_{n=m}^{N} [\lambda_n^*\bar{F}^{-1}(\lambda_n^*/\Lambda)-\lambda_n^*\CC_q] \dfrac{P_n(\vec{\lambda}^*)}{\sum_{n=m}^{N} P_n(\vec{\lambda}^*)}  \right)\\
    \underset{}{\leq} & \left( \sum_{n=0}^{m-1} P_n(\vec{\lambda}^*) \right) [\hat{\delta_s}\bar{F}^{-1}(\hat{\delta_s}/\Lambda)-\hat{\delta_s}\CC_s] + \left( \sum_{n=m}^{N} P_n(\vec{\lambda}^*) \right) [\hat{\delta_q}\bar{F}^{-1}(\hat{\delta_q}/\Lambda)-\hat{\delta_q}\CC_q]\\
    \underset{}{\leq} & \left( \sum_{n=0}^{m-1} P_n(\vec{\lambda}^*) \right) \tilde{\CR_s} + \left( \sum_{n=m}^{N} P_n(\vec{\lambda}^*) \right) \tilde{\CR}_q. \addrefnum \label{cc_s.geq.cc_q.upper}
\end{align*}
The first inequality holds from (\ref{cc_s.geq.cc_q.cr_n.realtion}) and the second inequality holds from the definitions of $\hat{\delta_s}$ and $\hat{\delta_q}$ and applying Jensen's inequality on each term. Using the definition of $\tilde{\CR}_s$ and $\tilde{\CR}_q$, we obtain the third inequality. Finally, notice that $\CC_s\geq \CC_q$ implies $\tilde{\CR}_s\leq \tilde{\CR}_q$ and from Lemma \ref{lemma.idle.probability.bound}, we have 
\[
\sum_{n=0}^{m-1} P_n(\vec{\lambda}^*)\geq P_s^\Lambda, \text{ which implies } P_q^\Lambda = 1- P_s^\Lambda \geq 1 - \sum_{n=0}^{m-1} P_n(\vec{\lambda}^*)= \sum_{n=m}^{N} P_n(\vec{\lambda}^*).
\]
Therefore, (\ref{cc_s.geq.cc_q.upper}) implies that
\begin{align*}
P_s^\Lambda  \tilde{\CR_s} + P_q^\Lambda \tilde{\CR}_q-g^*
\geq & P_s^\Lambda  \tilde{\CR_s} + P_q^\Lambda \tilde{\CR}_q - \left[ \left( \sum_{n=0}^{m-1} P_n(\vec{\lambda}^*) \right) \tilde{\CR_s} + \left( \sum_{n=m}^{N} P_n(\vec{\lambda}^*) \right) \tilde{\CR}_q \right]\\
=&(1-P_q^\Lambda)  \tilde{\CR_s} + P_q^\Lambda \tilde{\CR}_q - \left[ \left( 1- \sum_{n=m}^{N} P_n(\vec{\lambda}^*) \right) \tilde{\CR_s} + \left( \sum_{n=m}^{N} P_n(\vec{\lambda}^*) \right) \tilde{\CR}_q \right]\\
=& \left[ P_q^\Lambda - \sum_{n=m}^{N} P_n(\vec{\lambda}^*) \right]\left( \tilde{\CR}_q - \tilde{\CR}_s  \right) \geq 0. \addrefnum \label{cc_s.geq.cc_q.upper.geq.0}
\end{align*}
We next derive a lower bound for $g_T^*$:
\begin{align*}
    g_T^* \geq & g_T^{\tilde{\delta_s},\tilde{\delta_q}}(N) =\sum_{n=0}^{m-1} \CR_n(\tilde{\delta_s}) P_n(\pi_N^{\tilde{\delta_s}, \tilde{\delta_q}}) + \sum_{n=m}^{N-1} \CR_n(\tilde{\delta}_q) P_n(\pi_N^{\tilde{\delta_s}, \tilde{\delta_q}})  \\
  \underset{}{\geq} & \sum_{n=0}^{m-1} \CR_n(\tilde{\delta_s}) P_n(\pi_N^{\tilde{\delta_s}, \tilde{\delta_q}}) + \sum_{n=m}^{N-1} [\tilde{\delta}_q \bar{F}^{-1}(\tilde{\delta}_q/\Lambda) - \tilde{\delta_q}\CC_s] P_n(\pi_N^{\tilde{\delta_s}, \tilde{\delta_q}}) \\
  =& \tilde{\CR_s} P_s^T + [\tilde{\delta}_q \bar{F}^{-1}(\tilde{\delta}_q/\Lambda) - \tilde{\delta_q}\CC_s] P_q^T =\tilde{\CR_s} P_s^T + \tilde{\CR}_q P_q^T - \tilde{\delta_q} (\CC_s -\CC_q) P_q^T, \addrefnum \label{cc_s.geq.cc_q.lower}
\end{align*}
where the second inequality follows from (\ref{cc_s.geq.cc_q.cr_n.realtion}). Combining the bounds in (\ref{cc_s.geq.cc_q.upper.geq.0}) and  (\ref{cc_s.geq.cc_q.lower}), we have
\begin{align*}
    g^*-g_T^* \leq & P_s^\Lambda  \tilde{\CR_s} + P_q^\Lambda \tilde{\CR}_q - [\tilde{\CR_s} P_s^T + \tilde{\CR}_q P_q^T - \tilde{\delta_q} (\CC_s -\CC_q) P_q^T]\\
    =& \tilde{\CR_s}(P_s^\Lambda- P_s^T) + \tilde{\CR_q}(P_q^\Lambda -P_q^T) + \tilde{\delta_q} (\mathcal{C}_s-\mathcal{C}_q) P_q^T\\
    =& \tilde{\CR_s} [1-P_q^\Lambda- (1 - P_q^T -P_N^T )] + \tilde{\CR_q}(P_q^\Lambda -P_q^T) + \tilde{\delta_q} (\mathcal{C}_s-\mathcal{C}_q) P_q^T \\
    =& (\tilde{\CR_q} - \tilde{\CR_s}) (P_q^\Lambda - P_q^T) + \tilde{\delta_q} (\mathcal{C}_s-\mathcal{C}_q) P_q^T + \tilde{\CR}_s P_N^T,
\end{align*}
which completes the proof. \Halmos
\end{proof}

\begin{proof}{Proof of Proposition \ref{Prop.tp.optimal}}
When $N\to\infty$, we have $P_N^T\to 0$ since the system is stable under any policy. Since $\mathcal{C}_s=\mathcal{C}_q$ implies $\tilde{\CR}_s=\tilde{\CR}_q$,  Theorem \ref{Theorem.tp.perform} yields
\[
0\leq g^*-g_T^* \leq \tilde{\CR}_s P_N^T \to 0 \text{ as } N\to \infty. \Halmos
\]
\end{proof}
}

\section{Additional Numerical Results for Section \ref{subsec:necessary_condition}}

\begin{figure}[h]
	\centering
	\caption{Scatter Plots for the Conditions of Monotonicity}
	\begin{minipage}{0.45\textwidth}
		\centering
		\includegraphics[width=\columnwidth]{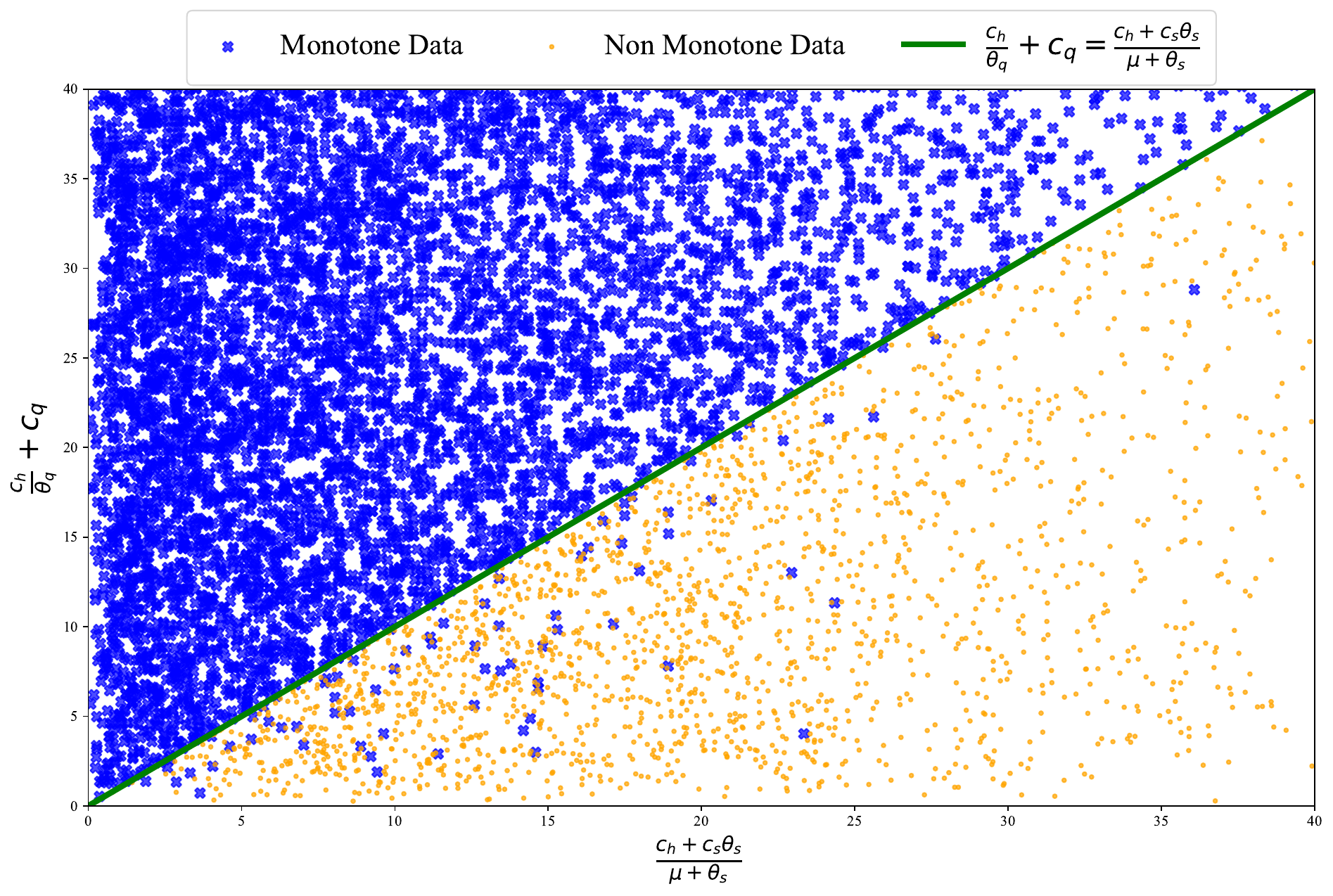}
		\subcaption{Evaluation distribution as $U(20,50)$}
        \label{Monotone_condition_uni_multi}
	\end{minipage}
	\begin{minipage}{0.45\textwidth}
		\centering
		\includegraphics[width=\columnwidth]{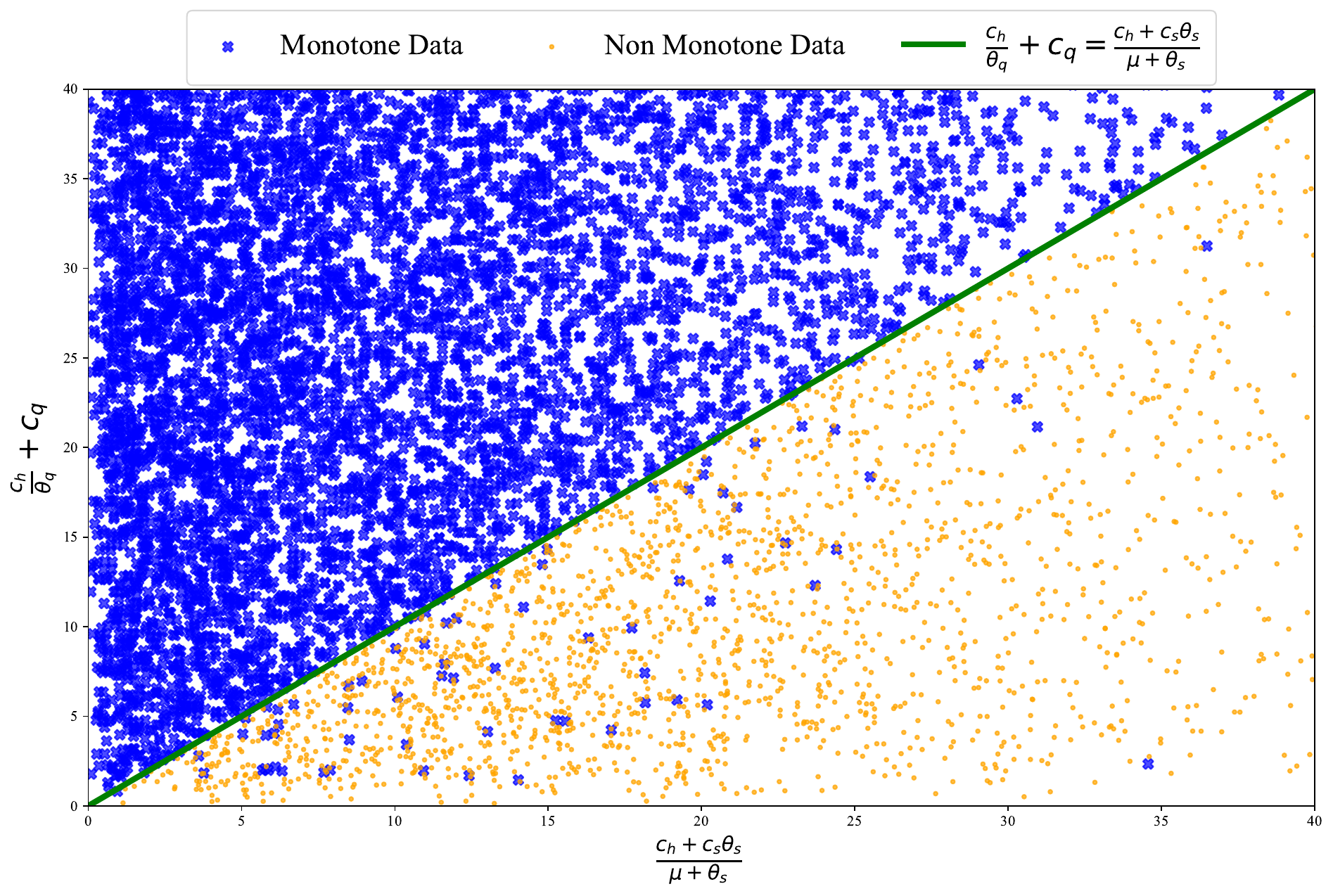}
		\subcaption{Evaluation distribution as $Exp(1/35)$}
        \label{Monotone_condition_exp_multi}
	\end{minipage}
\end{figure}

We also consider the case when the system has multiple servers. As for Figures \ref{Monotone_condition_unif} and \ref{Monotone_condition_exp}, we randomly select $\Lambda$, $\mu$, $\theta_s$, $\theta_q$, $c_s$, $c_q$, and $c_h$ independently from a uniform distribution with parameters 0 and 50. We choose the evaluation distribution to be either a uniform distribution with range [20,50] or an exponential distribution with mean 35. We let the number of server $m=4$ and the buffer size $N$ be a random integer selected from the set $\{5,6,\ldots,20\}$ with equal probability. We randomly generate 10,000 data points and the results are shown in Figures \ref{Monotone_condition_uni_multi} and \ref{Monotone_condition_exp_multi}.

\vspace{-0.4in}

Figure \ref{Monotone_condition_uni_multi} reveals that in uniform evaluation distributions with multiple servers ($m=4$), we find data points where optimal arrival rate policies decrease monotonically despite violating the condition $\frac{c_h+ c_s\theta_s}{\mu+\theta_s} \leq \frac{c_h}{\theta_q} + c_q$. Compared to Figure \ref{Monotone_condition_unif}, Figure \ref{Monotone_condition_uni_multi} exhibits this pattern more frequently (with more monotone data below the line $\frac{c_h+ c_s\theta_s}{\mu+\theta_s} = \frac{c_h}{\theta_q} + c_q$), further confirming that the condition $\frac{c_h+ c_s\theta_s}{\mu+\theta_s} \leq \frac{c_h}{\theta_q} + c_q$ is not necessary for monotone decreasing optimal policies.

Figure \ref{Monotone_condition_exp_multi} demonstrates that with exponential evaluation distribution and multiple servers, several data points exhibit monotone decreasing structure with $\frac{c_h+ c_s\theta_s}{\mu+\theta_s} > \frac{c_h}{\theta_q} + c_q$. Compared to Figure \ref{Monotone_condition_exp}, more data points fall below the line $\frac{c_h+ c_s\theta_s}{\mu+\theta_s} = \frac{c_h}{\theta_q} + c_q$, which further confirms that $\frac{c_h+ c_s\theta_s}{\mu+\theta_s} \leq \frac{c_h}{\theta_q} + c_q$ is not a necessary condition for the optimal arrival rate policy to be monotone and that this behavior is more pronounced for $m=4$ than that for $m=1$.

\section{Additional Numerical Results for Section \ref{sec.numerical.example}}

{We use the same procedure that generates Table \ref{Histogram_1_server} to randomly generate system parameters when the evaluation distribution is exponentially distributed with mean 35.}
{Note that since the exponential distribution has an infinite support, it does not fit the assumptions of our model. However, our numerical results demonstrate that our algorithms and heuristics remain robust and perform well even in this setting. In practice, distributions with unbounded support often can be effectively approximated by truncating the tail.}
Table \ref{Histogram_exp} shows the frequency of the performance ratios $R_S,\,R_C,$ and $R_T$ when the evaluation distribution is exponential.

{\begingroup\renewcommand{\baselinestretch}{0.6}\normalsize

\begin{table}[h]  
	\centering
	\caption{Performance ratios of heuristics with evaluation distribution as $Exp(35)$}
	\begin{minipage}{0.5\textwidth} 
		\centering
		\renewcommand{\arraystretch}{1.5}
		\begin{tabular}{|>{}c|>{}c|>{}c|>{}c|}  
			\hline
			{\tiny Performance ratio} & {\tiny Best static policy} & {\tiny Best cutoff-static policy} & {\tiny Best two-price policy}\\ \hline
			{\small $0\sim0.1$} & 0 & 0 & 0\\ \hline
			{\small $0.1\sim0.2$} & 0 & 0 & 0 \\ \hline
            {\small $0.2\sim0.3$} & 0 & 0 & 0 \\ \hline
            {\small $0.3\sim0.4$} & 0 & 0 & 0\\ \hline
            {\small $0.4\sim0.5$} & 0 & 0 & 0 \\ \hline
            {\small $0.5\sim0.6$} & 0 & 0 & 0 \\ \hline
            {\small $0.6\sim0.7$} & 0 & 0 & 0 \\ \hline
            {\small $0.7\sim0.8$} & 4 & 1 & 0 \\ \hline
            {\small $0.8\sim0.9$} & 46 & 6 & 0\\ \hline
            {\small $0.9\sim1$} & 9,950 & 9,993 & 10,000 \\ \hline
		\end{tabular}
        \vspace{0.1in}
		\subcaption{Full scale}
		\label{Histogram_exp.png}
	\end{minipage}%
	\hfill
	\begin{minipage}{0.5\textwidth} 
		\centering
		\renewcommand{\arraystretch}{1.5}
		\begin{tabular}{|>{}c|>{}c|>{}c|>{}c|}  
			\hline
			{\tiny Performance ratio} & {\tiny Best static policy} & {\tiny Best cutoff-static policy} & {\tiny Best two-price policy}\\ \hline
			{\small $0.90\sim0.91$} & 19 & 3 & 0\\ \hline
			{\small $0.91\sim0.92$} & 20 & 1 & 0 \\ \hline
            {\small $0.92\sim0.93$} & 37 & 23 & 0 \\ \hline
            {\small $0.93\sim0.94$} &52 & 44 & 0\\ \hline
            {\small $0.94\sim0.95$} & 61 & 68 & 0 \\ \hline
            {\small $0.95\sim0.96$} & 122 & 123 & 0 \\ \hline
            {\small $0.96\sim0.97$} & 233 & 259 & 0 \\ \hline
            {\small $0.97\sim0.98$} & 382 & 397 & 0 \\ \hline
            {\small $0.98\sim0.99$} & 861 & 882 & 1\\ \hline
            {\small $0.99\sim1$} & 8,163 & 8,193 & 9,999 \\ \hline
		\end{tabular}
        \vspace{0.1in}
		\subcaption{Zoom in}
		\label{Histogram_exp_90.png}
	\end{minipage}

	\label{Histogram_exp}
\end{table}
\endgroup}

From Table \ref{Histogram_exp.png}, we notice that the best static policy, best cutoff-static policy, and best two-price policy are performing better compared to the case when the evaluation distribution is uniform. Namely, all three heuristics achieve a performance ratio above 0.7 for all 10,000 data points. 
The best cutoff-static policy improves upon the best static policy but there is one case with a performance ratio under 0.8 (as compared to four cases for the best static policy). This shows that best cutoff-static policy is near optimal but not robust in some cases. The best two-price policy still has the best performance among all three heuristics, with performance ratios above 0.9 for all 10,000 data points.

In Table \ref{Histogram_exp_90.png}, we zoom into the last bar of Table \ref{Histogram_exp.png} to focus on the data where the performance ratio is over 0.9 to further demonstrate the advantage of the best two-price policy.
Table \ref{Histogram_exp_90.png}  zooms into the data with a performance ratio over 0.9. From Table \ref{Histogram_exp_90.png}, we observe that the best cutoff-static policy and the best two-price policy both perform better compared to the case when the evaluation distribution is uniform (see Table \ref{Histogram_1_server_90.png}). In particular, the best two-price policy has guaranteed 0.98 performance ratios for all 10,000 data points and has performance ratios above 0.99 in all but one of these cases. 
Therefore, the numerical experiments suggest the near optimality of the best cutoff-static and two-price policies and emphasize the robustness of the latter one. In practice, these two heuristics achieve near optimal performances and are easier to implement compared to the optimal dynamic policy.

\end{document}